\documentclass[10pt]{article}
\usepackage{geometry}
\usepackage{amssymb,amsmath,amsfonts,amsthm,mathtools,enumerate,color}
\usepackage{mathrsfs}

\usepackage[inline,shortlabels]{enumitem}

\usepackage[title,titletoc,header]{appendix}
\usepackage{graphicx}

\usepackage{fancyhdr}
\usepackage[titletoc]{appendix}
\usepackage{tikz}
\usetikzlibrary{arrows,positioning,shapes.geometric}
\usepackage{titlepic}
\usepackage{algorithm,algorithmic}

\usepackage{booktabs}

\usepackage{cancel}

\usepackage{float}
\usepackage{subcaption}
\DeclareMathOperator*{\argmin}{arg\,min}
\usepackage{indentfirst}
\usepackage{multicol}
\usepackage{url}

\graphicspath{ {./figures/} }
\usepackage[outdir=./figures/]{epstopdf} 

\usepackage{hyperref}
\hypersetup{
    colorlinks=true, 
    linktoc=all,     
    linkcolor=blue,  
}
\allowdisplaybreaks

\theoremstyle{plain}

\newtheorem*{lemma*}{Lemma}

\newtheorem*{theorem*}{Theorem}

\theoremstyle{definition}

\newtheorem{remark}{Remark}[section]

\newtheorem*{remark*}{Remark}
\newtheorem*{notation*}{Notation}

\numberwithin{equation}{section}

\def\cD{\mathcal{D}}
\def\cE{\mathcal{E}}
\def\cF{\mathcal{F}}

\def\cH{\mathcal{H}}

\def\cL{\mathcal{L}}

\def\cN{\mathcal{N}}
\def\cO{\mathcal{O}}
\def\cP{\mathcal{P}}

\def\cS{\mathcal{S}}

\def\cU{\mathcal{U}}

\def\sR{{\mathbb R}}
\def\Z{{\mathbb Z}}

\newcommand{\EE}{{\mathbb E}}

\newcommand{\PP}{{\mathbb P}}

\newcommand{\tr}{\textnormal{tr}}

\DeclareMathOperator*{\prox}{prox}

\usepackage{array}
\newcolumntype{P}[1]{>{\centering\arraybackslash}p{#1}}

\def\trans{\mathsf{T}}

\def\sE{{\mathbb{E}}}
\def\cD{\mathcal{D}}
\def\R{\mathbb R}    
\def\N{\mathbb N}  
\def\ee{\varepsilon}
   
\def\la{\langle} \def\ra{\rangle}
  
    \def\E{\mathbb E}
\def\d{\, \text{\rm{d}}} 
 \def\1{\lesssim}
 \def\a{\alpha}
\def\ms{\medskip}
 \def\p{\partial}
\def\t{\times}

\def\bb{\begin{equation}} \def\ee{\end{equation}}
\def\bbn{\begin{equation*}} \def\een{\end{equation*}}

\begin{document}

\title{
A fast iterative PDE-based algorithm for feedback controls of  nonsmooth mean-field control problems
}

\author{
Christoph Reisinger\thanks{
Mathematical Institute, University of Oxford, Oxford OX2 6GG, UK
 ({\tt christoph.reisinger@maths.ox.ac.uk, 
wolfgang.stockinger@maths.ox.ac.uk,
yufei.zhang@maths.ox.ac.uk})}
\and
Wolfgang Stockinger\footnotemark[1]
\and
Yufei Zhang\footnotemark[1]
}
\date{}

\maketitle

\noindent\textbf{Abstract.} 
We propose a PDE-based  accelerated gradient algorithm 
for optimal feedback controls of
McKean--Vlasov dynamics that involve mean-field interactions both in the state and action.
The method exploits a forward-backward splitting approach and
iteratively refines 
the approximate controls 
   based on  the  gradients of  smooth costs, 
   the proximal maps of  nonsmooth costs, and dynamically updated momentum parameters.
At each step, 
 the state dynamics is approximated via a particle system,
 and  the required gradient is evaluated through 
 a  coupled system of nonlocal linear PDEs. 
The latter is solved by finite difference approximation or neural network-based residual approximation,
depending on the state dimension. 
We present exhaustive numerical experiments for 
low and high-dimensional
mean-field control problems, including  sparse stabilization of  stochastic Cucker--Smale models,
which reveal that 
our algorithm captures important structures of the optimal feedback control
and achieves a robust performance with respect to parameter perturbation.

\medskip

\noindent
\textbf{Key words.} 
Controlled McKean--Vlasov diffusion, 
optimal gradient method, 
monotone scheme,
neural network,
sparse control,
stochastic Cucker--Smale model.

\ms
\noindent
\textbf{AMS subject classifications.} 
49N80,   60H35, 35Q93, 93A16

\medskip


\section{Introduction}

In this article, we propose a class of iterative  methods for solving
 mean-field  control (MFC) problems, 
where the state dynamics and 
cost functions depend upon the joint law of the state and the control processes.
Let $T>0$ be a given terminal time, $W=(W_t)_{t\in [0,T]}$ 
be an $n$-dimensional Brownian motion defined on the probability space 
$(\Omega,\cF,\PP)$, $\mathbb{F} = (\cF_t)_{t\in [0,T]}$ 
the natural filtration of $W$ augmented with an independent $\sigma$-algebra $\cF_0$,
and 
$\cH^2(\sR^k)$ be the set of admissible controls containing all 
 $\sR^k$-valued
square integrable $\mathbb{F}$-progressively measurable processes.
For a given $\cF_0$-measurable 
initial state $\xi\in L^2(\Omega;\R^d)$
and a control $\alpha\in \cH^2(\sR^k)$,
we consider the state process 
governed by the following controlled McKean--Vlasov diffusion:
\begin{align}\label{cforward}
\mathrm{d} X_t=
b(t,X_t,\alpha_t,\mathcal{L}_{(X_t,\alpha_t)})\, \mathrm{d} t
+\sigma(t,X_t,\alpha_t,\mathcal{L}_{(X_t,\alpha_t)})\, \mathrm{d} W_t, 
\; t\in [0,T]; 
\quad X_0=\xi,
\end{align}
where 
$b:[0,T]\t \sR^d\t \sR^k \t \cP_2(\sR^d\t \sR^k)\to \sR^d$ and
$\sigma:[0,T]\t \sR^d\t \sR^k \t \cP_2(\sR^d\t \sR^k)\to \sR^{d\t n}$
are  sufficiently regular    functions
such that \eqref{cforward}
admits a unique square integrable  solution $X^\a$.
The value function of the optimal control problem is defined by
\begin{equation}\label{cvalue}
J^\star(\xi)=\inf_{\alpha \in \mathcal{H}^{2}(\R^k)} J(\alpha;\xi),
\quad
\textnormal{with}
\;
J(\alpha;\xi)=\E \bigg[
\int_0^T \Big( f(t,X^{\alpha}_t,\alpha_t,\mathcal{L}_{(X^{\alpha}_t,\alpha_t)}) + \ell(\alpha_t) \Big) \, \mathrm{d}t+g(X^{\alpha}_T,\mathcal{L}_{X^{\alpha}_T})
\bigg],
\end{equation}
where  $f:[0,T]\t \sR^d\t \sR^k \t \cP_2(\sR^d\t \sR^k)\to \sR$, $g: \sR^d\t  \cP_2(\sR^d)\to \sR$ are 
 differentiable functions 
of at most quadratic growth,
and 
 $\ell: \mathbb{R}^{k} \to \mathbb{R}\cup\{\infty\}$ is a 
 proper, 
 lower semicontinuous and convex function.
 Above and hereafter, $\mathcal{L}_U$ denotes  the law of a  random variable $U$,
and $\cP_2(E)$
denotes 
the Wasserstein space of probability measures on  
the Euclidean space 
$E$ with finite second moment.

The above MFC  problem 
extends classical stochastic control problems
by  allowing the mean-field interactions  through the joint distribution of the state and control processes.
It describes large population equilibria of interacting individuals  controlled by a central planner,
and plays an important role in 
economics \cite{basei2017linear,acciaio2019extended,achdou2020mean},
production management \cite{gobet2019extended},
biology 
\cite{caponigro2013sparse, bailo2018optimal, albi2021moment}
and social interactions \cite{achdou2015system, achdou2016mean,   achdou2020mean}.
Moreover,  
the
extended real-valued 
  function  $\ell$ in \eqref{cvalue}
  includes
  important  examples such as
 characteristic functions of convex sets 
representing 
  control constraints \cite{acciaio2019extended, reisinger2020regularity}, 
  $\ell_1$-norm based regularizations   
used to induce sparsity or switching properties of the optimal control \cite{caponigro2013sparse, bailo2018optimal},
  and 
entropy regularizations  in  machine learning \cite{vsivska2020gradient,guo2021reinforcement}.

To solve \eqref{cforward}-\eqref{cvalue},
we aim to construct  an optimal (decentralized) feedback control, i.e.,
a sufficiently regular function $\phi^\star:[0,T]\t \sR^d\to \sR^k$ 
 such that 
 the corresponding controlled dynamics \eqref{cforward} (with $\alpha_t=\phi^\star(t,X_t)$)
admits a unique solution $X^{\phi^\star}$ and 
 $J^\star(\xi)=J(\phi^\star(\cdot, X^{\phi^\star}_\cdot),\xi)$.\footnotemark
 \footnotetext{
 $\phi^\star$ is called a decentralized control 
as it  acts explicitly only on the time and state variables,
while the dependence on the (deterministic) marginal laws of the optimal state and control processes
is implicit through the time dependence.}
 The existence of such feedback controls has been shown in 
\cite{reisinger2020regularity} (see also
\cite[Section 6.4]{carmona2018probabilistic} for special cases without  mean-field interactions through  control variables).
The main advantage of 
a feedback strategy is that 
implementing the optimal control
reduces to simple function evaluation at the  current state of the system.
Moreover, 
a feedback control allows us to interpret  the  mechanism of the optimal control.
This is particularly   important for
economics and social science \cite{miller2019explanation}, where one would like to 
understand the cause of a  decision, and explain its dependence on  the state dynamics and the objective function. 
A feedback control also allows us to analyze the failure of certain control policies  for fault diagnosis.

\paragraph{Existing  numerical methods for MFC problems and their limitations.}
As analytical solutions to 
optimal  feedback controls of \eqref{cforward}-\eqref{cvalue}
are rarely available, 
numerical schemes for solving such control problems become vital. 
Due to the nonlinear dependence on the marginal laws, 
it is difficult to
follow 
the classical dynamic programming (DP)  approach
for  constructing  optimal feedback controls.
The main challenge here 
is that 
to construct
 feedback controls,
DP requires us to find
derivatives of solutions to 
an  infinite-dimensional 
   Hamilton-Jacobi-Bellman (HJB) partial differential equation (PDE) 
  defined on
   $[0,T]\t \cP_2(\sR^d)$ (see e.g.,~\cite{pham2018bellman}),
  which is  computationally intractable.  
Hence, most existing numerical methods are  
based on the optimization interpretation of \eqref{cforward}-\eqref{cvalue}.

The most straightforward approach for 
solving  \eqref{cforward}-\eqref{cvalue} 
is to 
restrict the optimization  
over  feedback controls within  a prescribed parametric family,
i.e., the so-called   policy gradient (PG) method
(see e.g.,  \cite{carmona2019convergence}).
It
approximates the  optimal feedback control
in a parametric form  depending on weights $\theta$
(for instance a deep neural network),
and then 
 seeks the optimal approximation by performing gradient descent of $J$ with respect to the weights $\theta$
  based on simulated trajectories  
of the state process.
By exploiting an efficient neural network representation of the feedback control,
the PG method can be adapted to solve MFC problems with high-dimensional state processes.

However, 
the PG method   has two serious drawbacks,
especially for solving control problems with nonlinear dynamics and  nonsmooth costs
in \eqref{cforward}-\eqref{cvalue}.
Firstly, 
as the loss functional $J$ is  nonconvex and nonsmooth in the weights $\theta$ of the numerical feedback controls,
 there is no theoretical guarantee on the convergence of the PG method for 
 solving nonsmooth MFC problems.
  In practice, even though the PG method 
  may minimize the loss function reasonably well, 
  the resulting approximate feedback control
  often fails to capture important structural features
  of the optimal feedback control, 
  and consequently 
   lacks the capability to provide sufficient
 insights of the optimal decision;
  see  Figures \ref{fig:cs_2d_beta10}, \ref{fig:cs_6d_beta0_fb} and \ref{fig:value_conv_6d_l1}  
  in Section \ref{sec:CS_6d}
  where
 the PG method 
  ignores the temporal and spatial nonlinearity and  the sparsity of the optimal controls.
  Secondly, as the PG method computes 
approximate  feedback controls 
 based  purely on   sample   trajectories of the state dynamics, 
it cannot recover  optimal feedback controls
 outside the support of the optimal state process.
 Consequently, the performance of the approximate  feedback control
 in general 
 can be 
 very sensitive to perturbations of    the (random) initial state $\xi$  of \eqref{cforward};
 see Section \ref{SEC:PF} for details. 
This is undesirable for practical applications of MFC problems,
as the initial condition $\xi$ describes 
the  asymptotic regime of initial conditions of a large
number of players, and often cannot be observed exactly.
 
 Another approach is to solve  the optimality systems 
 arising from applying the Pontryagin Maximum Principle (PMP) to   \eqref{cforward}-\eqref{cvalue}.
Existing works consider special cases with neither nonsmooth costs nor mean-field interactions through control variables,
and design   numerical methods
based on either a probabilistic or deterministic formulation
(see e.g., \cite[Section 6.2.4]{carmona2018probabilistic}). 
The probabilistic formulation represents 
 an optimal control $\a^\star\in \cH^2(\sR^k)$ as 
$\a^\star_t=\hat{\a}(t, X^{\a^\star}_t, \mathcal{L}_{X^{\alpha^\star}_t},Y^{\a^\star}_t,Z^{\a^\star})$
 for $\d t\otimes \d \PP$-a.e.,
with $\hat{\a}$ being the pointwise minimizer of an associated Hamiltonian $\cH$,
and 
$(X^{\a^\star}, Y^{\a^\star},Z^{\a^\star})$ 
being the solution to 
 a coupled forward-backward stochastic differential equation (FBSDE) depending on $\hat{\a}$. 
The coupled FBSDE can  be solved by first
 representing the solution in terms of grid functions  \cite{angiuli2019cemracs,reisinger2020posteriori},
binomial trees \cite{angiuli2019cemracs} or neural networks  \cite{weinan2017deep, carmona2019convergence,germain2019numerical},
and then employing regression methods
to obtain the optimal approximation.
Similarly, 
the deterministic formulation represents 
 an optimal feedback control $\phi^\star$ as
$\phi^\star(t,x)=\bar{\a}(t, x, \mu^\star_t, (\nabla_x v)(t,x), (\textrm{Hess}_{x} v)(t,x))$
 for all $(t,x)\in [0,T]\t\sR^d$.
Here $\bar{\a}$ is the pointwise minimizer of an associated Hamiltonian $\mathscr{H}$ (possibly different from $\cH$),
and $(\mu^\star,v)$ satisfy a coupled  Fokker-Planck (FP)-HJB PDE system
depending on $\bar{\a}$,
which consists of a nonlinear FP equation 
for the marginal distribution $\mu^\star$ of the optimal state process,
and of a nonlinear HJB equation for the adjoint variable $v$.
The  FP-HJB system  can then be solved by finite difference methods 
as in   \cite{achdou2015system, achdou2016mean, achdou2020mean}
or by neural network methods as in \cite{carmona2021convergence_I}.

We observe, however, that the above PMP  approach suffers from the following limitations.
Firstly, 
the derivation of the optimality systems
relies heavily  on
the analytic expression of
  pointwise minimizers of the Hamiltonians, 
  which may not be available for general control problems.
More crucially, as pointed out in \cite{acciaio2019extended},
when there is a nonlinear
dependence on the law of the control,
the PMP   in general cannot be expressed in terms of a  pointwise minimization
of Hamiltonian
(see \cite{reisinger2020regularity} for a detailed investigation of this issue). 
These factors prevent us from applying the PMP  approach to solve 
\eqref{cforward}-\eqref{cvalue} 
  with general cost functions and control interactions. 
Secondly, similar to the PG method,
solving the   coupled FBSDE via regression
focuses  mainly along     trajectories of the optimal 
state, 
and consequently 
would result in an approximate  feedback control
that 
 is sensitive  to perturbation of the  initial state
(see Section \ref{SEC:PF}).
Finally,  solutions to the nonlinear FP equation 
in general only exist  
in the sense of distributions
and often admit temporal and spatial singularity. 
This creates  significant  numerical challenges,
especially when the diffusion coefficient $\sigma$ of \eqref{cforward} 
is degenerate (as in most kinetic models)
or  the initial state $\xi$ has a singular density;
see Section \ref{SEC:NUM} for concrete examples. 
In particular, in the high-dimensional setting, one may need to employ neural networks with  complex structures
to approximate such irregular solutions, 
which subsequently
  results in  challenging optimization problems for training the networks. 

\paragraph{Our contributions and related works.}
This paper proposes a class of iterative algorithms to construct  
optimal feedback controls for nonsmooth MFC problems  \eqref{cforward}-\eqref{cvalue}.

 \begin{itemize}
 
 \item
We construct a sequence of feedback controls $\phi^m: [0,T] \t \sR^d\to \sR^k$ whose realized  control processes minimize the functional $J$.
{\color{black} In the present setting, 
directly applying gradient descent   to either the  
stochastic formulation of the control problem over open-loop controls, or to  the deterministic reformulation  over feedback  controls 
has several critical deficiencies (see Sections  \ref{sec:GRADNM} and   \ref{sec:GRADM} for details).
To overcome these shortcomings, 
we propose a heuristic combination of the two formulations, 
and tailor  Nesterov's accelerated proximal gradient (NAG) method (also known as the Fast Iterative Shrinkage-Thresholding Algorithm \cite{beck2009fast})  to our problem.}
At each iteration with  given feedback  control,
we 
 evaluate  the   gradient of the smooth costs
at the corresponding realized control process, 
 and   update the feedback control   
 by incorporating the  gradient information,
 the proximal map of the nonsmooth cost $\ell$,
 and  an explicit dynamically updated momentum parameter.

The proposed accelerated proximal gradient  approach has the following advantages
in solving   \eqref{cforward}-\eqref{cvalue}:
(i) 
unlike the aforementioned PMP approach,
our algorithm  requires neither 
(pointwise) analytical minimization of the Hamiltonian 
nor deriving the optimality systems, 
and hence   can be applied to   MFC problems with general  mean-field interactions through the control variables;
(ii)  
our algorithm  shares the same  computational complexity  
as the gradient-based algorithms in \cite{pfeiffer2017numerical,archibald2020efficient,vsivska2020gradient, kerimkulov2021modified},
 but enjoys an accelerated 
 convergence rate 
 and can handle general convex
 nonsmooth costs,
including  $\ell_1$-regularizers and control constraints. 
In fact, such an accelerated gradient iteration is known to be an optimal first order (gradient) method (in the sense of complexity analysis) for minimizing finite-dimensional nonsmooth functions
\cite{beck2009fast};
 (iii)
 our method represents the control iterates in a feedback form
 (cf., \cite{archibald2020efficient,vsivska2020gradient, kerimkulov2021modified} which update controls  as stochastic processes), 
 and avoids the curse of dimensionality in the gradient evaluation
 as the number of iterations tends to infinity (see Section \ref{SEC:NAG} for details);
 (iv) {\color{black} 
  compared to  directly applying   NAG method to optimize \eqref{cvalue} over feedback controls, our algorithm avoids the necessity of evaluating the proximal map of the nonsmooth cost with respect to feedback controls. 
 This optimization enhances the efficiency of the algorithm; see Remark \ref{rem:referee2} for further details.}

\item
We present a practical implementation of the above accelerated gradient  algorithm
by combining Monte Carlo  and PDE approaches. 
At each iteration, the  state dynamics
with a given feedback control
is realized
by using 
 a particle approximation  and  Euler--Maruyama timestepping  scheme, 
and  the required gradient
is computed 
by solving  a  coupled system of nonlocal linear PDEs, whose coefficients depend on the empirical measure of the particle system.
The coupled PDE system 
is then solved 
with two  different approaches 
depending on the state dimension $d$ in 
\eqref{cforward},
in order to balance the efficiency and computation complexity.
In the low-dimensional setting (say $d\le 2$),
we discretize the coupled PDE system 
by a class of semi-implicit monotone finite difference approximations.
To accommodate the curse of dimensionality with large state dimension, 
 we  also propose 
   a residual approximation approach
   to 
 solve the coupled PDE system,  
 in which   the numerical solution is decomposed  into 
 a pre-determined candidate solution 
 and an unknown residual  term.
The computation of the residual term 
is addressed by a mesh-free method 
based on neural network approximation and 
 stochastic optimization algorithms, such as  the Stochastic Gradient Descent (SGD) algorithm or its variants.

The proposed algorithm combines the advantages of probabilistic and deterministic approaches. 
Firstly, 
the particle approximation 
 allows for efficient computation of
 the marginal distribution of the state process
 and avoids the
numerical challenge in
 solving a nonlinear FP equation
 (cf.,   \cite{achdou2015system, achdou2016mean, achdou2020mean,carmona2021convergence_I}).
 This is particularly relevant for high-dimensional MFC problems 
 with degenerate diffusion coefficients or irregular initial distribution
(see Section \ref{SEC:NUM}). 
Secondly, 
by exploiting the PDE formulation of the gradient evaluation,
our algorithm 
recovers  the optimal feedback control on the entire computational domain,
rather than merely along the trajectories of the optimal state 
(cf., the probabilistic methods in 
\cite{angiuli2019cemracs,carmona2019convergence,germain2019numerical,reisinger2020posteriori}).
This allows us to capture   important  structures of the optimal  control 
and achieve a robust performance with parameter uncertainty
(see Section \ref{SEC:NUM}).
As alluded to earlier,
such an accurate approximation of   the optimal feedback control  
is 
practically important for mathematical modelling and {fault diagnosis}  in engineering. 
Finally, instead of directly applying SGD to   high-dimensional PDE systems 
as in \cite{sirignano2018dgm, carmona2021convergence_I},
the proposed residual approximation approach leverages available efficient solvers
to compute  the dominant part of  solutions 
(see e.g., the Riccati-based solvers in \cite{gobet2019extended,albi2021moment}), 
and 
 employs  a small number of  SGD  iterations to fit 
 the  residual term.
This significantly 
accelerates the convergence of the algorithm for solving high-dimensional MFC problems
(see Figure \ref{fig:value_conv_6d_l1_ra}).

\item
We demonstrate the effectiveness of the algorithm
through extensive numerical experiments.
This includes 
a two-dimensional nonsmooth MFC problem
arising from 
portfolio liquidation  with trade crowding,
and 
a six-dimensional nonsmooth nonconvex MFC problem
arising from 
sparse consensus control   of stochastic Cucker--Smale models.
Our experiments show that
the resulting approximate feedback control 
correctly captures the temporal/spatial nonlinearity and  the sparsity
of the optimal control,
and achieves a robust performance in the presence  of
 initial state  perturbation.

\end{itemize}

The rest of the paper is organized as follows.
 Section \ref{SEC:NAG}
 describes our numerical
methodology, including the 
 accelerated proximal gradient  iteration,
 the particle system for the state process,
 and the PDE system for the gradient evaluation. 
We then 
propose a class of  finite difference approximations in Section \ref{SEC:FDM}
and neural network-based residual approximations in Section \ref{SEC:NN}
to solve the PDE systems. 
In Section \ref{SEC:NUM}, we present exhaustive numerical experiments for 
multidimensional nonsmooth nonconvex MFC problems,
which demonstrate that the proposed algorithm 
leads to more accurate and stable feedback controls
than the aforementioned PG method and the PMP method.

\section{Fast iterative PDE-based method for MFC problems}\label{SEC:NAG}

This section  proposes a class of Markovian accelerated proximal gradient methods for  MFC problems with nonsmooth running costs.

{\color{black} In the following subsections, we will incrementally introduce the main algorithm of the paper, in order to highlight the main ingredients and to stress which components (at least in special cases) are supported by theoretical results and which (at the present stage) are heuristic.\footnote{We thank an anonymous
referee for their comments that helped bring clarity to the presentation.}
In Section \ref{sec:GRADNM}, we will first formulate a gradient method with open-loop controls, without acceleration step and nonsmooth cost functional $\ell$. Section \ref{sec:GRADM} is dedicated to the formulation of a gradient method when the optimization is carried out over feedback functions. 

These two different approaches serve as motivation for Section \ref{sec:IPDE}, which introduces a heuristic iterative PDE based algorithm  to overcome certain deficiencies encountered in the first two formulations (both earlier sections assume $\ell \equiv 0$). Our main algorithm, which allows to incorporate nonsmooth cost and performs an additional moment step, is detailed in Section \ref{sec:FIPDE}. Finally, Section \ref{subsec:conv} discusses existing convergence results for these various gradient algorithms.  

\subsection{First motivating algorithm:   gradient descent in open-loop control}
\label{sec:GRADNM}
{\color{black}

This section considers MFC problems with regular running costs  ($\ell \equiv 0$ in \eqref{cvalue}),
and 
derives plain 
 gradient-descent iterations over open-loop controls.
 
If $\ell\equiv 0$,    \eqref{cforward}-\eqref{cvalue}  is  a   minimisation  problem
of $F:\cH^2(\R^k)\to \sR$
over the Hilbert space $\cH^2(\sR^k)$: 
\begin{equation}\label{eq:optimization}
\inf_{\alpha \in \mathcal{H}^{2}(\R^k)} F(\a),
\quad
\textnormal{
with
$F(\alpha)\coloneqq\E \bigg[
\int_0^T f(t,X^{\alpha}_t,\alpha_t,\mathcal{L}_{(X^{\alpha}_t,\alpha_t)}) \, \mathrm{d}t+g(X^{\alpha}_T,\mathcal{L}_{X^{\alpha}_T})
\bigg],
$
}
\end{equation}
where $X^\a\in \cS^2(\sR^d)$ is the state process controlled by $\a$ satisfying \eqref{cforward}. By 
\cite[Lemma 3.1]{acciaio2019extended}
and by assuming 
differentiability of $(b,\sigma, f)$, 
 $F $  is Fr\'{e}chet differentiable
with derivative
 $\nabla F:\cH^2(\R^k)\to \cH^2(\R^k)$
satisfying  
 for all $\a\in \cH^2(\sR^k)$,
 \bb\label{eq:F_gradient}
 {(\nabla F)(\a)_t}=
    (\partial_{a} H)(t,X^{\alpha}_t,\alpha_t, \mathcal{L}_{(X^{\alpha}_t,\alpha_t)}, Y^{\alpha}_t,Z^{\alpha}_t) 
 + \tilde{\mathbb{E}}[ (\partial_{\nu} H)(t,\tilde{X}^{\alpha}_t,\tilde{\alpha}_t, \mathcal{L}_{(X^{\alpha}_t,\alpha_t)}, \tilde{Y}^{\alpha}_t,\tilde{Z}^{\alpha}_t)(X^{\alpha}_t,\alpha_t)],
\ee
 where 
  $H:[0,T] \times \R^d \times \R^k \times  \mathcal{P}_2(\R^d \times \R^k)\times \R^d \times \R^{d \times n}\to \R$
is  the Hamiltonian defined by:
\begin{equation}\label{eq:mfcE_hamiltonian}
H(t,x,a,\eta,y,z)\coloneqq\left \langle b(t,x,a,\eta), y \right \rangle +\left \langle \sigma(t,x,a,\eta),z \right \rangle +f(t,x,a,\eta),
\end{equation}
 and $(Y^{\alpha},Z^{\alpha})$ are square integrable adapted  adjoint processes   such that 
for all $t\in [0,T]$,  
\begin{align}\label{adjointa}
\begin{split}
\mathrm{d}Y^{\alpha}_t&=-\Big((\partial_x H)(t,X^{\alpha}_t,\alpha_t,\mathcal{L}_{(X^{\alpha}_t,\alpha_t)},Y^{\alpha}_t,Z^{\alpha}_t) \\
&\qquad
+\tilde{\E}[(\partial_{\mu} H)(t,\tilde{X}^{\alpha}_t,\tilde{\alpha}_t,\mathcal{L}_{(X^{\alpha}_t,\alpha_t)},\tilde{Y}^{\alpha}_t,\tilde{Z}^{\alpha}_t)(X^{\alpha}_t,\alpha_t)]\Big)\,\d t
+Z^{\alpha}_t\,\d W_t,
 \\
 Y^{\alpha}_T&=(\partial_x g)(X^{\alpha}_T,\mathcal{L}_{X^{\alpha}_T})+\tilde{\E}[(\partial_{\mu} g)(\tilde{X}^{\alpha}_T,\mathcal{L}_{X^{\alpha}_T})(X^{\alpha}_T)].
 \end{split}
\end{align} 
Above and hereafter, following \cite{acciaio2019extended}, we 
use the tilde notation to denote an independent copy of a random variable.
Moreover, 
for a given function $h:\cP_2(\R^d\times \R^k)\to \R$
and a measure  $\eta\in \cP_2(\R^d\times \R^k)$
with marginals $\mu\in \cP_2(\R^d)$, $\nu\in \cP_2(\R^k)$,
 we  denote 
 by 
$\big((\p_\mu h)(\eta), (\p_\nu h)(\eta)\big)(\cdot):\R^n\times \R^k\to \R^n\times \R^k$
the partial L-derivatives of 
$h$ with respect to  the marginals;
see, e.g., \cite{acciaio2019extended, carmona2018probabilistic}
for detailed definitions.

Based on
the above interpretation,  
a gradient descent algorithm for optimal controls of \eqref{eq:optimization}
is given as follows:
 Let  $\alpha^{0}  \in \mathcal{H}^{2}(\mathbb{R}^{k})$
 be   the initial  guess of the optimal control,
 and $\tau>0$ 
  be a chosen stepsize. Consider the   sequence
  $(\alpha^m)_{m\in \N}\subset \cH^2(\sR^k)$ such that 
 \begin{align}\label{eq:gd_open}
 \alpha_t^{m+1}&= \alpha_t^m-\tau (\nabla F)(\alpha^{m})_t,
 \quad \forall m\in \N\cup\{0\},
\end{align}
{
$\mathrm{d}t \otimes \mathrm{d} \mathbb{P}$}-a.e.,
where $(\nabla F)(\alpha^{m})$  is defined by \eqref{eq:F_gradient}.
 

\paragraph{Pros and cons of Iteration \eqref{eq:gd_open}.} 
The main advantage of \eqref{eq:gd_open}
is that its convergence can be ensured given sufficient regularity of the functional 
$F:\cH^2(\sR^k)\to \sR$ in 
\eqref{eq:optimization}. 
Indeed, as a special case of Theorem 3.1 of \cite{beck2009fast} in conjunction with Remark 2.1 therein,
or by adapting  \cite[Theorem 2.1.14]{nesterov2003introductory}
to the present infinite-dimensional setting, one can deduce that, 
if $F:\cH^2(\sR^k)\to \sR$ is convex
and its derivative $\nabla F:\cH^2(\sR^k)\to \cH^2(\sR^k)$
 in \eqref{eq:F_gradient} is Lipschitz continuous,
then for sufficiently small stepsize,
the corresponding costs $(J(\alpha^m))_{m\in \N}$ of 
the sequence  $(\alpha^{m})_{m\in \N}$ 
converges to the optimal cost  $J^\star$ in \eqref{cvalue} (with $\ell \equiv 0$)
with the rate $\cO(m^{-1})$.
Such a regularity condition holds in particular 
for  commonly used linear-convex MFC problems
(see e.g., \cite{basei2017linear,acciaio2019extended,gobet2019extended,reisinger2020regularity}).
For nonconvex $F$, it is well-known that gradient-based algorithms in general can only find critical points.

However, it is difficult to implement the 
iteration \eqref{eq:gd_open} in practice. 
For each $m$,   since $\alpha^{m}$  
is a  non-Markovian stochastic process, 
the BSDE \eqref{adjointa} for  the   adjoint processes 
 $(Y^{m}, Z^m)$ is typically  non-Markovian.
This prevents us from 
  evaluating $(\nabla F)(\alpha^m)$
  and 
implementing the 
updates \eqref{eq:gd_open} in a pairwise sense
as for gradient-based algorithms for deterministic optimal control problems (see e.g., \cite{azmi2021optimal}).
To be more precise, let us initialize the iterates \eqref{eq:gd_open}
 with  $\a^0=0$. 
 Then one can express $(Y^1,Z^1)$ as 
 $(Y^0_t,Z^0_t)=(u(t,X^0_t),v(t,X^0_t))$  for some deterministic functions $u,v$ (often called decoupling fields), 
 and obtain the  stochastic processes $(Y^0,Z^0)$ by computing  the functions  $u,v$. 
Hence, one easily sees that $\a^1$ is a function of $X^0$, 
and consequently, the coefficients of the state dynamics \eqref{cforward} and adjoint equations \eqref{adjointa}
for  $(X^{1},Y^1,Z^1)$  would 
depend on both $X^{0}$ (through $\alpha^{1}$) and $X^{1}$. 
Repeating this process, 
one observes that 
$(Y^m,Z^m)$ are functions of time and 
 the enlarged system $(X^{0}, \ldots, X^{m})$,
 whose 
 computational complexity increases rapidly as the number of   iterations grows. 
 A similar difficulty has been observed  in \cite{bender2008time}
 for implementing (non-Markovian) Picard iterations to solve coupled FBSDEs.}

\subsection{Second motivating algorithm:   gradient descent in  feedback control}\label{sec:GRADM}
{\color{black}

In this section, we heuristically derive  a gradient descent  method 
for \eqref{eq:optimization}
by   restricting the optimisation over sufficiently regular feedback maps.

To this end, for each sufficiently regular feedback control  $\phi:[0,T] \times \mathbb{R}^d \to \mathbb{R}^k$,
let  $X^\phi $ be   the corresponding state process satisfying  the dynamics:
for all $t\in [0,T]$, 
 \begin{align}\label{sde_simple}
\mathrm{d} X_t=
b\left(t,X_t, \phi(t,X_t),\mathcal{L}_{(X_t,\phi(t,X_t))}\right)\, \mathrm{d} t
+\sigma\left(t,X_t, \phi(t,X_t),\mathcal{L}_{(X_t,\phi(t,X_t))}\right) \, \mathrm{d} W_t, 
\quad X_0=\xi.
\end{align}
Then 
\eqref{eq:optimization} 
(i.e.,  \eqref{cforward}-\eqref{cvalue} with
$\ell\equiv 0$) 
can be reformulated as minimizing 
\begin{align}\label{eq:cost_simple}
\widetilde{J}(\phi)\coloneqq 
\int_0^T \int_{\sR^d}  f\left(t,y,\phi(t,y),\mathcal{L}_{(X_t^{\phi},\phi(t,X_t^{\phi}))}\right)\, \mathcal{L}_{X_t^{\phi}}(\!\d y)  \, \mathrm{d}t+ \int_{\sR^d} g\left(y, \mathcal{L}_{X_T^{\phi}}\right)   \, \mathcal{L}_{X_T^{\phi}}(\!\d y)
\end{align}
over all sufficiently regular feedback controls $\phi$.
{
Recall that 
for any given feedback  map $\psi$, 
 the  derivative of $\widetilde J$  at $\phi$ 
 in the direction  $\psi$
(see e.g.,   \cite[Chapter 6.2.4]{carmona2018probabilistic})
is  
\begin{align}
 \label{eq:direction_derivative}
& \frac{\d  \widetilde{J}(\phi+\epsilon \psi )}{\d \epsilon}\bigg\vert_{\epsilon=0}
=\int_{0}^T \int_{\sR^d} 
\la 
\delta  \widetilde{J}_\phi  (t,y),  \psi(t,y)\ra  \, \mathcal{L}_{X_t^{\phi}}(\d y) \d t,
\end{align}
}
where 
$\delta \widetilde{J}_\phi: [0,T]\t \sR^d \to \sR^k$ is given by   
\begin{align}
\label{eq:gradient_feedback}
\begin{split}
&\delta \widetilde{J}_\phi  (t,x)
\\
& =(\partial_a \mathbb{H})\left(t,x,\phi(t,x),\mathcal{L}_{(X_t^{\phi},\phi(t,X^{\phi}_t))}, (\nabla_x u^{\phi})(t,x), (\textrm{Hess}_{x} u^{\phi})(t,x)\right)
 \\
 &\quad 
+ \tilde{\mathbb{E}}
\left[(\partial_{\nu}  \mathbb{H})\left(t,\tilde{X}_t^{\phi},\phi(t,\tilde{X}_t^{\phi}),\mathcal{L}_{(X_t^{\phi},\phi(t,X^{\phi}_t))}, (\nabla_x u^{\phi})(t,\tilde{X}_t^{\phi}), (\textrm{Hess}_{x} u^{\phi})  (t,\tilde{X}_t^{\phi})\right)(x,\phi(t,x))
\right], 
 \end{split}
\end{align}
  $\tilde{X}^\phi$  is  an independent copy of $\tilde{X}^{\phi}$,
 $\mathbb{H}:[0,T] \times \R^d \times \R^k \times  \mathcal{P}_2(\R^d \times \R^k)\times \R^d \times \R^{d \times d}\to \R$
is given by 
\begin{align*}
\mathbb{H}(t,x,a,\eta,y,z) \coloneqq\left \langle b(t,x,a,\eta), y \right \rangle +\frac{1}{2}\tr\left(  (\sigma\sigma^\trans)(t,x,a,\eta) z \right) +f(t,x,a,\eta),
\end{align*}
and 
$u^\phi: [0,T]\t \sR^d\to \sR$ satisfies the following linear PDE:
\begin{align}
\label{eq:gd_feedback_pde}
\begin{split}
(\partial_t u)(t,x) &= - \frac{\delta }{\delta \mu} \left[\int_{\sR^d}
\mathbb{H}(t,y,\phi(t,y),\mathcal{L}_{(X_t^{\phi},\phi(t,X_t^{\phi}))}, (\nabla_x u)(t,y), (\textrm{Hess}_{x} u)(t,y))\, \mathcal{L}_{X_t^{\phi}}(\!\d y)
 \right](x), \\
u (T,x) &= \frac{\delta }{\delta \mu} 
 \left[\int_{\sR^d}g(y,\mathcal{L}_{X_T^{\phi}})
\, \mathcal{L}_{X_T^{\phi}}(\!\d y)
 \right](x),
 \end{split}
\end{align}
with   
$\frac{\delta  }{\delta \mu}$  being the linear functional derivative with respect to the state law
$(\mathcal{L}_{X_t^{\phi}})_{t\in [0,T]}$
  (see \cite[Definition 5.43]{carmona2018probabilistic}).
  
{Now  
gradient descent updates of feedback controls
can be derived 
by selecting  a  function representation
of the directional derivative given  in \eqref{eq:direction_derivative}.
For instance, let $\phi^0$ be an initial feedback control, 
and $\tau>0$ be a stepsize.
For     all $m \in \mathbb{N} \cup \lbrace 0 \rbrace$,
suppose that   the law
$\mathcal{L}_{X_t^{\phi^m}} $  admits a  density 
$\frac{\d \mathcal{L}_{X_t^{\phi^m}} }{\d x}$
with respect to the Lebesgue measure for almost everywhere $t\in (0,T]$,
then  \eqref{eq:direction_derivative} can be written as:
\begin{align*}
& \frac{\d  \widetilde{J}(\phi^m+\epsilon \psi )}{\d \epsilon}\bigg\vert_{\epsilon=0}
=\int_{0}^T \int_{\sR^d} 
\left\la 
\delta  \widetilde{J}_{\phi^m}  (t,x) \frac{\d \mathcal{L}_{X_t^{\phi^m}} }{\d x}(x) ,  \psi(t,x)  \right \ra \d x \d t,
\end{align*}
which allows for interpreting   the function     $(t,x)\mapsto \delta \widetilde{J}_{\phi^m}(t,x)  \frac{\d \mathcal{L}_{X^{\phi^m}_t } }{\d x}(x) $ 
 as the gradient of $\widetilde{J}$  at $\phi^m$ with respect to         the canonical  inner product
on $L^2([0,T]\times \sR^d)$.
This suggests updating   $\phi^m$ by 
\begin{align}\label{eq:feedback_update}
\phi^{m+1}(t,x) &= \phi^{m}(t,x) -\tau \delta \widetilde{J}_{\phi^m} (t,x)\tfrac{\d \mathcal{L}_{X_t^{\phi^m}} }{\d x}(x),
 \quad (t,x) \in [0,T] \times \mathbb{R}^d,
\end{align}
which  is a direct application of gradient descent 
for   \eqref{eq:optimization} 
(or equivalently    \eqref{cforward}-\eqref{cvalue} with
$\ell\equiv 0$) 
 over feedback controls. 
 
%
}

\paragraph{Pros and cons of Iteration \eqref{eq:feedback_update}.}

 Iteration \eqref{eq:feedback_update} directly updates the feedback controls,
 and hence its complexity  
 at each iteration
 remains the same as the iteration proceeds.  
 It overcomes the drawback of  Iteration  \eqref{eq:gd_open}, 
 whose computational complexity increases as the number of iterations grows
 as illustrated at  
  the end of Section \ref{sec:GRADNM}.

However, Iteration \eqref{eq:feedback_update}   has several critical  deficiencies.  
{Compared with Iteration \eqref{eq:gd_open}, 
 the gradient direction in \eqref{eq:feedback_update} involves the density of the state process, which may not exist and can be computationally expensive to evaluate if the density function is irregular.  This irregularity arises in particular when the diffusion coefficient   of \eqref{sde_simple}   degenerates or when the initial state $\xi$ has a singular density.}
Moreover,  extending   Iteration  \eqref{eq:feedback_update} to general cost functionals \eqref{cvalue} with nonzero $\ell$   is more challenging. 
Observe that in the open-loop formulation,  
$\alpha\mapsto \sE\big[\int_0^T \ell(\alpha_t) \d t\big]$ 
in \eqref{cvalue} is    convex    due to the convexity of  $\ell$
and hence can be easily handled by the associated proximal map; see Section \ref{sec:FIPDE} for more details. 
However, 
due to 
the nonlinear dependence of $X^\phi$ on $\phi$,
 $\phi\mapsto \sE\big[\int_0^T \ell(\phi(t,X^\phi_t)) \d t\big]$ is typically nonconvex
 and evaluating the associated proximal map is computationally  expensive.

Furthermore, analysing the convergence of   \eqref{eq:feedback_update}
 is technically more challenging. This is because the mapping $\phi\mapsto \widetilde{J}(\phi)$ in \eqref{eq:cost_simple}  is typically nonconvex with respect to feedback controls, even for deterministic linear-quadratic problems  without mean-field interaction (see \cite[Proposition 2.4]{giegrich2022convergence} for a concrete example).
This implies that the convergence analysis of such PGMs is linked to analyzing nonasymptotic performance of gradient search for nonconvex objectives, which has always been one of the formidable challenges in optimization theory. 

In fact, it is unclear how to choose a suitable  function space to analyze Iteration \eqref{eq:feedback_update}.
By \eqref{eq:gradient_feedback} and \eqref{eq:feedback_update}, 
the regularity  of $\phi^{m+1}$ depends on the regularity of $\textrm{Hess}_{x} u^{\phi^m}$.
As $u^{\phi^m}$ satisfies \eqref{eq:gd_feedback_pde} whose coefficients involve $\phi^m$, 
it is expected that the regularity of $\textrm{Hess}_{x} u^{\phi^m}$ is controlled by the regularity of $\textrm{Hess}_{x} \phi^m$.
This suggests that 
 estimating the derivatives of $\phi^{m+1}$ requires bounds on higher order derivatives of $\phi^{m}$, and it is unclear how to close this norm gap.

}

\subsection{A heuristic modification  of  gradient methods in  feedback control}
\label{sec:IPDE}

{\color{black}
This section proposes a novel iterative method 
for MFC problems with regular running costs  ($\ell \equiv 0$ in \eqref{cvalue})
by combining the features of 
 Iterations  \eqref{eq:gd_open} and  \eqref{eq:feedback_update}.
On one hand, the algorithm directly updates 
the  feedback controls as Iteration \eqref{eq:feedback_update},   
and hence the   computational complexity   of the   adjoint variables do not change with respect to the number of iterations. 
On the other hand, it  heuristically replaces  the  
gradient direction  over feedback controls 
by the   functional derivative  in 
Iteration \eqref{eq:gd_open}. 
In particular, 
at the $m$-th iteration, we will update the present feedback map by evaluating the functional derivative of $J(\cdot;\xi_0)$ at the open-loop control $\a^{\phi^m}$ induced by the current policy $\phi^m$,
and then obtain the update direction based on a Markovian representation of the gradient.
The latter can be seen as a simplified version of Iteration  \eqref{eq:feedback_update} that does not require computing    the state process's density,    the derivative of feedback functions, or
the prox operator in the space of feedback functions in the case of non-smooth costs (see Remarks \ref{rem:referee} and \ref{rem:referee2}, respectively).

 We shall represent the control $\a^m$ as $\alpha^m_t=\phi^{m}(t,X^{m}_t)$ with some deterministic feedback function $ \phi^{m}:[0,T]\times \R^d\to \R^k$
at each iteration.
The update direction of $\phi^{m}$
is obtained via the functional derivative  in 
  \eqref{eq:gd_open}. 
More precisely, at the $m$-th iteration, 
given $\phi^m$, 
  consider the associated 
the controlled  state dynamics: $X_0 = \xi$, and for all $t\in [0,T]$, 
\begin{equation}\label{forward2}
\mathrm{d} X_t= b(t,X_t,\phi^{m}(t,X_t),\mathcal{L}_{(X_t,\phi^{m}(t,X_t))})\, \mathrm{d} t+\sigma(t,X_t,\phi^{m}(t,X_t),\mathcal{L}_{(X_t,\phi^{m}(t,X_t))})\, \mathrm{d} W_t.
\end{equation}
Let $X^m$ be the solution of \eqref{forward2};
we then seek the adjoint processes $(Y^m,Z^m)$ satisfying \eqref{adjointa} with $(X^\a, \a)$ replaced by 
$(X^m, \phi^m(\cdot, X^m))$,
in order to evaluate the gradient of $F$ at the corresponding control process $\phi^m(\cdot, X^m)$
as in \eqref{eq:F_gradient}.
The feedback structure of $\a^m$ implies that there exist deterministic decoupling fields 
$u^{m}:[0,T] \times \R^{d} \to \R^{d}$ and $v^{m}:[0,T] \times \R^{d} \to \R^{d\times n}$
such that
\bb\label{eq:YZ_m}
Y^m_t=u^m(t, X^m_t), \quad
Z^m_t=v^m(t, X^m_t),   
\quad 
\textnormal{$\mathrm{d}t \otimes \mathrm{d}\mathbb{P}$-a.e.}
\ee
Moreover, the functional derivative $\nabla F(\a^m)$ 
admits a Markovian representation  $(\nabla F)(\phi^m(\cdot, X^m))_t=(\nabla F)(\phi^{m})(t,X^m_t)$, $\mathrm{d}t \otimes \d \mathbb{P}$-a.e.,
with the function  $(\nabla F)(\phi^{m}):[0,T]\t \R^d\to \R^k$ defined by
\begin{align}\label{eq:F_gradient_fb}
\begin{split}
&(\nabla F)(\phi^{m})(t,x)
\\
&\coloneqq
 (\partial_{a} H)(t,x,\phi^{m}(t,x), \mathcal{L}_{(X^{m}_t,\phi^m(t, X^m_t))}, u^m(t,x),v^m(t,x)) 
 \\
&\quad 
+
\tilde{\E}[(\partial_{\nu} H)(t,\tilde{X}^{m}_t,\phi^m(t,\tilde{X}^{m}_t),\mathcal{L}_{(X^{m}_t,\phi^m(t,X^{m}_t))},
u^m(t,\tilde{X}^{m}_t),v^m(t,\tilde{X}^{m}_t))(x,\phi^m(t,x))],
\end{split}
\end{align}
where $H$ is defined in \eqref{eq:mfcE_hamiltonian}.
Hence, once the  functions $(u^m,v^m)$
are determined, one can compute the Markovian representation \eqref{eq:F_gradient_fb}
of the functional derivative $(\nabla F)(\phi^{m})$,
along with a gradient descent step can be applied to obtain
the updated      feedback map $\phi^{m+1}$; see   \eqref{eq:phi_psi2}.

To compute the decoupling fields $(u^m,v^m)$ in \eqref{eq:YZ_m} for each $m$, 
we connect them with solutions to some PDE system. 
In fact,  
the nonlinear Feynman--Kac formula in \cite{pardoux1998backward}
shows that, if $u^m$ is sufficiently smooth, 
then $v^m(t,x)=(\nabla_x u^m) (t,x)\sigma(t,x,\mathcal{L}_{X^{m}_t})$
for all $(t,x)\in [0,T]\t \R^d$, 
and $u^m$  
solves 
 the following  system of parabolic linear PDEs
 (depending on $\phi^m$ and the law of $X^m$):
\begin{align}\label{PDEsystem}
(\partial_t u)(t,x)+(\mathscr{L}^m u)(t,x)= -f^{m}(t,x,u,v),
\end{align}
with $u(T,x)=h^m(x)$
and $v(t,x)=(\nabla_x u) (t,x)\sigma^m(t,x)$
 for all $(t,x)\in [0,T]\t \R^d$, 
where
$\mathscr{L}^m$ is the (vector-valued)  differential operator  such that
 for each $\varphi=(\varphi_1, \ldots, \varphi_d)\in C^{1,2}([0,T]\times \R^d;\R^d)$, $i\in \{1,\ldots, d \}$ and $(t,x)\in [0,T]\times \R^d$,
\bb
\label{eq:cL_m}
 (\mathscr{L}^m \varphi)_i(t,x)=\frac{1}{2}\text{tr}\big(\sigma^{m}(t,x)\sigma^{m}(t,x)^\trans (\textrm{Hess}_x\varphi_i)(t,x)\big)+ \langle  b^{m}(t,x), (\nabla_x \varphi_i)(t,x)\rangle, 
\ee
 with
 $b^{m}:[0,T]\t \R^d\to \R^d$ and 
 $\sigma^{m}:[0,T]\t \R^d\to \R^{d\t n}$
 given by
\begin{align}\label{eq:b_sigma_m}
\begin{split}
b^{m}(t,x)&=b(t,x,\phi^{m}(t,x),\mathcal{L}_{(X^{m}_t,\phi^{m}(t,X^{m}_t) )}),  \\
\sigma^{m}(t,x) &=\sigma(t,x,\phi^{m}(t,x),\mathcal{L}_{(X^{m}_t,\phi^{m}(t,X^{m}_t) )}),
\end{split}
\end{align}
and the functions 
 $f^m$  
 and $h^m$  satisfy for all $(t,x)\in [0,T]\t \sR^d$
 (cf.~\eqref{adjointa}),
\begin{align}
f^{m}(t,x,u,v)
&= (\partial_x H)(t,x,\phi^m(t,x),\mathcal{L}_{(X^{m}_t,\phi^m(t,X^{m}_t))},u(t,x),v(t,x))
\label{eq:fm}
 \\
&\quad
+\tilde{\E}[(\partial_{\mu} H)(t,\tilde{X}^{m}_t,\phi^m(t,\tilde{X}^{m}_t),\mathcal{L}_{(X^{m}_t,\phi^m(t,X^{m}_t))},
u(t,\tilde{X}^{m}_t),v(t,\tilde{X}^{m}_t))(x,\phi^m(t,x))],
\nonumber
\\
h^m(x) 
&= (\partial_x g)(x,\mathcal{L}_{X^m_T})+\tilde{\E}[(\partial_{\mu} g)(\tilde{X}^m_T,\mathcal{L}_{X^m_T})(x)].
\label{eq:h}
\end{align}


Algorithm \ref{alg:NAG1} summarizes 
 the above gradient descent method for \eqref{cforward}-\eqref{cvalue} (with $\ell \equiv 0$),
 which will   be referred to as the iterative PDE-based (IPDE) method hereafter.
 
\begin{algorithm}[H]
\caption{Iterative PDE-based   method for MFC problems}\label{alg:NAG1}
\begin{algorithmic}[1]
\STATE 
\textbf{Input:} Choose the initial  feedback  control $\phi^{0}:[0,T] \times \R^{d} \mapsto \R^{k}$
  and stepsize $\tau>0$.  
\FOR{$m=0, 1 \ldots $} 
\STATE
Compute the law of the state process $X^m$ 
governed by \eqref{forward2}.
\STATE
Compute the decoupling fields $(u^m, v^m)$ of  $(Y^m,Z^m)$ 
by
 solving \eqref{PDEsystem} based on the law of $X^m$.
\STATE
Update the controls such that for all $(t,x)\in [0,T]\t \R^d$, 
\begin{align}\label{eq:phi_psi2}
\begin{split}
 \phi^{m+1}(t,x)&= \phi^m(t,x)-\tau 
 (\nabla F)(\phi^{m})(t,x)
\end{split}
\end{align}
 with
$(\nabla F)(\phi^{m})$  defined by \eqref{eq:F_gradient_fb}.
\ENDFOR

\end{algorithmic}
\end{algorithm}  
  
\begin{remark}\label{rem:referee}
For the avoidance of doubt, the update \eqref{eq:phi_psi2}   is neither a gradient step
for the stochastic formulation of the problem using a Markovian representation of the open-loop controls,
nor is it a gradient step with respect to the feedback map. 
It is not the former 
because the update \eqref{eq:phi_psi2}  is not equivalent to the corresponding gradient step \eqref{eq:gd_open}, since subtracting   feedback controls is not equivalent to subtracting the corresponding stochastic processes.
It is not the latter 
because
the update direction does not 
{involve the density of
$ { {X_t^{\phi^m}} }$},
and 
 $(u^m,v^m)$ in \eqref{eq:F_gradient_fb}
is different 
from $( \nabla_x u^{\phi^m}, \textrm{Hess}_{x} u^{\phi^m})$ in 
\eqref{eq:gradient_feedback}, 
as they are solutions to different PDEs. 
To see this, assume for simplicity that all coefficients are independent of the measure components.
Then  the PDE \eqref{PDEsystem} for $u^m$ simplifies to 
\begin{align}\label{PDEsystem_noMF}
(\partial_t u)(t,x)+(\mathscr{L}^m u)(t,x)= -(\partial_x H)(t,x,\phi^m(t,x), u(t,x),(\nabla_x u) (t,x)\sigma^m(t,x) ).
\end{align}
On the other hand, formally  differentiating the PDE  \eqref{eq:gd_feedback_pde}  for   $u^{\phi^m}$ with respect to $x$
suggests that 
$\nabla_x u^{\phi^m}$ satisfies the following equation:
\begin{align}
\label{eq:gd_feedback_pde_noMF}
\begin{split}
(\partial_t u)(t,x) +(\mathscr{L}^m u)(t,x) &= -  
(\p_x \mathbb{H}^m)(t,x,   u(t,x), (\nabla_{x} u)(t,x)),  
 \end{split}
\end{align}
where 
$\mathbb{H}^m:[0,T]\times \sR^d\times \sR^d\times\sR^d\to \sR$ is defined by
\begin{align*}
\mathbb{H}^m(t,x,y,z) \coloneqq\left \langle b(t,x,\phi^m(t,x)), y \right \rangle +\frac{1}{2}\tr\left(  (\sigma\sigma^\trans)(t,x,\phi^m(t,x)) z \right) +f(t,x,\phi^m(t,x)).
\end{align*}
It is clear from  the chain rule   that 
$\p_x \mathbb{H}^m$ consists of both the term $ \partial_x H$ in \eqref{PDEsystem_noMF}
and an additional term involving the derivative $\partial_x \phi^m$.

The above observation   indicates  that  the modified gradient update 
\eqref{eq:phi_psi2}
imposes weaker regularity requirements  on the iterate $  \phi^m$ than
the vanilla gradient descent  \eqref{eq:feedback_update}, 
as  it avoids evaluating  
the density of ${X^{\phi^m} }$ and 
the derivatives of $  \phi^m$ 
in the  computation of the gradient  $(\nabla F)(\phi^m)$. 
This difference in gradient evaluation  not only
enhances the stability of \eqref{eq:phi_psi2} (by preventing numerical instability when
computing an irregular density and 
 differentiating the feedback controls), 
but is also essential for analyzing the convergence of \eqref{eq:phi_psi2}; see Section \ref{subsec:conv}
for   details.
\end{remark}

  Concrete numerical methods for solving \eqref{PDEsystem} 
  (and consequently the functions  $(u^m, v^m)$)
will be given in 
Sections \ref{SEC:FDM} and  \ref{SEC:NN}. Here, we point out the following three features of the PDE system \eqref{PDEsystem},
which are crucial in the design of numerical methods:
(i) 
As the  Hamiltonian  $H$ defined in \eqref{eq:mfcE_hamiltonian}
is affine 
in the components $y$ and $z$,
the function $f^m$ is 
{linear} in $u$ and $v$,
and hence  \eqref{PDEsystem} is a linear PDE system. 
(ii)
Due to the measure dependence  in   $f^m$,
\eqref{PDEsystem} is nonlocal 
in the sense that the value of the solution $u^m$ at each point 
evolves based on the weighted average of other values of  $u^m$ 
with respect to the marginal laws $(\cL_{X^m_t})_{t\in [0,T]}$
 of the  process $X^m$.
(iii)
Even though the $i$-th component of the differential operator $\mathscr{L}^m$ 
in \eqref{eq:cL_m} only involves the $i$-th component of the solution $u^m$, 
the  function $f^m$ in 
 \eqref{PDEsystem} in general results in a coupling among  all components of $u$ and their gradients $v$.
In the special cases  where
the diffusion coefficient $\sigma$ of \eqref{cforward}   depends only on time
  (see Section \ref{SEC:NUM} for concrete examples),  
  the  function $f^m$ is independent of $v$ and hence 
the system  \eqref{PDEsystem}   
is only coupled through the solution  $u^m$.

 We emphasize  that,  
a key feature of the IPDE method
is that it 
 computes  the functions  $(u^m, v^m)$ 
for each $(t,x)\in [0,T]\t \sR^d$ based on a PDE formulation,
in contrast to the pure data-driven algorithms in
\cite{angiuli2019cemracs,carmona2019convergence,germain2019numerical,reisinger2020posteriori},
which solve
$(u^m, v^m)$ merely along the trajectories of $X^m$.
Hence, the IPDE method leads to a more accurate approximation 
of the optimal feedback control,
especially
outside the support of the optimal state process
of \eqref{cforward}-\eqref{cvalue}. 
In particular,
the IPDE method
is capable of
recovering  important structural features of the optimal  control,
and 
the performance of the  approximate feedback control 
is robust with respect to perturbation of model parameters;
see    Section \ref{SEC:NUM} for a  detailed comparison between the (F)IPDE method and several  data-driven algorithms.
}

\subsection{A heuristic acceleration of    gradient methods with nonsmooth cost}\label{sec:FIPDE}

{\color{black} 

This section extends Algorithm \ref{alg:NAG1} 
in Section \ref{sec:IPDE} to MFC problems 
\eqref{cforward}-\eqref{cvalue} with   nonsmooth costs and further  accelerates the algorithm  convergence   by incorporating  an additional moment step. 
To handle the nonsmooth term $\ell :\sR^k\to \sR\cup\{\infty\}$, 
for each $\tau>0$, we consider   the following    proximal map 
 $\prox_{\tau\ell}:\sR^k\to \sR^k$: 
  $$\textrm{prox}_{\tau\ell}(a)=\arg\min_{z\in \sR^k}\left(\frac{1}{2}|z-a|^2+\tau\ell(z)\right),
  \quad a\in \sR^k,$$
  which is well-defined as $\ell$ is proper, lower semicontinuous and convex. 
  For many practically  important nonsmooth functions $\ell$,
the proximal function $\prox_{\tau \ell}$ can be 
 evaluated 
 either analytically
 (see e.g.,~\cite[Chapter 6]{beck2017first}) or approximately by efficient numerical methods
 (see \cite{schmidt2011convergence} and the references therein).
  The acceleration 
 of \eqref{eq:phi_psi2}
  is   inspired by Nesterov's accelerated proximal gradient (NAG) method 
{\color{black}
(also known as 
the Fast Iterative Shrinkage-Thresholding
Algorithm  \cite{beck2009fast}}) tailored
to the present setting.
The algorithm is summarized as follows
and 
   will   be referred to as the fast iterative PDE-based (FIPDE) method.

\begin{algorithm}[H]
\caption{Fast Iterative PDE-based method for MFC problems}\label{alg:NAG2}
\begin{algorithmic}[1]
\STATE 
\textbf{Input:} Choose the initial  feedback  control $\phi^{0}:[0,T] \times \R^{d} \mapsto \R^{k}$
  and stepsize $\tau>0$. Set $\psi^{0}=\phi^{0}$.  
\FOR{$m=0, 1 \ldots $} 
\STATE
Compute the law of the state process $X^m$ 
governed by \eqref{forward2}.
\STATE
Compute the decoupling fields $(u^m, v^m)$ of  $(Y^m,Z^m)$ 
by
 solving \eqref{PDEsystem} based on the law of $X^m$.
\STATE
Update the controls such that for all $(t,x)\in [0,T]\t \R^d$, 
\begin{subequations}
\label{eq:phi_psi}
\begin{alignat}{1}
 \phi^{m+1}(t,x)&= \text{prox}_{\tau\ell} \left(\psi^m(t,x)-\tau 
 (\nabla F)(\psi^{m})(t,x) \right), 
 \label{eq:phi_psi_prox} \\
  \psi^{m+1}(t,x)&=\phi^{m+1}(t,x) + \frac{m}{m+3}(\phi^{m+1}(t,x)-\phi^m(t,x)),
   \label{eq:phi_psi_nag}
\end{alignat}
\end{subequations}
 with
$(\nabla F)(\psi^{m})$  defined by \eqref{eq:F_gradient_fb}.
\ENDFOR

\end{algorithmic}
\end{algorithm}

 The ratio $\frac{m}{m+3}$
in \eqref{eq:phi_psi_nag}
is often referred to as the  momentum parameter,
whose selection adheres to the general guideline provided by \cite{beck2009fast}.\footnotemark
 \footnotetext{
 {\color{black}
By   \cite[Theorem 4.4]{beck2009fast}, 
  NAG method (or FISTA)  
   with  the momentum step 
 $
  y^{m+1}=x^{m+1} +\frac{ \theta_{m+1}-1}{\theta_{m+2}}(x^{m+1}-x^m)
$  converges  for convex optimisation problems, 
 provided  that 
$\theta_{m+1} \ge 1$
and 
$\theta_m^2\ge \theta_{m+1}^2 -\theta_{m+1}$ 
for all $m \ge 0$.
For all 
$m\in \N$, \eqref{eq:phi_psi_nag} corresponds to   
$\theta_{m}= {(m+1)}/{2}$,
while 
 \cite{beck2009fast} 
 chooses
{
$\theta_{m+1}\ge 1$ 
}
such that
$\theta_m^2= \theta_{m+1}^2 -\theta_{m+1}$.}
}

\begin{remark}
\label{rem:referee2}

Further to Remark \ref{rem:referee} for  Algorithm \ref{alg:NAG1},
Algorithm \ref{alg:NAG2} is a heuristic combination  of  the NAG method applied to 
  the stochastic formulation of the control problem over open-loop controls, and to  the deterministic reformulation  over feedback  controls 
  discussed in Section \ref{sec:GRADM}.
The computation of the prox
in \eqref{eq:phi_psi_prox} 
 is the one for the formulation of the problem with stochastic processes.
It aims to overcome the difficulty in computing a suitable prox for feedback maps  
caused by the lack of convexity
of 
  $\phi\mapsto \mathbb{E}\big[\int_0^T \ell(\phi(t,X^\phi_t))\d t\big]$    (see   the end of Section \ref{sec:GRADM}).   
In particular, 
we have   
  exploited
   the structure of the nonsmooth functional $G: \alpha \mapsto \mathbb{E}\big[\int_0^T \ell(\alpha_t)\d t\big] $
and explicitly expressed 
the  proximal   of  $G:\cH^2(\R^k)\to \sR\cup\{\infty\}$  
  via a pointwise composition of the proximal   of $\ell$.
 The addition/subtraction 
 of controls 
 in \eqref{eq:phi_psi_prox}  and \eqref{eq:phi_psi_nag} 
 is the addition/subtraction of feedback maps and is not equivalent to the addition/subtraction of the corresponding processes.
 \end{remark}

 In practice, we  represent the feedback maps $(\phi^m,\psi^m)_{m\in \N}$
  and the decoupling fields $(u^m,v^m)_{m\in \N}$ in suitable parametric forms, 
  whose precise choices depend on 
   the dimension of the problem and the numerical methods used to solve the PDE system
   \eqref{PDEsystem}
   (see Sections \ref{SEC:FDM} and \ref{SEC:NN} for details).
   Given a parameterized 
 feedback map $\psi^{m}$,
 the law of the controlled state process $X^m$ 
in Step 3 
of Algorithm \ref{alg:NAG2}
can be approximated 
by  a particle method  and an Euler--Maruyama discretization of  
\eqref{forward2}.
For instance, 
let 
 $N \in \N$ be the number of particles
and  $\{0=t_0<\ldots<t_M=T\}$
a partition of  $[0,T]$ with time stepsize $\Delta t=T/M$ for some  $M \in \mathbb{N}$.
Then we consider the discrete-time interacting particle system $(X^{l,N}_t)_{t\in [0,T]}$, $l=1,\ldots, N$,  
such that 
$X^{l,N}_0=\xi^l$,   for   $j =  0, \ldots, M-1 $ and $t\in [t_j,t_{j+1})$,
$X^{l,N}_t=X^{l,N}_{t_j}$, and 
\begin{align}\label{eq:forward_samples}
X^{l,N}_{t_{j+1}} = X^{l,N}_{t_{j}} + b(t_j,X^{l,N}_{t_j},\psi^m(t_j,X^{l,N}_{t_j}),\mu^{X,\psi^m}_{N,t_j}) \Delta t + \sigma(t_j,X^{l,N}_{t_j},\psi^m(t_j,X^{l,N}_{t_j}),\mu^{X,\psi^m}_{N,t_j}) \Delta W_j^{l},
\end{align}
with $\mu^{X,\psi^m}_{N,t_j}(\mathrm{d}x)= \frac{1}{N} \sum_{l=1}^{N} \delta_{(X^{l,N}_{t_j},\psi^m(t_j,X^{l,N}_{t_j}))}(\mathrm{d}x) $,
where 
$(\xi^{l})_{l =  1, \ldots, N}$
and 
 $(W^{l})_{l =  1, \ldots, N}$ are   independent copies of  $\xi$ and $W$, respectively,
  $\Delta W_j^{l} = W^{l}_{t_{j+1}} - W^{l}_{t_{j}}$ for all $j,l$, and 
  $\delta_x$ denotes the  Dirac measure supported at $x$ for all $x\in \R^d$.

Compared with the deterministic methods in 
  \cite{achdou2015system, achdou2016mean, achdou2020mean,carmona2021convergence_I},
\eqref{eq:forward_samples} 
 avoids the
numerical challenge in
 solving the nonlinear FP equation,
 and allows for efficient computation of
 the marginal distribution of the state process,
 especially for MFC problems \eqref{cforward}-\eqref{cvalue}
 with
degenerate diffusion coefficient $\sigma$ 
or initial state $\xi$ with singular density. 
In fact, 
it is well-known that, for sufficiently regular   feedback maps $\psi^{m}$,
 the empirical measure $\mu^{X,\psi^m}_N$ converges to the law of $(X^m,\psi^{m}(\cdot, X^m))$ 
 in the Wasserstein metric as $M,N$ tend to infinity;
see also \cite{bossy1997stochastic,bao2021first} for the convergence rates of 
\eqref{eq:forward_samples} in terms of $M$ and $N$.
For our numerical experiments in Section \ref{SEC:NUM}, we choose  sufficiently large  $M$ and $N$ such that
the presented results are not influenced by those choices.

As numerical approximations of   \eqref{forward2} 
are relatively well-understood, 
in the subsequent sections, we  focus on the numerical approximations of 
the PDE system  \eqref{PDEsystem} 
with given 
parameterized 
 feedback function $\psi^m$
 and  empirical approximations $\mu^{X,\psi^m}_{N,\cdot}$ of the  law of $(X^m,\psi^{m}(\cdot, X^m))$. 
In particular, we shall 
propose a class of monotone schemes in Section \ref{SEC:FDM}
for  \eqref{PDEsystem}  with  spatial  dimension $d\le 3$,
and 
neural network-based schemes in Section \ref{SEC:NN}
for the high-dimensional setting.  
}

 \subsection{Convergence of Algorithms \ref{alg:NAG1} and \ref{alg:NAG2}}\label{subsec:conv}

{\color{black}

We begin by recalling the discussion of convergence for the NAG method 
with open-loop controls at the end of Section \ref{sec:GRADNM}. 
Furthermore, as shown in 
Theorem 4.4 of \cite{beck2009fast},
including a momentum step after \eqref{eq:gd_open}
yields an improved  convergence rate 
 $\cO(m^{-2})$  for minimising  smooth convex   functionals $J$,
 which  is the optimal convergence rate of
 gradient-based algorithms
(see \cite{beck2009fast} and references therein).
The same convergence rate also holds when convex nonsmooth costs are incorporated through 
an application of a proximal operator. Additionally, if
the controls are updated 
with approximate proximal operators,
one   can still recover the same convergence rate,
 provided that the errors made 
 in  the calculation of the proximal operator
 decrease at appropriate rates \cite{schmidt2011convergence}. 
For nonconvex functionals,
it is  known that the momentum step 
can produce highly oscillatory solutions.  
Additional monitoring steps  
are typically necessary to 
 ensure  a sufficient descent of 
the function values and the convergence  of the iterates to critical  points 
  \cite{beck2009fast_denoise,li2015accelerated, ghadimi2016accelerated, li2017convergence}. 

Unfortunately, existing  convergence results of the NAG method 
cannot be applied to Algorithm  \ref{alg:NAG2}. 
As pointed out in Remarks \ref{rem:referee} and \ref{rem:referee2},
Algorithm \ref{alg:NAG2} is not a direct application of the NAG method (or   FISTA) to the stochastic formulation of the problem \eqref{cforward}-\eqref{cvalue} nor to the deterministic reformulation  discussed in Section \ref{sec:GRADM}. 
The gradient direction  and the proximal map in \eqref{eq:phi_psi_prox} 
are computed using   the stochastic formulation over open-loop controls,
whereas the  updated policy is obtained forming an addition of feedback functions which, in general, is different to the corresponding combination of open-loop controls.
Hence, although 
our numerical results in Section \ref{SEC:CS}
indicate that  
 Algorithm \ref{alg:NAG2}   performs reasonably well
even  for nonconvex control problems,
a    convergence analysis of Algorithm \ref{alg:NAG2}   remains an open question. 

Note that 
a comparable discrepancy between two formulations arises even in the absence of the momentum step \eqref{eq:phi_psi_nag}  (i.e., $\psi^m=\phi^m$ for all $m$).
In    the   recent work \cite{reisinger2022linear},
such an algorithm has been analyzed via a contraction argument 
 for specific drift-controlled nonsmooth control problems without mean-field interaction. 
We identify conditions under which  the  algorithm 
generates  uniform Lipschitz feedback controls $(\phi^m)_{m\in \N}$, and further establish that the associated control processes 
converge   to a critical point of the  functional $J$ in \eqref{cvalue}.
The convergence result    holds    when the diffusion coefficient is uncontrolled, and in addition
one of the following five cases is satisfied: 
\begin{enumerate*}[(i)]
\item time horizon $T$ is small;
\item running cost is sufficiently convex in control;
\item costs depend weakly on state;
\item control affects state dynamics weakly;
\item state dynamics is strongly dissipative.
\end{enumerate*}
Note that
the conditions allow for nonlinear state dynamics with degenerate noise, and   nonconvex and nonsmooth cost functions.
We refer the reader to \cite{reisinger2022linear} 
for rigorous statements of the convergence results.

The main result in \cite{reisinger2022linear} shows that  in the case without mean-field interaction and without the momentum step, 
convergence of Algorithm  \ref{alg:NAG2}  to stationary points can   be proven, 
despite the aforementioned mismatch appearing in our formulation. 
This in particular shows the convergence of     Algorithm \ref{alg:NAG1} without mean-field interaction.
The proof  exploits  regularity estimates for the decoupling fields $(u^m,v^m)_{m\in\N}$ in \eqref{eq:YZ_m}
and does not require the convexity of the cost landscape. We conjecture that similar convergence results also hold in the mean-field setting, whose rigorous proof  require establishing the  regularity of the decoupling fields in the measure component, and is left for future research. 
 However, at the present stage, it is unclear if this mismatch when applying an additional momentum step \eqref{eq:phi_psi_nag}  in the policy update can be overcome, even under the assumptions of \cite{reisinger2022linear}.

}}
\section{Implementation of the FIPDE method via finite differences}\label{SEC:FDM}

In this section, we discuss the practical implementation of the FIPDE method in Section \ref{SEC:NAG} 
for low-dimensional state dynamics \eqref{cforward}.
In particular, 
we shall 
propose a class of  semi-implicit  monotone finite difference approximations for solving   \eqref{PDEsystem}
at the $m$-th NAG iteration,
which  achieves an efficient performance in terms of the computation time
if \eqref{cforward} has  a spatial dimension 
 $d\le 2$.

Throughout this section, 
we focus on the  $m$-th NAG iteration with $m\in \N\cup\{0\}$
and assume the
feedback function $\psi^m$
 and   particle approximation $(X^{l,N})_{l=1,\ldots, N}$
 of $X^m$ 
 are  given
 (see the discussion below Algorithm \ref{alg:NAG2}).
Then, 
we need
to solve the following  nonlocal parabolic PDE system 
 (cf.~\eqref{PDEsystem}):
 \begin{align}\label{eq:PDEsystem_em}
(\partial_t u)(t,x)+(\mathscr{L}^m_N u)(t,x)= -f^{m}_N(t,x,u,v),
\end{align}
with $u(T,x)=h^m(x)$
and $v(t,x)=(\nabla_x u) (t,x)\sigma^m(t,x)$
 for all $(t,x)\in [0,T]\t \R^d$,
 where the operator $\mathscr{L}^m_N$ (resp.~the function $f^m_N$) is defined similar to \eqref{eq:cL_m}
 (resp.~\eqref{eq:fm}), but depends on  the empirical measure
 $\mu^{X,\psi^m}_{N,t}(\mathrm{d}x)= \frac{1}{N} \sum_{l=1}^{N} \delta_{(X^{l,N}_{t},\psi^m(t,X^{l,N}_{t}))}(\mathrm{d}x)$
instead of the law $\cL_{(X^m_t,\psi^m(t,X^m_t))}$ 
for all $t\in [0,T]$.
%
To simplify the presentation,
we shall 
focus on the  uniform spatial grid   $\{x_k\}_k=h\Z^d$ on $\R^d$
 with mesh size $h>0$
  and a time partition $\{t_j\}_{j=0}^M$ with 
time stepsize $\Delta t=T/M$ for $M\in \N$, 
but similar schemes can be designed  for unstructured nondegenerate  grids as well.

We start by introducing  a semi-implicit  timestepping approximation to \eqref{eq:PDEsystem_em}.
Observe that
the $i$-th component of the differential operator 
$\mathscr{L}^m_N$ 
depends only on  the $i$-th component of the solution $u$,
and all nonlocal and coupling terms appear in $f^m_N$.
Hence, we shall adopt implicit timestepping for the  local operator $\mathscr{L}^m_N$ and explicit timestepping for $f^m_N$, which leads to the following (backward) time discretization of  \eqref{eq:PDEsystem_em}:
$U^M(x)=h^m(x)$ for all $x\in\sR^d$ and for all $j=1,\ldots, M$, $x\in \R^d$,
$$
\frac{U^{j}(x)-U^{j-1}(x)}{\Delta t}
+(\mathscr{L}^m_N U^{j-1})(x)= -f^{m}_N(t_{j},x,U^j,V^j),
$$
where 
for all $x\in \sR^d$,
$U^j(x)$   is  the  approximation of $u(t_j,x)$
and $V^j(x)=(\nabla_x U^j) (x)\sigma^m_N(t_j,x)$.
Note that 
the implicit timestepping for $\mathscr{L}^m_N$
  enables us to enjoy a less restrictive stability condition than 
that for fully explicit schemes,
while the explicit timestepping for $f^m_N$
avoids solving the dense system resulting from the mean-field terms,
and allows us to solve for each component of $U^j$ independently, given $U^{j+1}$.

We proceed to perform spatial discretization of $\mathscr{L}^m_N$. 
Note that  
the $i$-th component of $\mathscr{L}^m_N$ 
depends only on  the $i$-th component of $u$,
and all  components of 
$\mathscr{L}^m_N$ 
have the same coefficients (see \eqref{eq:cL_m} and \eqref{eq:b_sigma_m}).
Then, 
as  shown in \cite{biswas2010difference}, 
one can construct monotone and consistent approximations of $\mathscr{L}^m_N$
such that for any 
$\varphi=(\varphi_1, \ldots, \varphi_d)\in C^{2}( \R^d;\R^d)$, $i=1,\ldots, d$,  $j=0,\ldots, M$, $k\in \Z^d$,
$
L^{m,j}_{N,h}[\varphi]_{k}=(L^{m,j}_{N,h}[\varphi]_{1,k},\ldots, L^{m,j}_{N,h}[\varphi]_{d,k})
$
satisfies 
\begin{align}
\label{eq:L_m_h}
\begin{split}
L^{m,j}_{N,h}[\varphi]_{i,k}&=
\sum_{q\in \Z^d} a^{m,N}_{h,j,q,k}[\varphi_i(x_q)-\varphi_i(x_k)], 
\\
|
(\mathscr{L}^m_N\varphi)_i(t_j,x_k)
&-
L^{m,j}_{N,h}[\varphi]_{i,k}|\to 0,
\quad \textnormal{as $h\to 0$,}
\end{split}
\end{align}
with coefficients $a^{m,N}_{h,j,q,k}\ge 0$ for all $q,k,j$. 
The  precise construction of such numerical approximations 
depends on the  structures of the coefficients $b^m$ and $\sigma^m$.
In particular, one can adopt the  standard finite difference schemes  in \cite{barles1991convergence, kushner2001numerical} if the diffusion coefficient  ${\sigma}^m({\sigma}^m)^\trans$ is diagonally dominant, and use the semi-Lagrangian scheme in \cite{debrabant2013semi} for general cases.
We refer the reader to Section \ref{SEC:NUM} for more details.

It remains to discretize the term $f^m_N(t,x,u,v)$.
By \eqref{eq:mfcE_hamiltonian} and \eqref{eq:fm},
we have  
$$
f^{m}_N(t,x,u,v)
=
f^{m,\textrm{re}}_N(t,x,u)
+
f^{m,\textrm{ex}}_N(t,x,u),
$$
with the terms  
\begin{align*}
f^{m,\textrm{re}}_N(t,x,u)
&\coloneqq
 \big((\p_x b)(t,x,\psi^m(t,x),\mu^{X,\psi^m}_{N,t})\big)^\trans u(t,x) 
 +(\p_x f)(t,x,\psi^m(t,x),\mu^{X,\psi^m}_{N,t})
 \nonumber\\
&  \quad 
+
\sE^{\mu^{X}_{N,t}}
\Big[
 \big((\p_\mu b)(t,\cdot,\psi^m(t,\cdot),\mu^{X,\psi^m}_{N,t})(t,\psi^m(t,x))\big)^\trans u(t,\cdot)
 \\
&\quad +(\p_\mu f)(t,\cdot,\psi^m(t,\cdot),\mu^{X,\psi^m}_{N,t})(t,\psi^m(t,x))
 \Big],
 \nonumber
 \\
f^{m,\textrm{ex}}_N(t,x,u)
&\coloneqq 
\Big((\nabla_x \sigma)(t,x,\psi^{m}(t,x),\mu^{X,\psi^m}_{N,t})\Big)^\trans
\Big( (\nabla_x u)(t,x) \sigma(t,x,\psi^{m}(t,x),\mu^{X,\psi^m}_{N,t})\Big)
  \\
&  \quad 
+
\sE^{\mu^{X}_{N,t}}
\Big[
 \Big((\p_\mu \sigma)(t,\cdot,\psi^m(t,\cdot),\mu^{X,\psi^m}_{N,t})(t,\psi^m(t,x))\Big)^\trans 
\Big(  (\nabla_x u)(t,\cdot) \sigma(t,\cdot,\psi^{m}(t,\cdot),\mu^{X,\psi^m}_{N,t})\Big)
  \Big],
  \nonumber
 \end{align*}
 where  
$\sE^{\mu^{X}_{N,t}}[\varphi(\cdot)]\coloneqq \frac{1}{N} \sum_{l=1}^{N}\varphi(X^{l,N}_{t})$ for given function $\varphi:\sR^d\to \sR^d$,
{
and 
$f^{m,\textrm{ex}}_N$ uses 
$v(t,x)=(\nabla_x u) (t,x)\sigma^m(t,x)$.}
Note that  
if 
 $\sigma$ 
  is independent of 
the state variable and the marginal law of the state variable
(see 
Section \ref{SEC:NUM}
for concrete examples),
then
$f^{m,\textrm{ex}}_N(t,x,u)\equiv 0$
and hence 
$f^m_N$ is independent of the gradients of $u$.

Now let  
$(U_k^j)_{k\in \Z^d}$ 
be a  discrete approximation
  of  $u(t_j,\cdot)$  on the grid  $\{x_k\}_{k\in \Z^d}$.
  We approximate $f^{m,\textrm{re}}_N$ 
  by  replacing $u(t_j,\cdot)$ 
  with the monotone interpolation of $(U_k^j)_{k\in \Z^d}$:
  \begin{align}
  \begin{split}
f^{m,\textrm{re}}_{N,h}(t_j,x_k, U^j)
\label{eq:f_m_N_h_re}
&\coloneqq
 \big((\p_x b)(t_j,x_k,\psi^m(t_j,x_k),\mu^{X,\psi^m}_{N,t_j})\big)^\trans U^j_k
 +(\p_x f)(t_j,x_k,\psi^m(t_j,x_k),\mu^{X,\psi^m}_{N,t_j})
 \\
&  \quad 
+
\sE^{\mu^{X}_{N,t_j}}
\Big[
 \big((\p_\mu b)(t_j,\cdot,\psi^m(t_j,\cdot),\mu^{X,\psi^m}_{N,t})(t_j,\psi^m(t_j,x_k))\big)^\trans \textbf{i}_h[U^j](\cdot)
 \\
&\quad +(\p_\mu f)(t_j,\cdot,\psi^m(t_j,\cdot),\mu^{X,\psi^m}_{N,t_j})(t_j,\psi^m(t_j,x_k))
 \Big],
\end{split}
\end{align}
where $\textbf{i}_h$   is the   piecewise linear/multilinear interpolation operator
such that for all $\varphi=(\varphi_1, \ldots, \varphi_d):h\Z^d\to \R^d$,
$i=1,\ldots, d$, $x\in \R^d$,
$$
(\textbf{i}_h[\varphi])_i(x)=\sum_{k\in \Z^d} \varphi_i(x_k)\omega_k(x;h), \quad
$$
with the standard ``tent functions" $\{\omega_k\}_k$  satisfying $0\le \omega_k(x;h)\le 1$, $\sum_k\omega_k=1$, $\omega_k(x_j;h)=\boldsymbol{\delta}_{kj}$\footnotemark
\footnotetext{Here $\boldsymbol{\delta}_{kj}$ is the  Kronecker delta.}
 and $\textrm{supp}\, \omega_k\subset \{x\in \R^d\mid |x-x_k|\le 2h\}$.
 To approximate $f^{m,\textrm{ex}}_N$, we observe that 
 \begin{align}
\begin{split}
 f^{m,\textrm{ex}}_N(t,x,u)
=
\sum_{i,l=1}^d
\Big(
b^{m}_{N,il}(t,x) (\p_{x_l} u_i)(t,x) 
+\sE^{\mu^{X}_{N,t}}
\big[
c^{m}_{N,il}(t,x,\cdot) (\p_{x_l} u_i)(t,\cdot) 
\big]
\Big),
\end{split}
 \end{align} 
 for some functions $b^m_{N,il}:[0,T]\t \sR^d\to \sR$, 
$c^m_{N,il}(t,x,\cdot):[0,T]\t \sR^d\t \sR^d\to \sR$
 depending explicitly on 
$\sigma, \nabla_x \sigma$, $ \nabla_\mu \sigma$, $\psi^m$ and $\mu^{X,\psi^m}_{N,t}$.
Based on the signs of $b^m_{N,il}$ and $c^m_{N,il}$,
we  discretize
the terms
 $b^{m}_{N,il} (\p_{x_l} u_i)$
 and $c^{m}_{N,ij}(t,x,\cdot) (\p_{x_j} u_i)(t,\cdot)$ 
  by the  upwind finite difference schemes:
  \begin{align*}
  \mathfrak{D}^{m,j,b}_{N,h, i,l}[u(t_j,\cdot)](x)
 & =
b^{m,+}_{N,il}(t_j,x)
\frac{u_i(t_j,x+e_l h)-u_i(t_j,x)}{h}
+
b^{m,-}_{N,il}(t_j,x)
\frac{u_i(t_j,x-e_l h)-u_i(t_j,x)}{h},
\\
  \mathfrak{D}^{m,j,c}_{N,h,i,l}[u(t_j,\cdot)](x,y)
&  =
c^{m,+}_{N,il}(t_j,x,y)
\frac{u_i(t_j,y+e_l h)-u_i(t_j,y)}{h}
+
c^{m,-}_{N,il}(t_j,x,y)
\frac{u_i(t_j,y-e_l h)-u_i(t_j,y)}{h},
\end{align*}
where $\{e_l\}_{l=1}^d\subset\sR^d$ is the standard basis of $\sR^d$,
and 
$b^{\pm}=\max(\pm b,0)$  for any $b\in \sR$.
 In practice, 
  to 
evaluate 
 $ \mathfrak{D}^{m,j,c}_{N,h,i,l}[u(t_j,\cdot)](x,\cdot)$
 on the particles $(X^{l,N}_{t_j})_{l=1,\ldots, N}$,
    we shall replace 
  the grid function 
  $(U_k^j)_{k\in \Z^d}$ 
by its monotone interpolant, 
which leads to the following approximation:
 \begin{align}
\begin{split}
 f^{m,\textrm{ex}}_{N,h}(t_j,x_k,U^j)
   \label{eq:f_m_N_h_ex}
\coloneqq
\sum_{i,j=1}^d
\Big(
\mathfrak{D}^{m,j,b}_{N,h, i,l}[U^j](x_k)
+\sE^{\mu^{X}_{N,t}}
\big[
\mathfrak{D}^{m,j,c}_{N,h, i,l}\big[\textbf{i}_h[U^j]\big](x_k,\cdot)
\big]
\Big).
\end{split}
 \end{align}

Therefore, the  fully discrete scheme of \eqref{eq:PDEsystem_em} reads as: 
$U^M_k=h^m(x_k)$ for all $k\in \Z^d$,
and for all $j=1,\ldots, M, k\in \Z^d$, 
\bb\label{eq:full_discrete}
U^{j-1}_k-\Delta t {L}^{m,j-1}_{N,h}[U^{j-1}]_k=U^{j}_k+\Delta t \big(f^{m,\textrm{re}}_{N,h}(t_j,x_k,U^j)+ f^{m,\textrm{ex}}_{N,h}(t_j,x_k,U^j)),
\ee
with 
${L}^{m,j-1}_{N,h}$  defined as in 
\eqref{eq:L_m_h}, 
$f^{m,\textrm{re}}_{N,h}$ defined as in 
\eqref{eq:f_m_N_h_re}
and 
$f^{m,\textrm{ex}}_{N,h}$ defined as in 
\eqref{eq:f_m_N_h_ex}.
As the scheme adopts explicit timestepping for the gradient of $u$ {but is implicit in the second order terms}, 
we can set  $\Delta t=\cO(h)$ for numerical stability.

After obtaining the  discrete solution  $(U^j_k)_{j, k}$, 
we follow \eqref{eq:phi_psi}  
to
 update the feedback controls $\phi^{m+1}$ and $\psi^{m+1}$,
 which requires us to evaluate 
 $(\nabla F)(\psi^{m})(t_j, x_k)$  
 for all $j=0,\ldots M$ and $k\in \Z^d$.
 Observe from 
 \eqref{eq:fm} and  \eqref{eq:F_gradient_fb} 
that 
 $f^m$ and  $(\nabla F)(\psi^{m})$  
 have  similar structures, 
except from the fact  that 
 $f^m$ depends on $(\p_x H, \p_\mu H)$,
and  $(\nabla F)(\psi^{m})$ depends on $(\p_a H, \p_\nu H)$.
 Hence, 
 one can 
construct an analogue  approximation of 
$(\nabla F)(\psi^{m})(t_j, x_k)$  
by replacing 
  $(\p_x b, \p_\mu b,\p_x \sigma, \p_\mu \sigma)$ in 
   $ f^{m,\textrm{re}}_N+ f^{m,\textrm{ex}}_N$
   with $(\p_a b, \p_\nu b,\p_a \sigma, \p_\nu \sigma)$.
   This enables us to evaluate $\phi^{m+1}$ and $\psi^{m+1}$ at all grid points
  $ (t_j, x_k)$,
  and  subsequently to
obtain a particle approximation of $X^{m+1}$
for the next NAG iteration based on \eqref{eq:forward_samples} and 
the monotone interpolant  of $\psi^{m+1}$ over the grids.

\section{Implementation of the FIPDE method via residual approximation}\label{SEC:NN}
Despite  
the finite difference approximation in Section \ref{SEC:FDM}
being very effective in solving low-dimensional MFC problems,
it 
 cannot be applied to 
MFC problems with high-dimensional state processes,
due to unaffordable computational costs. 
In this section, we shall propose a neural network-based implementation of the FIPDE method, 
where 
at each NAG iteration,
we compute an approximate solution to  \eqref{PDEsystem}
by minimizing 
 a proper  residual  over a family of neural networks.

As in Section \ref{SEC:FDM}, 
we focus on the $m$-th NAG iteration
and  seek a vector-valued function $u:[0,T]\t \sR^d\to \sR^d$ satisfying   the following nonlocal parabolic PDE system:
 \begin{align}\label{eq:PDEsystem_em_nn}
(\partial_t u)(t,x)+(\mathscr{L}^m_N u)(t,x)= -f^{m}_N(t,x,u,v),
\end{align}
with $u(T,x)=h^m(x)$
and $v(t,x)=(\nabla_x u) (t,x)\sigma^m(t,x)$
 for all $(t,x)\in [0,T]\t \R^d$.
The   operator 
$\mathscr{L}^m_N$ 
and the function $f^m_N$
 (which is affine and nonlocal in $u,v$)
  are defined as in   \eqref{eq:PDEsystem_em}, 
  which depend on 
a  given feedback control   $\psi^m$
 (represented by a multilayer neural network)
and
   the empirical measures
  of a given particle approximation $(X^{l,N})_{l=1,\ldots, N}$
 of the state process  $X^m$
 for the present NAG iteration.

In the following, we shall reformulate \eqref{eq:PDEsystem_em_nn} 
into an empirical risk minimization problem over multilayer neural networks,
which is then solved by using stochastic gradient descent (SGD) algorithms;
see e.g., the Deep Galerkin Method (DGM) in \cite{sirignano2018dgm}.
However, instead of directly applying DGM to \eqref{eq:PDEsystem_em_nn}, 
we shall consider an acceleration method by first decomposing the solution $u$ into:
\bb\label{eq:decomposition}
u(t,x)=\bar{u}(t,x)+\tilde{u}(t,x), 
\quad (t,x)\in [0,T]\t \sR^d,
\ee 
where 
$\bar{u}$ is an approximate solution to \eqref{eq:PDEsystem_em_nn}
computed by some efficient numerical methods,
and  
$\tilde{u}$ is a  residual correction of $\bar{u}$ based on neural networks. 
By computing   the dominant part of the solution $u$    efficiently 
and merely applying  SGD  algorithms for
 the  small residual term,
we can  
obtain an accurate and stable approximation of $u$
with a small number of SGD iterations,
and subsequently 
reduce the total computation time  for solving  PDE systems at all NAG iterations.

In general,
the numerical solver of $\bar{u}$ should be designed in a problem  dependent way. 
For many practical MFC problems
(see, e.g.,~\cite{gobet2019extended,albi2021moment}), 
we can first  linearize the dynamics  \eqref{cforward} around the  target states 
and approximate the cost functions \eqref{cvalue} by suitable quadratic costs.
This leads to a linear-quadratic (LQ) approximation of the MFC problem \eqref{cforward}-\eqref{cvalue},
and the  approximate solution  $\bar{u}$ of \eqref{eq:PDEsystem_em_nn} can be chosen 
as the solution of the  (matrix-valued) differential  Riccati equations
for the resulting LQ MFC problem.
Consequently, the decomposition \eqref{eq:decomposition}
can be viewed as a neural network-based  nonlinear correction to the (suboptimal) linear feedback control. 
We refer the reader to Section \ref{sec:CS_6d}
for more details on  the LQ approximation of   MFC problems with nonsmooth costs.

Given the approximate solution $\bar{u}:[0,T]\t \sR^d\to \sR^d$, we see  from \eqref{eq:PDEsystem_em_nn}
that  the residual term $\tilde{u}:[0,T]\t \sR^d\to \sR^d$ satisfies 
 for all $(t,x)\in [0,T]\t \R^d$,
 \begin{align}\label{eq:PDEsystem_residual}
(\partial_t \tilde{u})(t,x)+(\mathscr{L}^m_N \tilde{u})(t,x)= -\tilde{f}^{m}_N(t,x,\tilde{u},\tilde{v}),
\quad \tilde{u}(T,x)=h^m(x)-\bar{u}(T,x),
\end{align}
where 
$\tilde{v}(t,x)=(\nabla_x \tilde{u}) (t,x)\sigma^m(t,x)$, and 
$$
\tilde{f}^{m}_N(t,x,\tilde{u},\tilde{v})
=f^{m}_N(t,x,\tilde{u}+\bar{u},\tilde{u}+\bar{v})
+(\partial_t \bar{u})(t,x)+(\mathscr{L}^m_N \bar{u})(t,x).
$$
We then extend the residual based method for scalar PDEs
in \cite{sirignano2018dgm,ito2021neural}
 to  the coupled PDE system   \eqref{eq:PDEsystem_residual}.
 In particular, let 
 $\cD\subset \sR^d$ be the chosen computational domain and 
 $\cN_u=\{{u}^\theta:[0,T]\t \sR^d\to \sR^d\mid \theta\in \sR^p\}$
 be a family of multilayer neural networks with some prescribed  architectures
 and sufficiently smooth activation functions. Then
 we 
 seek the optimal neural network in $\cN_u$ to approximate $\tilde{u}$
by minimising the following  loss function over the weights   $\theta$ :
\begin{align}\label{eq:cE}
\begin{split}
\cE(\theta)
&=
\|(\partial_t \tilde{u}^\theta)(\cdot,\cdot)+(\mathscr{L}^m_N \tilde{u}^\theta)(\cdot,\cdot)+\tilde{f}^{m}_N(\cdot,\cdot,\tilde{u}^\theta,(\nabla_x \tilde{u}^\theta) \sigma^m)\|^2_{[0,T]\t \cD, \nu_1}
\\
&\quad 
+\eta_1\| \tilde{u}^\theta(T,\cdot) -(h^m(\cdot)-\bar{u}(T,\cdot))\|^2_{\cD, \nu_2}
+\eta_2\| \tilde{u}^\theta(\cdot,\cdot) \|^2_{[0,T]\t \p\cD,  \nu_3}.
\end{split}
\end{align}
Here, $\nu_i$, $i=1,2,3$, are some given probability measures on $[0,T]\t \cD $, $\cD$ and $[0,T]\t \p\cD$ respectively, 
and $\eta_1,\eta_2> 0$ are 
  some given weighting parameters  (possibly different among all NAG iterations)
  introduced to balance  
the interior residual  and the residuals of the boundary data.
Note that $\tilde{u}^\theta$ takes values in $\sR^d$, and  \eqref{eq:cE} contains residuals of all components. 

In practice, 
the loss function \eqref{eq:cE}
can be minimized by using SGD based on a sequence of  
mini-batches of 
pseudorandom points or quasi-Monte Carlo points.
More precisely, for the $j$-th SGD iteration with $j\in \N$,
we first  generate $N_{\textrm{in}}$ points from $[0,T]\t \cD$,
$N_{\textrm{ter}}$ points from $\cD$
and $N_{\textrm{bdy}}$ points from $[0,T]\t\p \cD$ 
according to the measures $\nu_1$, $\nu_2$ and $\nu_3$, 
respectively, 
then evaluate the following empirical loss with the current weights $\theta_j\in \sR^k$: 
\begin{align}\label{eq:E_empirical}
\begin{split}
\cE_{\textrm{em}}(\theta_j)
&=
\frac{1}{N_{\textrm{in}}}\sum_{i=1}^{N_{\textrm{in}}}
|(\partial_t \tilde{u}^{\theta_j})(t_i,x_i)+(\mathscr{L}^m_N \tilde{u}^{\theta_j})(t_i,x_i)+\tilde{f}^{m}_N(t_i,x_i,\tilde{u}^{\theta_j},(\nabla_x \tilde{u}^{\theta_j}) \sigma^m)|^2
\\
&\quad 
+
\frac{\eta_1}{N_{\textrm{ter}}}\sum_{i=1}^{N_{\textrm{ter}}}
| \tilde{u}^{\theta_j}(T,x_i) -(h^m(x_i)-\bar{u}(T,x_i))|^2
+
\frac{\eta_2}{N_{\textrm{bdy}}}\sum_{i=1}^{N_{\textrm{bdy}}}
| \tilde{u}^{\theta_j}(t_i,x_i) |^2,
\end{split}
\end{align}
and finally obtain the updated weights $\theta_{j+1}=\theta_{j+1}-\tau_j (\nabla_\theta\cE_{\textrm{em}})(\theta_j)$ with a stepsize $\tau_j>0$.
The above  SGD iterations are performed until 
the  index $j$ meets the maximum iteration number
or 
an accuracy tolerance is satisfied. 

Finally, we discuss the implementation of Step 5 in Algorithm \ref{alg:NAG2}
based on
 an approximate solution  $\tilde{u}^\theta$ 
 to \eqref{eq:PDEsystem_residual}
 and the current  feedback control $\psi^m$, where
 both are  
  represented by a multilayer neural network. 
By  \eqref{eq:decomposition},
  the solution $u^m$ to \eqref{eq:PDEsystem_em_nn}
  is now approximated by 
  ${u}^\theta\coloneqq \bar{u}+\tilde{u}^\theta$, 
  which leads to a pointwise approximation 
$(\nabla F)(\psi^{m},\theta)$  
of   the function 
$(\nabla F)(\psi^{m})$
by replacing $u^m$ in \eqref{eq:F_gradient_fb}
with $u^\theta$.
Note that the required derivative $\nabla _x {u}^\theta$ for approximating the function $v^m$ 
can be computed analytically if  $\bar{u}$ and $\tilde{u}^\theta$ are differentiable.
However,  in contrast to the grid-based representation of feedback controls 
in Section \ref{SEC:FDM},
the  neural network representation of the feedback control
prevents us from obtaining the updated control $\psi^{m+1}$ 
in 
\eqref{eq:phi_psi}
via simple operations applied to the parameters.
This is due to the fact that 
a linear combination of  multilayer neural networks 
with nonlinear activation functions
in general can only be expressed as a neural network with more complicated architectures
(see e.g., \cite[Lemma A.1]{reisinger2020rectified}).
Hence, exactly following the update rules in \eqref{eq:phi_psi}
would require us to save all  $(\psi^m)_{m\in \N\cup\{0\}}$ and $(u^m)_{m\in \N\cup\{0\}}$,
which increases both the memory requirements and  the computational cost 
for Steps 3 and 4 of Algorithm \ref{alg:NAG2}
if the NAG iteration $m$ is large. 
To overcome this difficulty, 
we shall 
follow    \eqref{eq:phi_psi}  approximately 
and obtain the updated controls 
as follows:
\begin{align}\label{eq:supervised}
\begin{split}
 \phi^{m+1}
 &\in \argmin_{\phi\in\cN_\phi} 
 \|\phi -\big(\text{prox}_{\tau\ell} (\psi^m -\tau 
 (\nabla F)(\psi^{m},\theta)  )\big)\|^2_{[0,T]\t \cD, \nu_1},
 \\
 \psi^{m+1}
 &\in \argmin_{\psi\in\cN_\psi} 
 \|\psi-\big(\phi^{m+1}  + \tfrac{m}{m+3}(\phi^{m+1} -\phi^m)\big)
 \|^2_{[0,T]\t \cD, \nu_1},
\end{split}
\end{align}
where $\cN_\phi$ and $\cN_\psi$
are families of multilayer neural networks with prescribed architectures and  unknown parameters. 
These  supervised learning  problems can  be easily solved
by performing gradient descent 
based on  sample  points from $[0,T]\t \cD$.

\section{Numerical experiments}\label{SEC:NUM}

In this section, we demonstrate the effectiveness of the FIPDE scheme through numerical experiments. We present two MFC problems with nonsmooth optimization objectives:
a portfolio liquidation problem with trade crowding and  transaction  costs
in Section \ref{SEC:PF},
and 
sparse consensus control of two-dimensional and six-dimensional stochastic Cucker--Smale models
in Section \ref{SEC:CS}.
For both examples, 
we  benchmark  the FIPDE method 
with existing pure data-driven algorithms,
including
the empirical regression method 
in \cite{lemor2006rate}
and the neural network-based policy gradient method
in \cite{carmona2019convergence},
which  lead to approximate 
 feedback control
merely along  trajectories of the optimal state process.
Our experiments show that
compared with  pure data-driven approaches,
 the  FIPDE method leads to a global approximation of the optimal feedback control,
and achieves a robust performance in terms of model perturbations.

\subsection{Portfolio liquidation with trade crowding and  transaction costs}\label{SEC:PF}
In this section, we consider 
a portfolio liquidation problem  
where a large number  of market participants try to liquidate
their positions on the same asset
by a given
 terminal time
$T>0$ (see \cite{basei2017linear,acciaio2019extended}),
while 
taking  into account 
  (possibly nonsmooth) execution costs
  and
the permanent price impact caused by their trading actions.
The cooperative equilibrium leads to a MFC problem
 for 
a representative agent.

Let $(\alpha_t)_{t\in [0,T]}$
be the  trading speed
chosen by the representative agent, 
the state dynamics of the MFC problem 
is given by:
for all $t\in [0,T]$,
\begin{align*}
\d Q_t = \alpha_t \, \mathrm{d}t, \quad
\mathrm{d}S_t &= \lambda \mathbb{E}[\alpha_t] \, \mathrm{d}t + \sigma \, \mathrm{d}W_t, 
\end{align*}
where $Q = (Q_t)_{t \in [0,T]}$ is the inventory process
with  a \textit{random}  initial state $Q_0$ representing
the initial inventories for all participants,
$S = (S_t)_{t \in [0,T]}$ is the  price process
with a deterministic initial state 
$S_0 = s_0\ge 0$,
and 
 $\lambda \mathbb{E}[\alpha_t]$ with $\lambda\ge 0$ 
represents the permanent market impact
on the asset price
 due to 
the  trading of all 
  participants. Here, the processes $Q$, $S$ and $W$ are all one-dimensional.
  
 The objective of the agent  is then to
minimize the following cost functional
\bb \label{eq:price_impact}
J(\alpha;(S_0,Q_0))=\mathbb{E}\left[ 
 \int_{0}^{T} \Big( \alpha_t S^\a_t + 
 (Q^\a_t)^2 +
 k_1|\a_t |^2+k_2|\a_t|
 \Big) \, \mathrm{d}t - Q^\alpha_T(S^\a_T-\gamma Q^\a_T)
 \right],
\ee
over all possible  trading speeds $\a$, 
where 
$(Q^\a_t)^2$ penalizes the current inventory,
  $k_1,k_2\ge 0$,
  and $Q^\a_T (S^\a_T- \gamma  Q^\a_T )$ with $\gamma\ge 0$ 
 is the liquidation value of the remaining inventory at terminal time.
Note that    the nonsmooth term $k_2|\a|$ models 
proportional execution costs such as the bid-ask spread, the fees paid to the venue,
and/or a stamp duty
(see, e.g., \cite{gueant2016financial}).
 In the following, we shall perform experiments with 
different choices of $Q_0$ and $k_2$ 
(whose values will be specified later)
while fixing the other 
model  parameters as:
$T=1$, 
$s_0=2$, $\lambda=0.5$, $\sigma=0.7$, $\gamma=0.5$
and  $k_1=1$.

We  
initialize
Algorithm \ref{alg:NAG2}
with stepsize $\tau=1/6$
and  initial guess $\phi^0= 0$.
At the $m$-th NAG iteration, 
given an approximate control strategy 
$\psi^m:[0,T]\times \R^2\to \R$
with associated state  processes $(S^{m},Q^m)$,
we consider the following decoupled system of PDEs
(cf.,~\eqref{PDEsystem}):
for all $(t,s,q)\in [0,T)\t \sR\t \sR$,
\begin{subequations}\label{eq:pde_price}
\begin{alignat}{2}
(\p_tu_1)(t,s,q)+( \mathscr{L}^m u_1)(t,s,q)
&=
- \psi^m(t,s,q),
&&
\quad  
u_1(T,s,q)=-q, 
\\
(\p_tu_2)(t,s,q)+( \mathscr{L}^m u_2)(t,s,q)
&=
-2q,
&&
\quad 
u_2(T,s,q)=-s + 2\gamma q, 
\end{alignat}
\end{subequations}
with 
$\mathscr{L}^m$
such that 
 for each $\varphi\in C^{1,2}([0,T]\times \R^2;\R)$ and $(t,s,q)\in [0,T]\times \R\t \R$,
 \bb
 \label{eq:L_portfolio}
( \mathscr{L}^m \varphi)(t,s,q)
=
\tfrac{1}{2}\sigma^2(\p_{ss}\varphi)(t,s,q)
+\lambda \mathbb{E}[\psi^m(t,S^m_t,Q^m_t)](\p_{s}\varphi)(t,s,q)
 + \psi^m(t,s,q)(\p_{q}\varphi)(t,s,q).
 \ee
 
 We implement 
 the
FIPDE algorithm 
with the monotone scheme 
\eqref{eq:full_discrete}
to solve the system \eqref{eq:pde_price}.
 We first
approximate the expectations in 
\eqref{eq:L_portfolio}
by 
empirical averages over
$N=10^4$ particle approximations
of $(S^{m},Q^m)$
generated by \eqref{eq:forward_samples}
with $\Delta t=1/50$,
localize the equation on the  domain  $ \cD=[-2,6]\times [0,4]$
and impose    boundary conditions as the terminal condition,
i.e., $(u_1(t,s,q),u_2(t,s,q))=(-q, -s + 2\gamma q)$
for all $(t,s,q)\in [0,T]\times \p\cD$.
Our experiments with larger computational domains
indicate that 
this domain truncation and boundary condition
 leads to a  negligible  domain truncation error.
Then we 
 construct 
an implicit first-order monotone scheme
\eqref{eq:full_discrete}
for the localized  system \eqref{eq:pde_price}
by discretizing
the first-order
and second-order 
 derivatives in \eqref{eq:L_portfolio}
via   the  upwind finite difference
and
 the central difference,
 respectively,
 The chosen time stepsize is $\Delta t=1/50$,
and the spatial mesh sizes $h_s=8/50$, $h_q=4/50$. 

As a benchmark for the  FIPDE scheme, we 
also implement an empirical regression (EMReg) method 
to solve \eqref{eq:price_impact},
where for each NAG iteration, we solve 
  the adjoint BSDE \eqref{adjointa} by projecting the decoupling fields $(u^m,v^m)$ of $(Y^m,Z^m)$
on prescribed vector spaces of  basis functions, 
and evaluating
the coefficients  by performing regressions based on  simulated trajectories of  $(S^{m},Q^m)$
(see, e.g.,~\cite{lemor2006rate}).
In particular, 
we  partition 
 the computation domain    $ \cD=[-2,6]\times [0,4]$
with the same meshsize 
$h_s=8/50$ and $h_q=4/50$ as that of the FIPDE method,
and choose the indicator functions of all subcells as the basis functions in the regression. 
We also employ 
\eqref{eq:forward_samples} with the same parameters $N=10^4$ and $M=50$
as those for the FIPDE method
to generate state trajectories, in order to ensure a fair comparison.
{
Other implementation details are given in Appendix \ref{appendix:portfolio}.
}

\begin{figure}[!ht]
    \centering
    
       \begin{subfigure}{.31\textwidth}
        \centering
        \includegraphics[width=\linewidth,height=4.4cm]{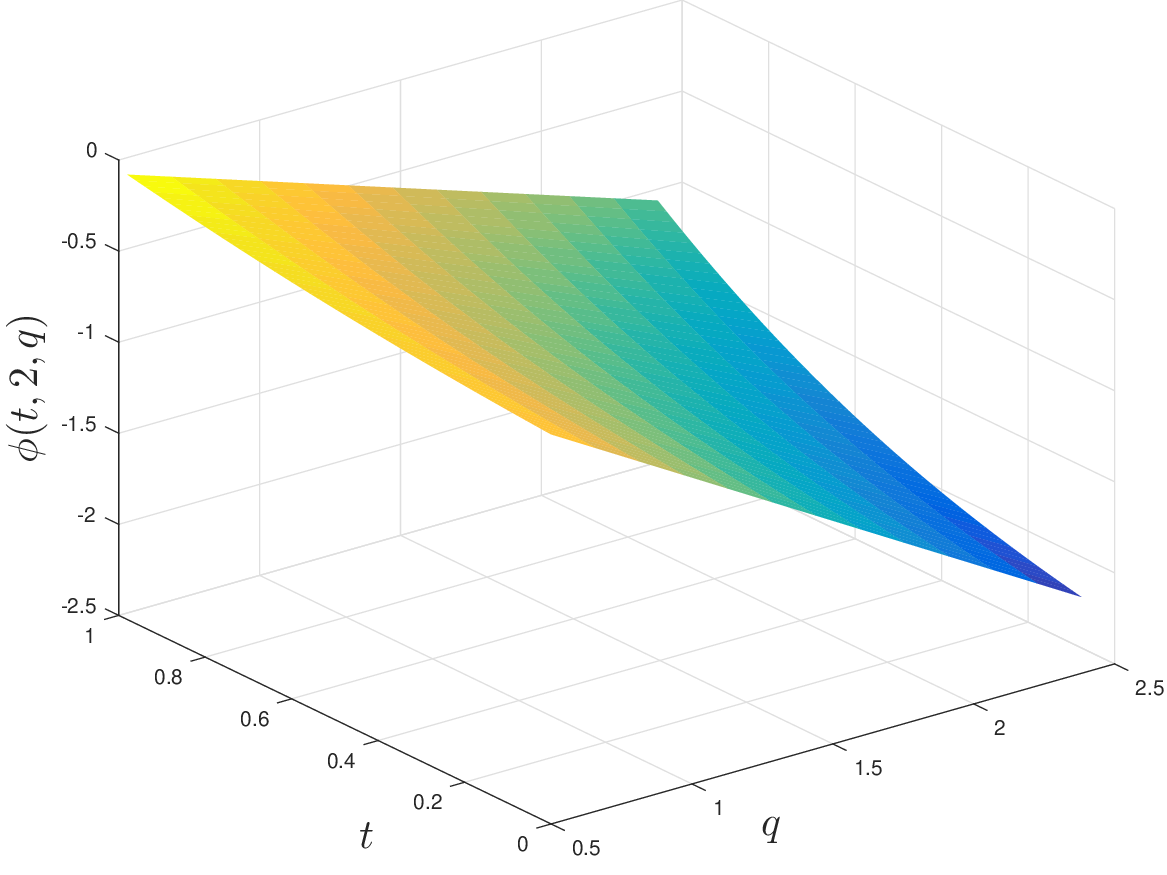}
        \caption{Analytic solution}
    \end{subfigure}\quad
    \begin{subfigure}{.31\textwidth}
        \centering
         \includegraphics[width=\linewidth,height=4.4cm]{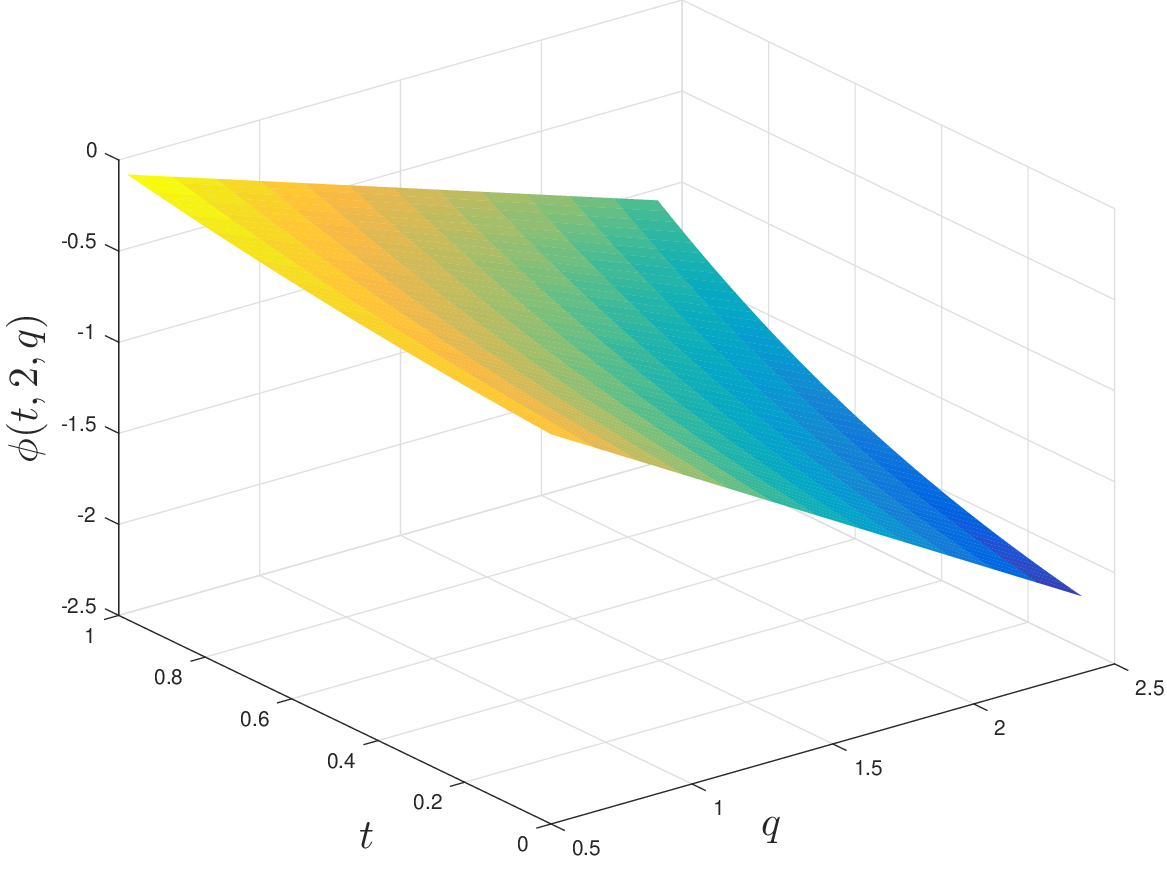}
        \caption{The FIPDE scheme}
    \end{subfigure}\quad
    \begin{subfigure}{.31\textwidth}
        \centering
         \includegraphics[width=\linewidth,height=4.4cm]{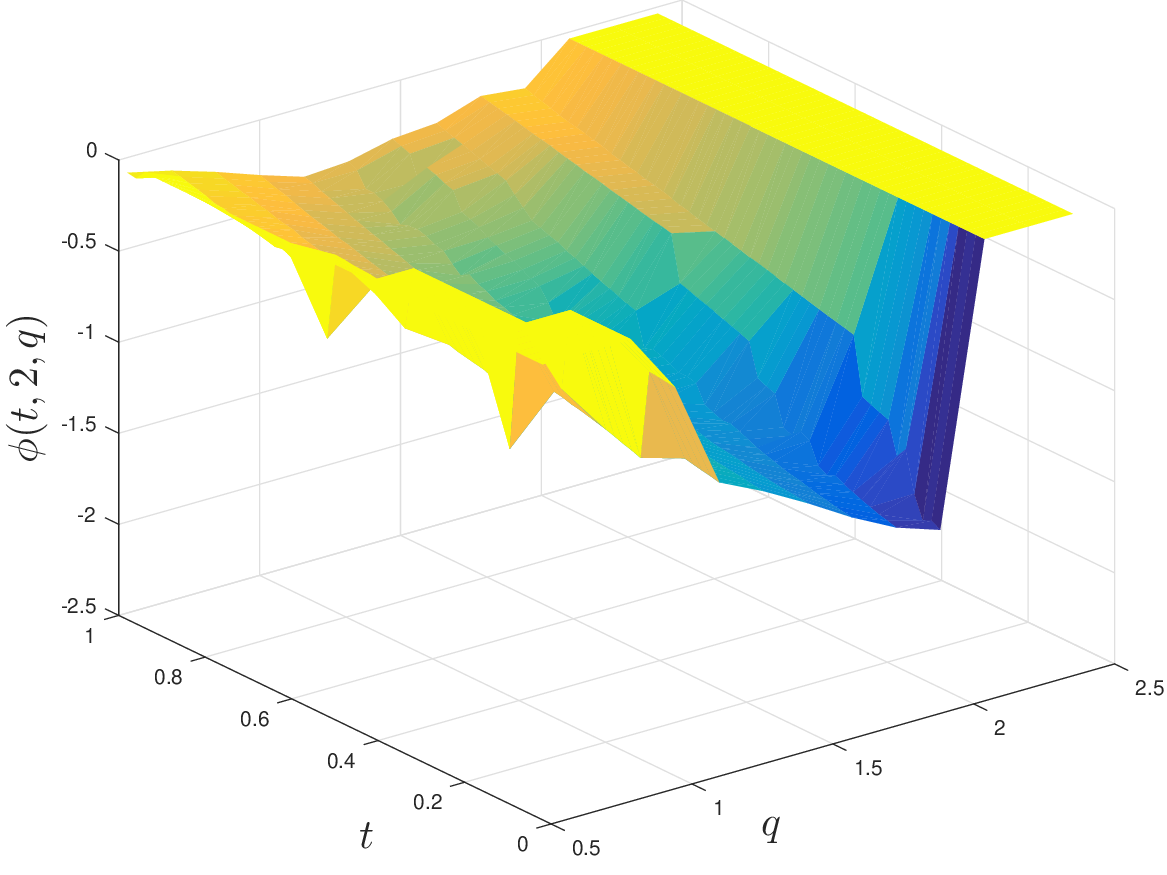}
        \caption{The EMReg  scheme}
        \end{subfigure}
 
    \caption{
    Feedback controls 
    of the LQ optimal liquidation problem with 
    $Q_0\sim\mathcal{U}(1,2)$
    obtained by using different methods. }
    \label{fig:linear_control}
 \end{figure}

We first examine the performance of the FIPDE 
and  EMReg schemes
for solving the linear-quadratic (LQ)
MFC problem \eqref{eq:price_impact}
with 
$Q_0\sim \mathcal{U}(1,2)$
(the uniform distribution on $(1,2)$) and
 $k_2=0$.
Extending the arguments  in  \cite{yong2013linear}
to \eqref{eq:price_impact}
yields  that 
for any square-integrable initial condition $Q_0$,
the optimal feedback control of \eqref{eq:price_impact} is of the form:
\bb\label{eq:optimal_fb_price}
\phi^\star(t,s,q)=a_t q+b_t\sE[Q^\star_t],
\quad \forall (t,s,q)\in [0,T]\t \sR\t \sR,
\ee
where $Q^\star$ is the optimal inventory process,
and $a_t,b_t:[0,T]\to \sR$ are solutions to some Riccati equations
depending explicitly on $T,\lambda, \gamma$ and $k_1$, but independent of $s_0$ and $Q_0$.
Despite the fact that
both  the FIPDE and EMReg schemes 
 achieve 
less than  $1\%$ absolute error
for value function approximations
within 5  NAG iterations
  (see Figure \ref{fig:price_value_conv} in Appendix \ref{appendix:portfolio}
 for more details),
  these two methods generate 
 qualitatively different feedback controls.
 Figure \ref{fig:linear_control} compares  the exact feedback control \eqref{eq:optimal_fb_price}
 and  the approximate  feedback controls
obtained by the FIPDE  and  EMReg schemes after 20 NAG iterations,
for which we evaluate  the feedback strategies 
at 
$s=2$ and  $(t,q)\in [0,1]\t [0.5,2.5]$.
One can clearly observe that 
 the approximate feedback control from the FIPDE scheme is in almost exact agreement with the analytic solution
  on the entire computational domain.
  In contrast,
  the EMReg method produces a much more irregular control strategy,
which is only accurate on a certain part  of the domain.
Recall that 
we initialize the NAG iteration in Algorithm \ref{alg:NAG2} with  $\phi^0=0$,
and Figure \ref{fig:linear_control} (right) indicates that   the EMReg method does not update the initial guess 
for $q>2$.
In fact, 
as  the EMReg  method
approximates the adjoint processes
by performing regression 
based on the simulated  trajectories of  
the   state process,
 the EMReg  method
 not only suffers more from 
statistical errors,
but also
cannot recover  the exact feedback control
beyond the support of the simulated trajectories. 
As we shall see soon, 
such a local approximation property makes these pure data-driven approaches unstable with respect to model perturbation.

We proceed to investigate the robustness of 
feedback controls obtained by the FIPDE scheme 
and the EMReg scheme
with respect to perturbation of the random initial condition $Q_0$,
which  models the initial inventories 
of all market participants.
As it is practically difficult to obtain the precise distribution of  $Q_0$,
it is desirable for  the feedback control to be stable in terms of 
the uncertainty in the initial law (see, e.g.,~\cite{germain2021deepsets}).
In particular, we first 
obtain
 approximate feedback controls
by 
applying
the FIPDE  and EMReg schemes 
to solve  \eqref{eq:price_impact}
with $Q^\textrm{pre}_0\sim \mathcal{U}(1,2)$, 
and then 
examine the performance of 
these feedback controls
 on perturbed models with 
$Q_0\sim \mathcal{U}(Q_{0,\min},Q_{0,\max})$ for   $Q_{0,\min}\in [0.5,1.5]$ and $Q_{0,\max}\in [1.5, 2.5]$.

\begin{figure}[!ht]
    \centering
     \includegraphics[  width=0.35\columnwidth,height=4.8cm]{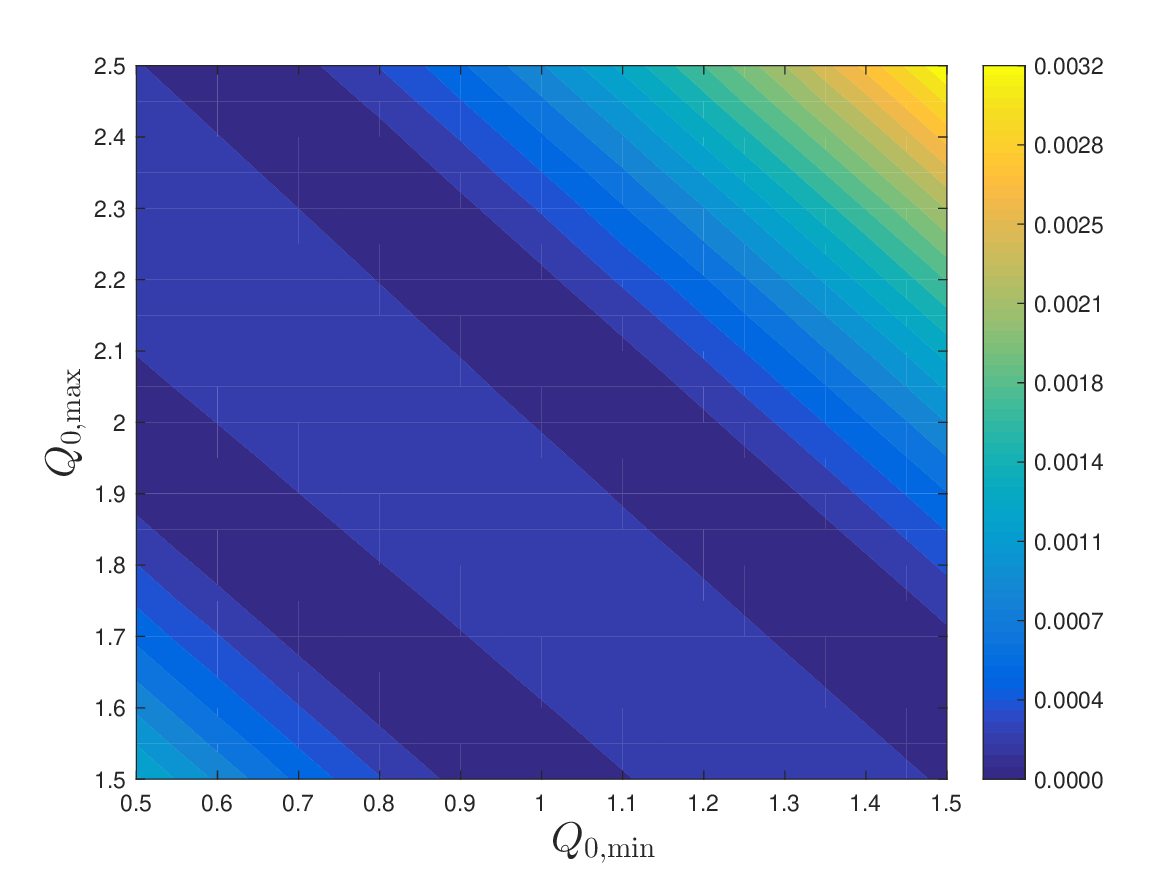}
       \hspace{1cm}
     \includegraphics[  width=0.35\columnwidth,height=4.8cm]{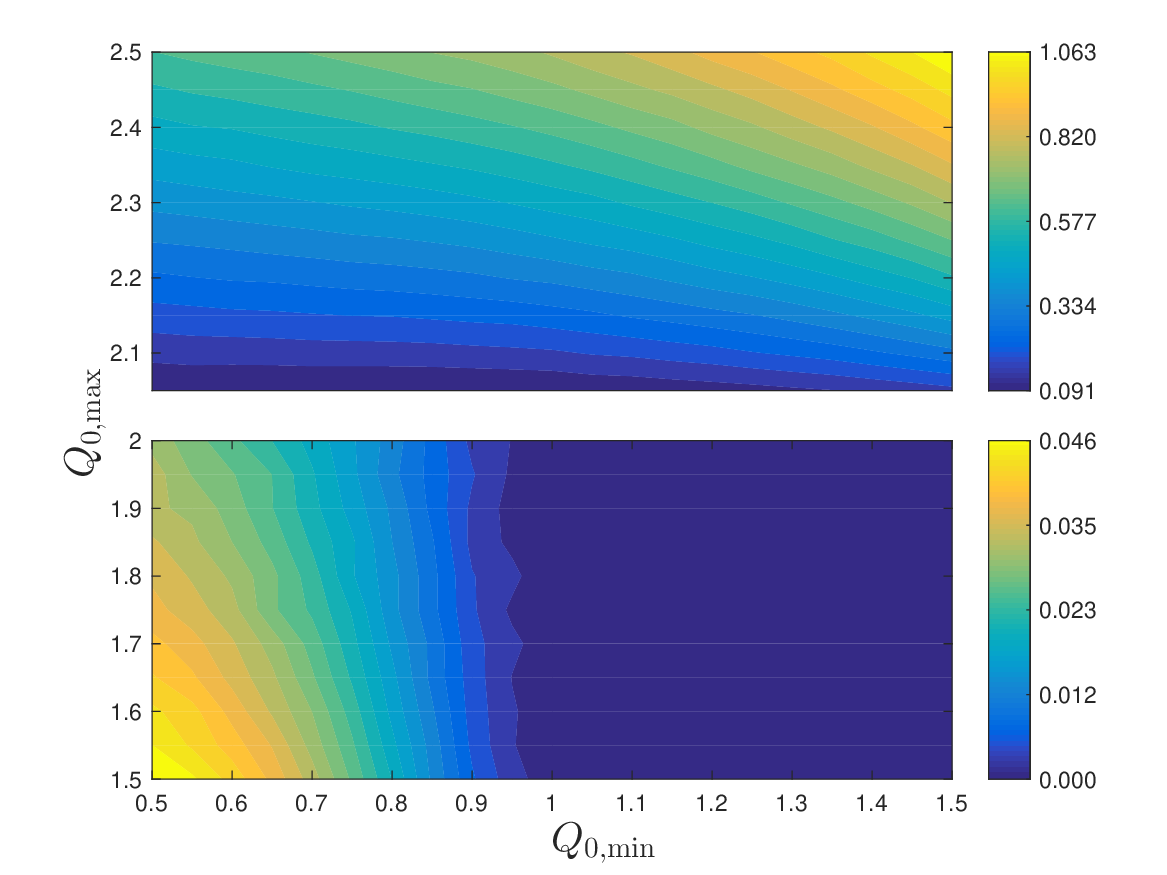} 
     
       \caption{
       Absolute performance gaps of  
       precomputed feedback controls  
from
the FIPDE scheme (left) and 
       the EMReg scheme (right)
       on perturbed models with   $Q_0\sim \mathcal{U}(Q_{0,\min},Q_{0,\max})$.
    }
    \label{fig:robustness_linear_abs}
 \end{figure}

Figure \ref{fig:robustness_linear_abs}  illustrates 
 the performance   of these  precomputed controls on models with different  $Q_0$ 
 in terms of the absolute performance gap
$|J(\phi^\textrm{pre})-J^\star|$, 
where  
$J(\phi^\textrm{pre})$  
is the expected cost 
 of 
 a precomputed  feedback control $\phi^\textrm{pre}$
on  the perturbed model,
 and 
$J^\star$
is the optimal cost of the perturbed model.  
One can clearly observe from Figure \ref{fig:robustness_linear_abs} (left)
that 
the precomputed feedback control
from 
the FIPDE scheme 
is  very robust with respect to perturbation of $Q_0$,
as it yields  extremely small   performance gaps
for all perturbed  models.
Moreover,
the absolute performance gaps 
remain almost constants
 for different $Q_0$ with the same expectation.
This is because
the  feedback control from
the FIPDE scheme 
 captures  the precise spatial 
dependence of the optimal feedback control
of a perturbed problem
(i.e., the function $\phi^\star$ in \eqref{eq:optimal_fb_price}),
while it keeps the measure   dependence 
the same
as that for the unperturbed model. 
As 
the optimal feedback control
of the LQ MFC problem \eqref{eq:price_impact}
 depends  on the 
law of  $Q_0$ only 
 through its first moment,
 the absolute performance gap 
 $|J(\phi^\textrm{pre})-J^\star|$
 for the FIPDE scheme is purely determined by 
the perturbation  of $\sE[Q_0]$.

In contrast, 
as shown in Figure \ref{fig:robustness_linear_abs} (right),
the approximate feedback control from the EMReg scheme 
is very sensitive to  perturbations of $Q_0$,
where
  absolute performance gaps 
are typically a few magnitudes larger than
those of the FIPDE scheme.
This
phenomenon is more pronounced
if the support of the perturbed $Q_0$ 
is not contained by   the interval $[1,2]$,
i.e., the support of the original $Q_0$.
This is because
the  EMReg scheme
(and other pure data-driven algorithms)
computes  approximate feedback controls
 merely based on trajectories of the original state process,
and consequently the resulting   feedback control 
performs poorly along the trajectories 
outside the support of the original system.
In particular,
since 
it is critical to execute  
a proper  strategy for large values of $Q_t$
in the present liquidation problem,
  a slight perturbation of $Q_{0,\max}$ will significantly worsen the  performance of 
the precomputed control.

\begin{figure}[!ht]
    \centering
     \includegraphics[  width=0.35\columnwidth,height=4.8cm]{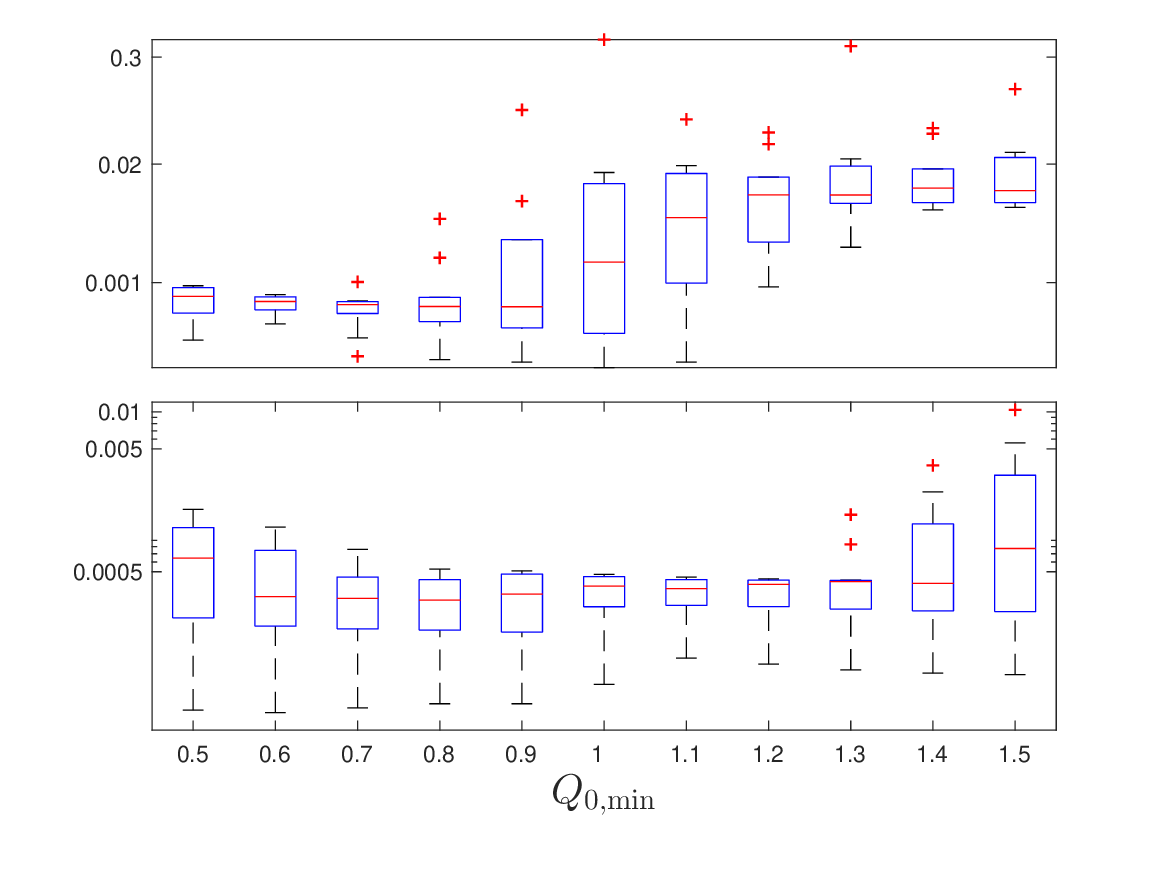}
       \hspace{1cm}
     \includegraphics[  width=0.35\columnwidth,height=4.8cm]{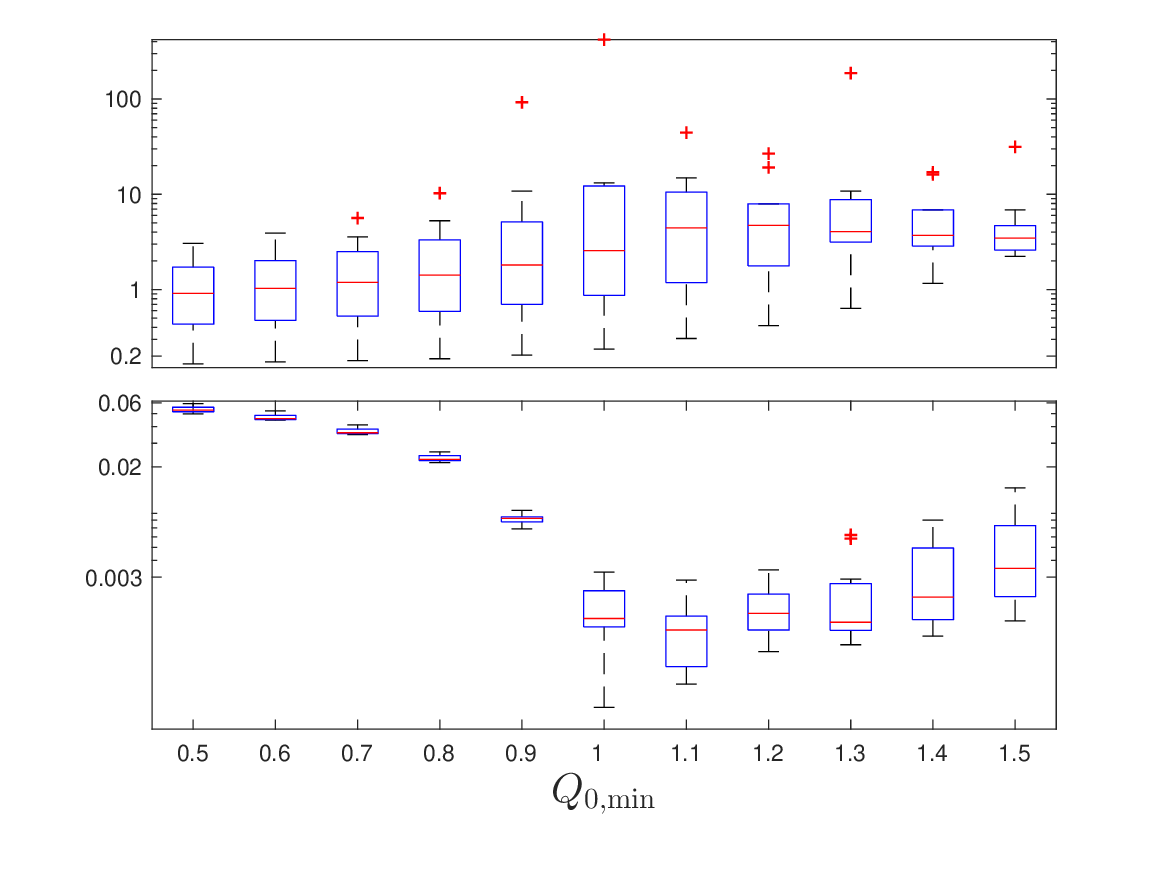} 
     
       \caption{
Relative performance gaps of 
	precomputed feedback controls  
       on perturbed models       with   $Q_0\sim \mathcal{U}(Q_{0,\min},Q_{0,\max})$;
       from left to right:  results for
       the FIPDE and EMReg scheme;
	from top  to bottom: 
	   	distributions of
	relative errors for fixed $Q_{0,\min}$
	and varying   $Q_{0,\max}\in [2,2.5]$ 
	and $Q_{0,\max}\in [1.5,2]$,
	where on each box, 
the central line 
 is the median,
 the edges of the box are 
the 25th and 75th percentiles,
the  whiskers extend  to  non-outliers extreme data points, 
and  outliers are plotted individually.  
    }
    \label{fig:robustness_linear_rel}
 \end{figure}

The improved robustness of the FIPDE scheme over the EMReg scheme 
can be better analyzed by the 
relative performance gaps 
$|J(\phi^\textrm{pre})-J^\star|/|J^\star|$,
whose distributions 
(with  fixed $Q_{0,\min}$ and varying $Q_{0,\max}$)
are summarized by the box plots
in Figure \ref{fig:robustness_linear_rel}.
We can clearly see that 
the precomputed control of the FIPDE scheme 
achieves a relative error of less than 2\%  
on most perturbed models
(Figure \ref{fig:robustness_linear_abs} (left)), 
while 
 the precomputed control 
of the EMReg scheme
 will typically lead to  
a relative error ranging from $20\%$ to $10^3\,\%$
(Figure \ref{fig:robustness_linear_rel} (top-right)).\footnotemark
\footnotetext{
We have ignored 
 the
 outliers  (marked as plus signs)
in  Figure \ref{fig:robustness_linear_rel} (top),
which resulted from evaluating  
relative errors
with   optimal costs   $J^\star$ 
very close to zero. 
}

Finally, we  examine  the performance of the FIPDE scheme for solving 
the nonsmooth  MFC  problem
\eqref{eq:price_impact}
with
 $k_1=k_2=1$ in \eqref{eq:price_impact}.
 We 
 initialize Algorithm \ref{alg:NAG2}
 with  $\tau=1/6$ and $\phi^0=0$,
 and for each NAG iteration,  discretize 
the PDE system \eqref{eq:pde_price}
by using the same monotone scheme as that for the LQ case.
Note that for all $\tau>0$, the proximal operator of $\tau |\cdot|$ is given by
$\text{prox}_{\tau |\cdot|}(x)= \textrm{sgn}(x)(|x|-k_2\tau)$ for all $|x| \geq k_2\tau$
and $\text{prox}_{\tau \ell}(x)=0$ for all $|x| \leq k_2\tau$.

\begin{figure}[ht]
    \centering
         \includegraphics[  width=0.3\columnwidth,height=4.5cm]{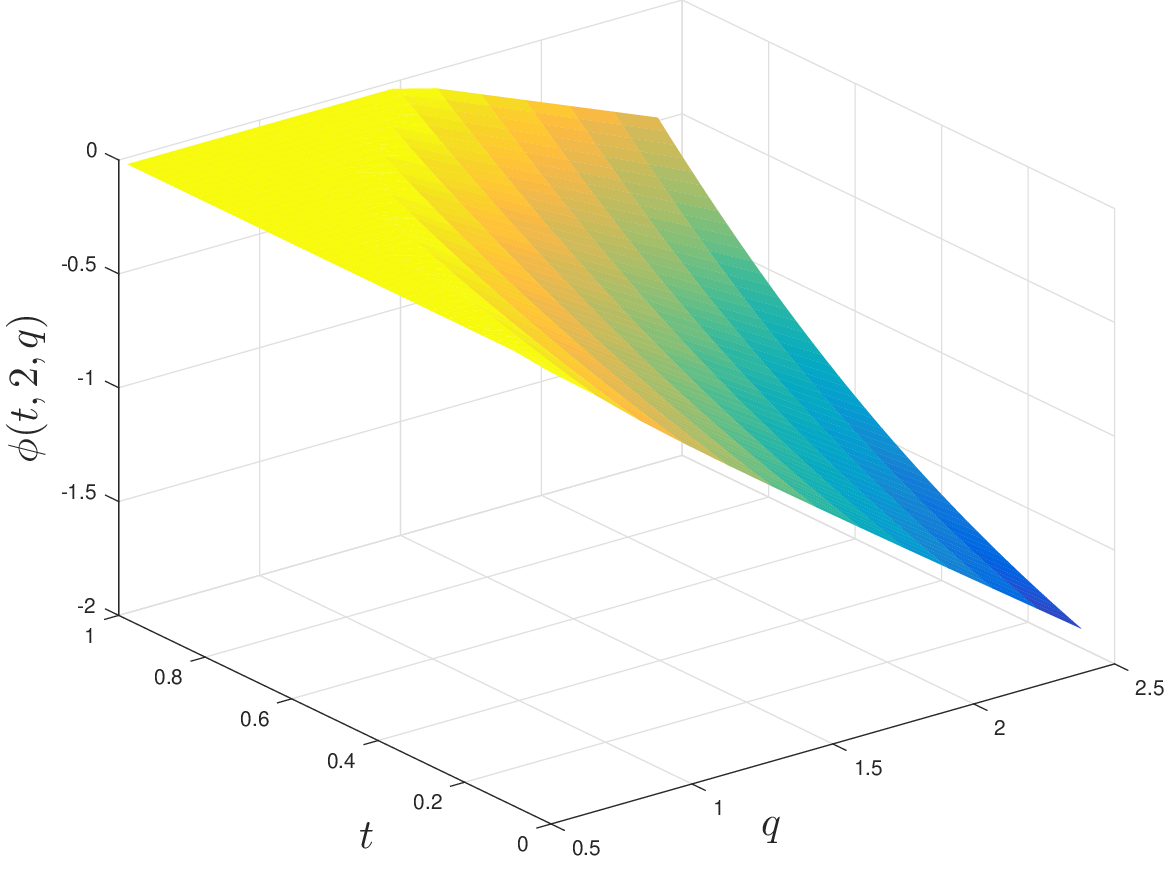}
  \hspace{1.5cm}
       \includegraphics[  trim=0 5 0 20, clip, width=0.35\columnwidth,height=4.5cm]{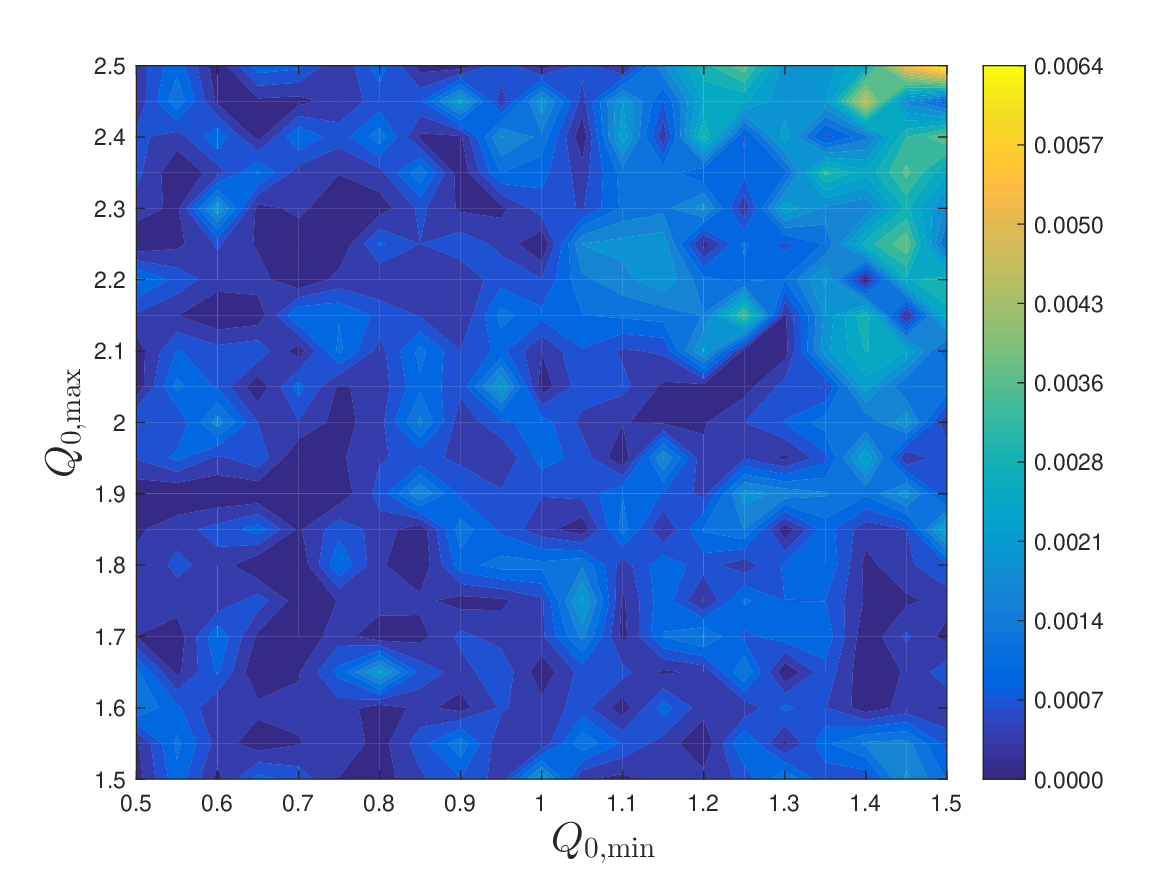}

    \caption{
  Numerical results 
  of the FIPDE scheme 
  for the nonsmooth optimal liquidation problem with 
$k_2=1$; 
from left to right: 
approximate feedback control
 for $Q_0\sim \mathcal{U}(1,2)$,
 and 
absolute performance gaps of the precomputed  feedback  control
  on perturbed models with    $Q_0\sim \mathcal{U}(Q_{0,\min},Q_{0,\max})$.
 }
    \label{fig:nonlinear_price}
 \end{figure}

We carry out  the FIPDE scheme with 20 NAG iterations,
whose convergence is shown in Figure \ref{fig:price_value_conv}
in Appendix \ref{appendix:portfolio};
note that in fact 4 NAG iterations are sufficient to approximate the value function accurately.
The resulting  approximate feedback control 
is  independent of the  variable $s$
but nonlinear in the  variable $q$.
Figure    \ref{fig:nonlinear_price} (left)
presents the approximate  feedback controls
 for $s=2$ and  $(t,q)\in [0,1]\t [0.5,2.5]$,
which clearly shows that the nonsmooth cost 
enhances the sparsity of the optimal strategy,
especially near the terminal time.
We then analyze the robustness of
 the FIPDE scheme 
 by exercising the feedback control
$\phi^\textrm{pre}$
computed  with  $Q_0\sim \mathcal{U}(1,2)$ 
on perturbed models with  $Q_0\sim \mathcal{U}(Q_{0,\min},Q_{0,\max})$.
Due to the absence of  analytic solution for the nonsmooth MFC problem, 
we  compare 
 the expected cost
  $J(\phi^\textrm{pre})$
  of $\phi^\textrm{pre}$ on the perturbed model 
  against the numerical approximation  $J(\phi^\textrm{pert})$
  of the optimal cost of the perturbed model,
where 
 $\phi^\textrm{pert}$
is   the feedback control
obtained by   the FIPDE scheme
with the perturbed $Q_0$.
Figure \ref{fig:nonlinear_price} (right)
depicts  the absolute performance gaps 
$|J(\phi^\textrm{pre})-J(\phi^\textrm{pert})|$
for different $Q_0$.
Due to the nonsmooth cost function,
 the optimal feedback controls 
of \eqref{eq:price_impact} 
depend nonlinearly on the law of $Q_0$,
and hence 
 the  performance gaps 
are no longer
 constant along the diagonals
(cf., Figure \ref{fig:robustness_linear_abs} (left)).
  However, 
 the absolute performance gaps remain small
 for all perturbations,
 which demonstrates the robustness of the 
FIPDE scheme in the present nonsmooth setting.

\subsection{Sparse consensus control of stochastic Cucker--Smale models}\label{SEC:CS}
%

In this section, we study optimal control of 
the  multidimensional stochastic mean-field Cucker--Smale (C-S) dynamics (see, e.g., \cite{nourian2011mean} and \cite[Chapter 4]{carmona2018probabilistic}),
where the controller  aims to enforce consensus emergence of an interactive particle system via 
(possibly sparse) external intervention. 
In its general form this leads to a nonsmooth nonconvex MFC problem.

Given a terminal time $T>0$ and a $d$-dimensional adapted control strategy
 $(\a_t)_{t \in [0,T]}$,
 we consider
 the following 
 $2d$-dimensional
 controlled 
 dynamics with additive noise,
which can be viewed as the large population limit
 of the  finite-particle model studied in \cite{cucker2007emergent}:
 for all $t \in [0,T]$,
  \begin{align}
  \label{eq:cs}
 \d x_t=v_t\,\d t, \quad  
 \d v_t=
 \bigg(\int_{\R^d \times \R^d}
\kappa(x_t,v_t,x',v')\,\mathcal{L}_{(x_t,v_t)}(\d x',\d v') + \a_t\bigg)\,\d t+\sigma \, \d W_t, 
 \end{align}
with initial state $ (x_0,v_0)\in L^2(\Omega;\R^d \times \R^d)$,
where
$W$ is a $d$-dimensional Brownian motion defined on a filtered probability space 
$(\Omega, \mathcal{F}, \{\mathcal{F}_t\}_{t\in [0,T]},\mathbb{P})$,
$\sigma\in \R^{d \times d}$
and $\kappa:\R^d \times \R^d \times \R^d \times \R^d \to \R^d$ is  given by
$$
\kappa(x,v,x',v')=\frac{K(v'-v)}{(1+|x-x'|^2)^\beta},
\quad \textnormal{with some  $\beta,K\ge 0$.}
$$
It is well-known that
for uncontrolled deterministic models (with $\sigma=0$ and $\a\equiv 0$), 
a time-asymptotic flocking behaviour (i.e., all trajectories of the velocity process 
tend to the same value as $t\to \infty$)
only appears for $\beta\le 1/2$ or for specific initial conditions if $\beta>1/2$;
see \cite{cucker2007emergent, carrillo2010asymptotic, bailo2018optimal}
and references therein.
Moreover,
as shown in 
\cite{cucker2007emergent,nourian2011mean},
 even when
 the  deterministic counterpart 
exhibits a time-asymptotic flocking
behaviour,
the additive noise may
 prevent the emergence of flocking
in the stochastic model.

The aim of the controller is to either induce
consensus on models   that would
otherwise diverge, or 
to accelerate the flocking for 
an initial configuration
that would naturally self-organise.
More precisely, 
for  given constants $\gamma_1,\gamma_2\ge 0$,
we consider minimizing the  following cost functional
\begin{align}\label{eq:cs_loss}
J(\a;(x_0,v_0))
=\mathbb{E} \left[ \int_{0}^{T} 
\Big( |v_t-\mathbb{E}[v_t]|^2+
\gamma_1 |\a_t|^2+ \gamma_2 |\a_t|_1 \Big) \, \mathrm{d}t + |v_T-\mathbb{E}[v_T]|^2
 \right]
\end{align}
over all adapted control processes $\a$ taking values in $\R^d$,
where  $|\cdot|_1$ is the $\ell_1$-norm of a given vector.
Note that in general,
 the  cost functional $J$ is neither convex nor smooth in the control process $\alpha$,
due to the nonlinear 
 interaction  kernel $\kappa$ in \eqref{eq:cs}
and the $\ell_1$-norm in \eqref{eq:cs_loss}.
However,  in the special case with
 $\beta=\gamma_2=0$ and $\gamma_1>0$,
 \eqref{eq:cs}-\eqref{eq:cs_loss}  is a LQ MFC problem
 whose  optimal feedback control can be found via 
 Riccati equations (see, e.g.,~\cite{yong2013linear,basei2017linear}).

In the subsequent two sections, we 
 demonstrate the effectiveness of 
the FIPDE method (i.e., Algorithm \ref{alg:NAG2})
to solve  \eqref{eq:cs}-\eqref{eq:cs_loss}
with
different choices of $\beta, d$ and $\gamma_2$.
For the  $m$-th NAG iteration, 
given the approximate feedback control 
$\psi^m:[0,T]\times \R^d \times \R^d \to \R^d$
with associated 
 state processes $(x^m_t,v^m_t)_{t\in[0,T]}$,
the FIPDE method   seeks functions 
${u}_1,{u}_2:[0,T]\times \R^d \times \R^d \to \R^d$
satisfying the following 
coupled system of  $2d$-dimensional  parabolic PDEs:
  for all $(t,x,v)\in [0,T]\times \R^{d} \times \R^d$,%
\begin{subequations}\label{eq:pde_cs}
\begin{alignat}{2}
( \partial_t u_1+\mathscr{L}^m u_1)(t,x,v)& =-f^m_1(t,x,v,u_2),
&&
\quad  
{u}_1(T,x,v)=0, 
\\
( \partial_t u_2+\mathscr{L}^m u_2)(t,x,v)&=-f^m_2(t,x,v,u_1,u_2),
&&
\quad 
{u}_2(T,x,v)=2(v-\mathbb{E}[v^m_T]),
\end{alignat}
\end{subequations}
with  the operator $\mathscr{L}^m$ and the source terms $f^m_1,f^m_2$ satisfying 
for all
 $\varphi\in C^{1,2}([0,T]\times \R^{2d};\R)$,
\begin{align}\label{eq:cs_Lf}
\begin{split}
(\mathscr{L}^m \varphi)(t,x,v)&=
\tfrac{1}{2}
\tr\big(\sigma \sigma^\trans (\p_{vv}\varphi) (t,x,v)\big)
+
v^\trans (\p_x \varphi) (t,x,v)
\\
&\quad 
+
\big( {\EE} [  \kappa(x,v,{x}^m_t,{v}^m_t)  ] 
+\psi^m(t,x,v)\big)^\trans (\p_v \varphi) (t,x,v),
\\
f^m_1(t,x,v,u_2)
&=
 {\EE} [ (\partial_x \kappa)(x,v,{x}^m_t,{v}^m_t)  ]u_2(t,x,v)
+{\EE}[(\p_{{x}'} \kappa)({x}^m_t,{v}^m_t,{x},{v}) u_2(t,x^m_t,v^m_t)],
\\
f^m_2(t,x,v,u_1,u_2)
&=
u_1(t,x,v) + 2(v-\mathbb{E}[v^m_t]) 
+  {\EE} [ (\partial_v \kappa)(x,v,{x}^m_t,{v}^m_t)  ] u_2(t,x,v)
\\
&\quad
+
{\EE}[(\p_{{v}'} \kappa)({x}^m_t,{v}^m_t,{x},{v}) u_2(t,x^m_t,v^m_t)],
\end{split}
\end{align}
where for each $w\in \{x,x',v,v'\}$,
$\p_w \kappa$
denotes 
 the  Jacobian matrix  such that 
$(\p_w \kappa)_{ij}=
{\p_{w_i} \kappa_j}$.
In practice, 
we approximate 
the expectations in the coefficients of  
\eqref{eq:cs_Lf}   by the corresponding empirical averages 
over particle approximations of $(x^m,v^m)$, which 
are generated by \eqref{eq:forward_samples}
with sufficiently large $N,M\in \N$.
The PDE systems with approximated coefficients will be solved by 
using the finite difference approximation in Section \ref{SEC:FDM} or the residual approximation method 
in Section \ref{SEC:NN},
depending on the problem dimension  $2d$;
see 
 Sections \ref{sec:CS_2d} and \ref{sec:CS_6d}, respectively,
 for more details.

We shall also compare the  FIPDE method with 
 the neural network-based policy gradient  (NNPG) method,
 which is 
a pure data-driven algorithm 
proposed in \cite{carmona2019convergence}
to solve MFC problems.
The NNPG method considers minimizing 
\eqref{eq:cs_loss} over a class of 
 feedback controls 
represented by 
multilayer  neural networks
with weights $\theta$, 
 and   obtains the optimal weights 
by applying  gradient descent   algorithms based on simulated trajectories  
of the state process \eqref{eq:cs}
(see Sections \ref{sec:CS_2d} and \ref{sec:CS_6d}
for more details).
As we shall see soon,
the proposed FIPDE method 
achieves more  accurate and interpretable feedback controls than
the NNPG method, in both the low-dimensional and high-dimensional settings.

\subsubsection{Consensus control of two-dimensional C-S   models}
\label{sec:CS_2d}

This section
studies  
 \eqref{eq:cs}-\eqref{eq:cs_loss} 
 for two-dimensional C-S models
 with 
different communication rates 
$\beta$.
In particular, we shall perform experiments with the following model parameters: 
$T=d=K=1$, 
$\sigma=\gamma_1=0.1$, 
$\gamma_2=0$,
$\beta\in \{0,10\}$, 
and 
$(x_0,v_0)$ follows  the 
 two-component  Gaussian mixture distribution
 with mixture weights $(0.5,0.5)$, 
mean $(1.2, 1.8)$  and covariance $0.01\mathbb{I}_2$ for  component 1,
and mean $(1.8, 1.2)$  and covariance $0.01\mathbb{I}_2$ for component 2,
where $\mathbb{I}_2$ is the $2\t 2$ identity matrix.

 We 
 now discuss the implementation details of the 
FIPDE algorithm 
with  finite difference approximation
for the present two-dimensional setting.
The FIPDE algorithm 
is initialized with $\tau=1/6$ and $\phi^0\equiv 0$. 
At the $m$-th NAG iteration, 
given an approximate feedback control $\psi^m$,
we first generate particle approximations
$(x^{m,l},v^{m,l})_{l=1,\ldots, 10^4}$
of the state process $(x^{m},v^m)$ with time stepsize $\Delta t= 1/50$
(i.e., \eqref{eq:forward_samples} with $N=10^4$ and $M=50$), 
and  replace the   expectations in 
\eqref{eq:cs_Lf}
by the 
empirical averages over particles.
The PDE system with approximate coefficients 
is then localized on the  domain  $ \cD=[0,5]\times [0,4]$
{\color{black}
with     boundary conditions being the terminal values,
i.e., $(u_1(t,x,v),u_2(t,x,v))=(0, 2(v-\tfrac{1}{N}\sum_{l=1}^Nv^{m,l}_T))$
for all $(t,x,v)\in [0,T]\times \p\cD$.}
We  further 
 construct 
a semi-implicit first-order monotone scheme
\eqref{eq:full_discrete}
for the localized  system \eqref{eq:pde_cs}
by 
adopting
 the time stepsize $\Delta t=1/50$ and mesh sizes $h_x=5/50$, $h_v=4/50$,
and discretizing 
the first-order derivatives in \eqref{eq:cs_Lf}
via the  upwind finite difference
and
  the second-order derivative in \eqref{eq:cs_Lf}
via the central difference.

\begin{figure}[!ht]
    \centering
    \includegraphics[  width=0.3\columnwidth,height=4.4cm]{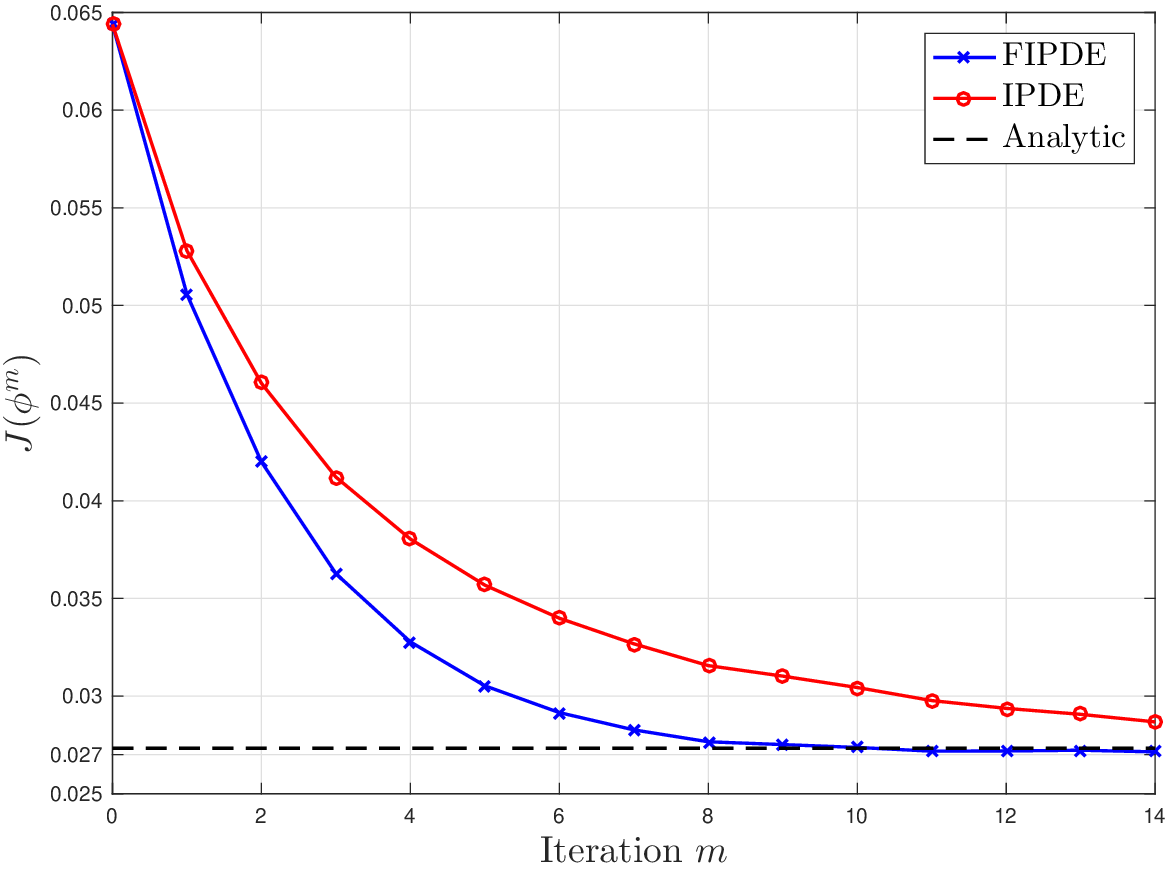}
         \hspace{1.5cm}
       \includegraphics[ trim=30 5 40 20, clip, width=0.32\columnwidth,height=4.5cm]{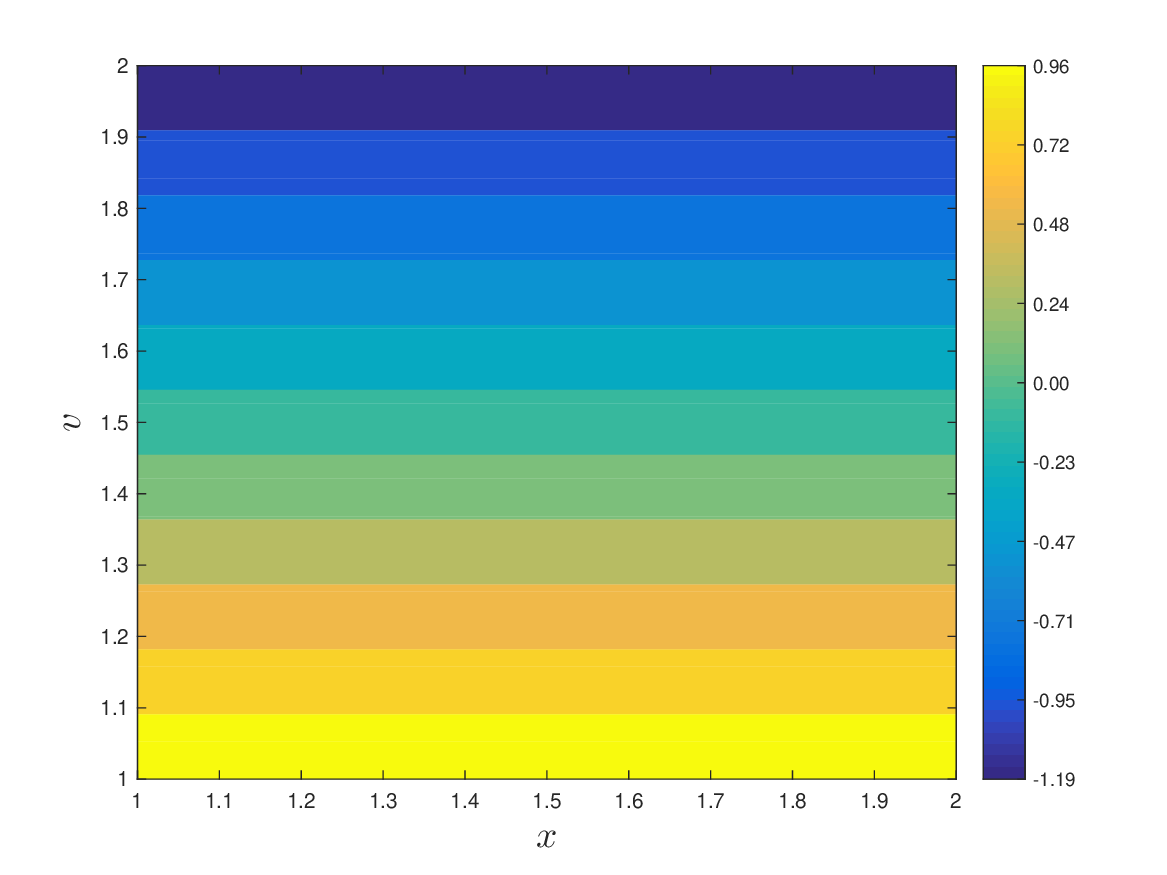}

    \caption{
  Numerical results 
  of the FIPDE scheme 
  for  two-dimensional  C-S model with 
$\beta=0$; 
from left to right: 
convergence of the FIPDE method 
in terms of NAG iterations,
and 
the approximate feedback control at $t=0$.
 }
    \label{fig:cs_2d_beta0}
 \end{figure}

For the LQ MFC problem  \eqref{eq:cs}-\eqref{eq:cs_loss} with $\beta=0$,
we compare the approximate feedback control  from the FIPDE scheme with
the optimal feedback control
given by
\bb\label{eq:cs_2d_fb}
\phi^\star(t,x,v)=-\tfrac{a_t}{2\gamma_1}(v-\sE[v^\star_t]), \quad (t,x,v)\in [0,T]\t \sR\t \sR,
\ee
where $a:[0,T]\t \sR$ satisfies  $a'_t-2Ka_t-\frac{1}{2\gamma_1}a_t^2+2=0$ with $a_T=2$,
and $v^\star$ is the optimal velocity process
 (see, e.g., \cite{yong2013linear,basei2017linear}).
 The blue line in 
 Figure \ref{fig:cs_2d_beta0}  (left) presents the expected costs of the approximate feedback controls
 obtained by the FIPDE method at all NAG iterations, 
 which converge exponentially  to 
   the optimal cost 
   in terms of the  number of NAG iterations;
a linear regression of the data shows 
 the absolute error  for the $m$-th iteration is 
  of the magnitude $\cO(0.55^m)$.
  {\color{black}
The    FIPDE method  
is     benchmarked 
 against 
the iterative PDE-based (IPDE) algorithm Algorithm \ref{alg:NAG1} 
  (i.e., $\psi^m=\phi^m$ for all $m$),
  whose  convergence  is shown by 
the red line in  Figure \ref{fig:cs_2d_beta0}  (left).
Note that the FIPDE and IPDE methods take a similar time 
  to perform one gradient descent iteration,
  but the FIPDE method 
  requires much fewer iterations to achieve high accuracy.
This shows that the momentum step indeed accelerates 
the algorithm convergence, 
and justifies the terminology ``fast" in FIPDE.
  }
We further depict  the approximate feedback control $\phi$ (generated by the last NAG iteration)
with $t=0$ and $(x,v)\in [1,2]^2$
in  Figure \ref{fig:cs_2d_beta0}  (right),
which approximates the optimal  feedback
 control $\phi^\star(0,\cdot,\cdot)$ (see \eqref{eq:cs_2d_fb}) 
 with a relative error of $1.7\%$ in  the $L^2$-norm. 
One can clearly observe that the approximate feedback control 
captures the affine structure of
 $\phi^\star$ {in $v$, constant in $x$}.

Then we proceed to study 
the (nonconvex) MFC problem  \eqref{eq:cs}-\eqref{eq:cs_loss} with $\beta=10$,
where  the uncontrolled velocity process
does not exhibit
a  flocking behaviour;
see
Figure \ref{fig:cs_state_2d} 
in Appendix \ref{appendix:cs}.
In the following, we  apply the FIPDE method (with the same discretization parameters as above) 
to design feedback controls
and induce consensus for such models.
We shall further 
 benchmark the performance of 
 the FIPDE method
 against a neural network-based policy gradient (NNPG) method
 (see Algorithm 1 in \cite{carmona2019convergence}),
 whose implementation details are  given in Appendix \ref{appendix:cs}. 
In particular, we  
 seek an approximate feedback control among a family of neural networks 
 whose   flexibility is sufficient to capture  the nonlinearity of the optimal feedback control
 in $(t,x,v)$,
 and obtain the optimal 
 neural network  approximation
 by running
the   Adam   algorithm
 \cite{kingma2014adam} with sufficiently many iterations.

\begin{figure}[!ht]
    \centering
         \includegraphics[  width=0.3\columnwidth,height=4.4cm]{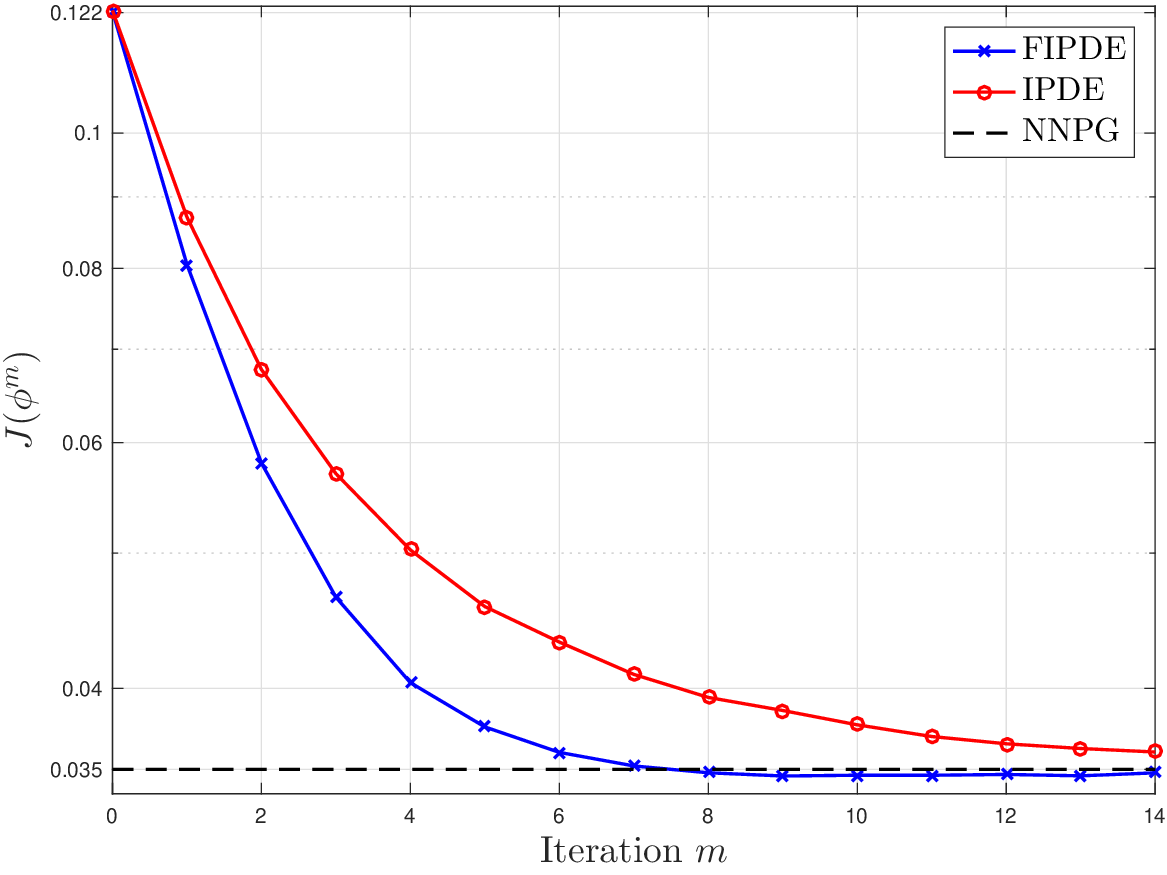}
         \;
	 \includegraphics[ trim=30 5 40 20, clip, width=0.32\columnwidth,height=4.5cm]{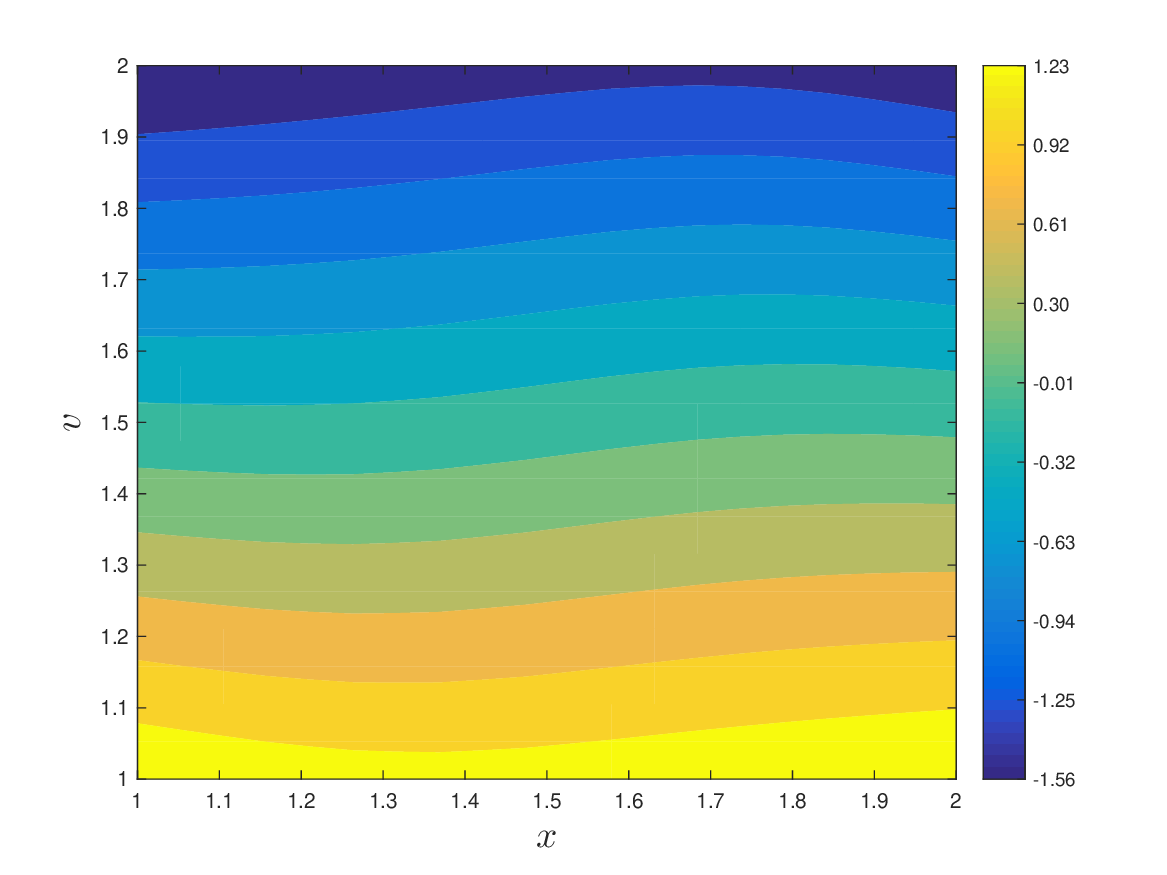}
	   \;
       \includegraphics[ trim=30 5 40 20, clip, width=0.32\columnwidth,height=4.5cm]{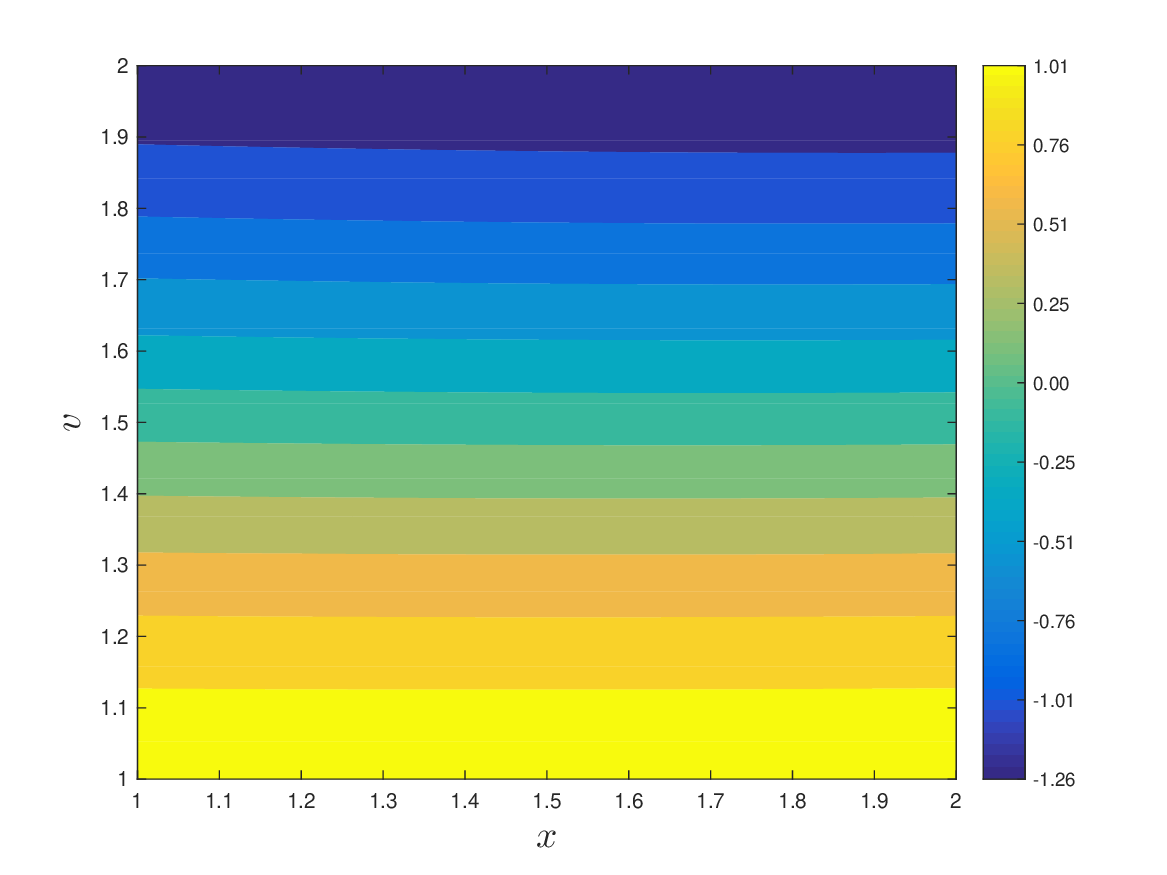}

    \caption{
  Numerical results 
  of the FIPDE   and
 NNPG methods
  for  two-dimensional  C-S model with 
$\beta=10$; 
from left to right: 
convergence of the FIPDE scheme 
in terms of NAG iterations,
the 
feedback control of the FIPDE method at $t=0$,
and 
the 
feedback control of the NNPG method at $t=0$.
 }
    \label{fig:cs_2d_beta10}
 \end{figure}

Figure \ref{fig:cs_2d_beta10}
 exhibits the approximate  feedback controls 
 and the associated costs 
 from the FIPDE and NNPG methods. 
 From Figure \ref{fig:cs_2d_beta10} (left),
we see that
despite the nonconvexity of the control problem, 
 the FIPDE method  converges exponentially 
 in the value function approximation  as the iteration index $m$ tends to infinity,
with the absolute  error $\cO(0.5^m)$. 
{\color{black}
Compared with the IPDE method,
the momentum step in  Algorithm \ref{alg:NAG2} 
accelerates  the algorithm's convergence also 
in the present nonconvex setting.}
We refer the reader to Figure  \ref{fig:cs_state_2d}  in Appendix \ref{appendix:cs}
 for the flocking behaviour of the controlled velocity process 
obtained by the FIPDE method.

More importantly, 
the FIPDE method generates  a more  interpretable 
approximate feedback control 
compared to the NNPG method,
which helps us understand the mechanism of the optimal control process. 
Note that 
for the uncontrolled C-S model   \eqref{eq:cs} with  $\beta>0$,
the interaction kernel 
  $\kappa$  indicates that 
  a particle's velocity is largely affected by the velocities of particles 
 in the nearest neighborhood.
 Consequently, for a given particle whose velocity is above the average velocity,
 the further it is away from the population with small velocity,
 the less internal attraction exists for the particle's velocity 
to the  average velocity,
and the stronger external intervention is required 
to induce consensus
(and similarly for particle whose  velocity is below the average velocity).
 Such a nonlinear dependence of the optimal feedback control
on the variable $x$
is correctly captured by the approximate feedback control $\phi$ of  the FIPDE method,
whose values at $t=0$ are depicted in Figure \ref{fig:cs_2d_beta10} (middle).
Recall that  at $t=0$, 
the average  velocity is $1.5$,
and  particles with 
velocity below  and above  $1.5$ 
  cluster around the points $x=1.8$ and $x=1.2$, respectively. 
Hence,
the absolute magnitude of the function 
$x\mapsto \phi(0,x, v)$ 
 is minimized near $1.8$ if  $v>1.5$
and near  $1.2$ if   $v<1.5$,
which confirms our theoretical understanding of the model. 
By contrast,
the NNPG method produces a less interpretable approximate feedback control
with no clear dependence on the variable $x$, as shown in Figure \ref{fig:cs_2d_beta10} (right).

\subsubsection{Sparse consensus control of six-dimensional C-S models}
\label{sec:CS_6d}
This section examines the performance of 
the FIPDE method
(i.e., Algorithm \ref{alg:NAG2})
for solving
 \eqref{eq:cs}-\eqref{eq:cs_loss}
 with possibly nonsmooth costs
 in a six-dimensional setting, 
where  the corresponding  PDE systems   \eqref{eq:pde_cs}
are solved 
by using the residual approximation approach    introduced in Section \ref{SEC:NN}.
We  carry out numerical experiments  with  
$d=3$, $T=K=1$, $\sigma=0.1 \mathbb{I}_3$, $\gamma_1=0.3$,
$(x_0,v_0)\sim \cU([0,1]^6)$  (the uniform distribution on $[0,1]^6$)
but  different  choices of   $\beta, \gamma_2\ge 0$.
For the sake of presentation,
we shall only specify the neural network architectures 
used in our  computation,
and refer the reader to  Appendix \ref{appendix:cs}
for the detailed  implementation of the empirical risk minimization problem 
\eqref{eq:E_empirical}.

Let us start with the LQ MFC problem by choosing $\beta=\gamma_2=0$.
We  initialize Algorithm \ref{alg:NAG2} with $\phi^0=0$ and $\tau=1/5$,
and adopt   feedforward neural networks with the same architecture
to approximate each component of 
$\phi^m, \psi^m:[0,T]\t \sR^6\to \sR^3$
 and 
$u^m:[0,T]\t \sR^6\to \sR^6$
for all NAG iterations.
More precisely, 
for each $m$,
we have  $\phi^m=(\phi^\theta_1,\ldots \phi^\theta_3)$,
$\psi^m=(\psi^\theta_1,\ldots \psi^\theta_3)$ and $u^m=(u^\theta_1,\ldots u^\theta_6)$,
where $\phi^\theta_i,\psi^\theta_i, u^\theta_i$ are some
fully-connected networks with the sigmoid activation function, depth 2 (1 hidden layer),
and the dimensions of the input, output and hidden layers being 
7, 1 and  20, respectively.
We then  apply the residual approximation method  
to solve   \eqref{eq:pde_cs} on the spatial  domain $\cD=[-1,2]^6$.
Note that  
approximating  each  component
 of the solutions individually
 allows us to
capture the  heterogeneity among   components
with  
 shallow networks,
and consequently makes the empirical residuals relatively easy to optimize.
Moreover,  we set $\bar{u}=0$ on $[0,T]\t \cD$ 
 for the decomposition \eqref{eq:decomposition},
 since in the present LQ case, 
 there is no obvious candidate solution to  \eqref{eq:pde_cs},
 except the exact solution based on Riccati equations.

We carry out the FIPDE method with 15 NAG iterations
and summarize
the numerical results in   Table \ref{table:architecture_impact_6d}.
It is clear that 
  the approximate controls and their associated costs 
  of the FIPDE method converge rapidly to the optimal control process
  and the optimal cost function obtained by  Riccati equations. 
Moreover, 
a sensitivity analysis of  
  the neural  network architectures 
    indicates  that 
the convergence of the FIPDE method is very robust with respect to
 the network depth $L$,  
 the dimension $H$ of  hidden layers, and the (smooth) activation functions.
Figure  \ref{fig:cs_v_6d_lq}  in Appendix \ref{appendix:cs}
depicts 
 the trajectories of the uncontrolled and controlled velocity processes,
 which demonstrates that the approximate feedback control of the FIPDE method
 effectively accelerates the emergence of the flocking behaviour.

\begin{table}[!ht]
 \renewcommand{\arraystretch}{1.05} 
\centering
\caption{ Impact of the  network architecture  on the convergence  of the FIPDE method.
Shown are: (a) the absolute error  $|J(\phi^m)-J^\star|$ of the numerical value function $J(\phi^m)$,
and (b) the absolute error  $ (\mathbb{E}[ \sup_{t \in [0,T]} |\alpha_t^m - \alpha_t^\star |^2])^{\frac{1}{2}}$ of the numerical control process $\alpha^m$.
}
\label{table:architecture_impact_6d}
\begin{tabular}[t]{@{}c cccccc c@{}}
\toprule
 & NAG Itr $m$ & 1 & 3 & 6 & 9 & 12 & 15 \\ \midrule

$L=2$, $H=20$, $\varrho=\textrm{sigm}$ & 
(a) &
0.0401    & 0.0121  &   0.0050   & 0.0023 &   0.0010 
&	0.0007
 \\
& (b) & 0.2650 & 0.1549  & 0.1015 & 0.0632 &  0.0400 
&	0.0283
\vspace{1mm}
\\

$L=4$, $H=20$, $\varrho=\textrm{sigm}$ & 
(a) &
0.0396    &    0.0122    &  0.0045    & 0.0019    &    0.0008    &    0.0005    
\\
& (b) &
0.2676     & 0.1667     & 0.1082     & 0.0721     & 0.0470    & 0.0346     
\vspace{1mm}

\\

$L=2$, $H=40$, $\varrho=\textrm{sigm}$ &
(a) &
0.0395    &    0.0127    &     0.0053    &   0.0024    &   0.0012    &    0.0009     
\\
& (b) &
0.2638     &  0.1594     &  0.1015     &  0.0663     &  0.0412     &  0.0300     
\vspace{1mm}

\\
$L=2$, $H=20$, $\varrho=\tanh$ &
(a) &
0.0412 &     0.0197  &     0.0086  &     0.0041  &     0.0022  &     0.0013  
\\
& (b) &
0.2659 &     0.1936  &    0.1261  &    0.0860  &    0.0557  &    0.0346  
\\
\bottomrule
\end{tabular}
\end{table}%

\begin{figure}[!ht]
    \centering
         \includegraphics[  width=0.3\columnwidth,height=4.4cm]{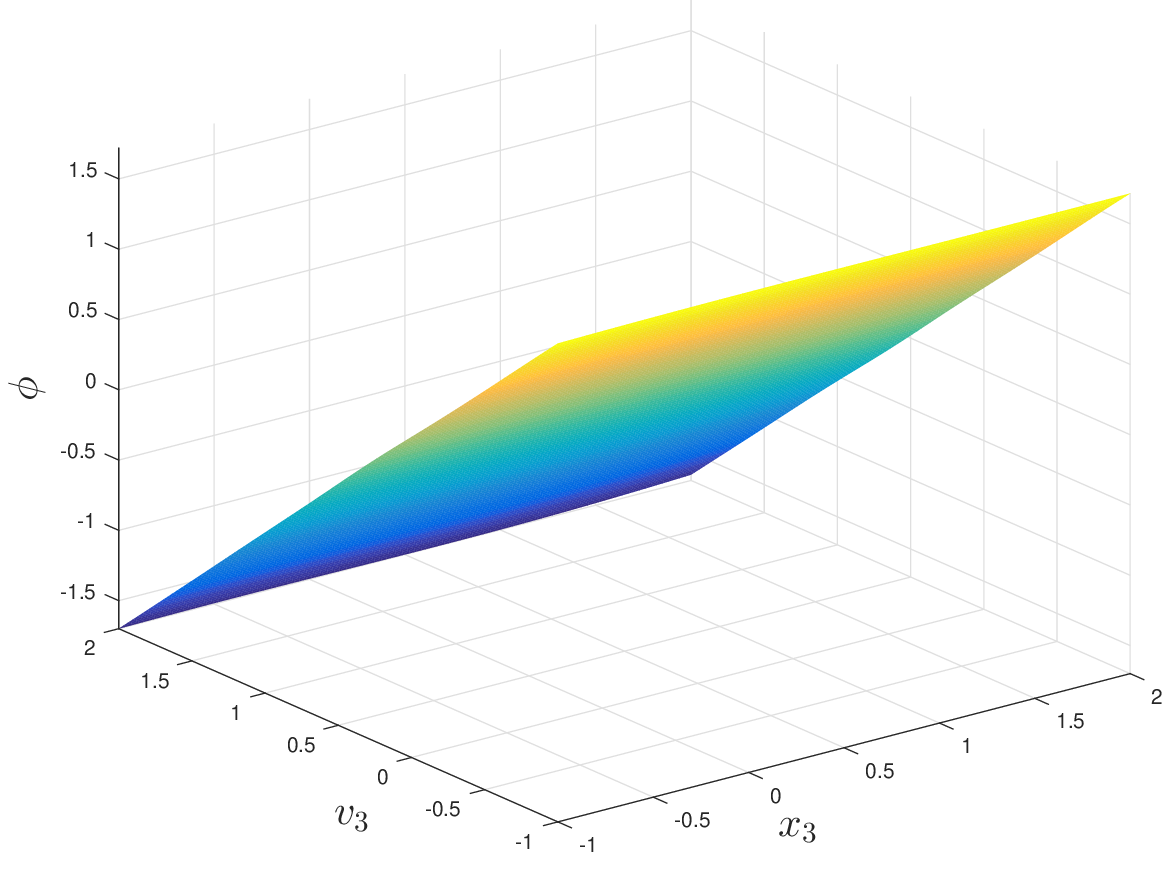}
         \quad
	 \includegraphics[ width=0.3\columnwidth,height=4.4cm]{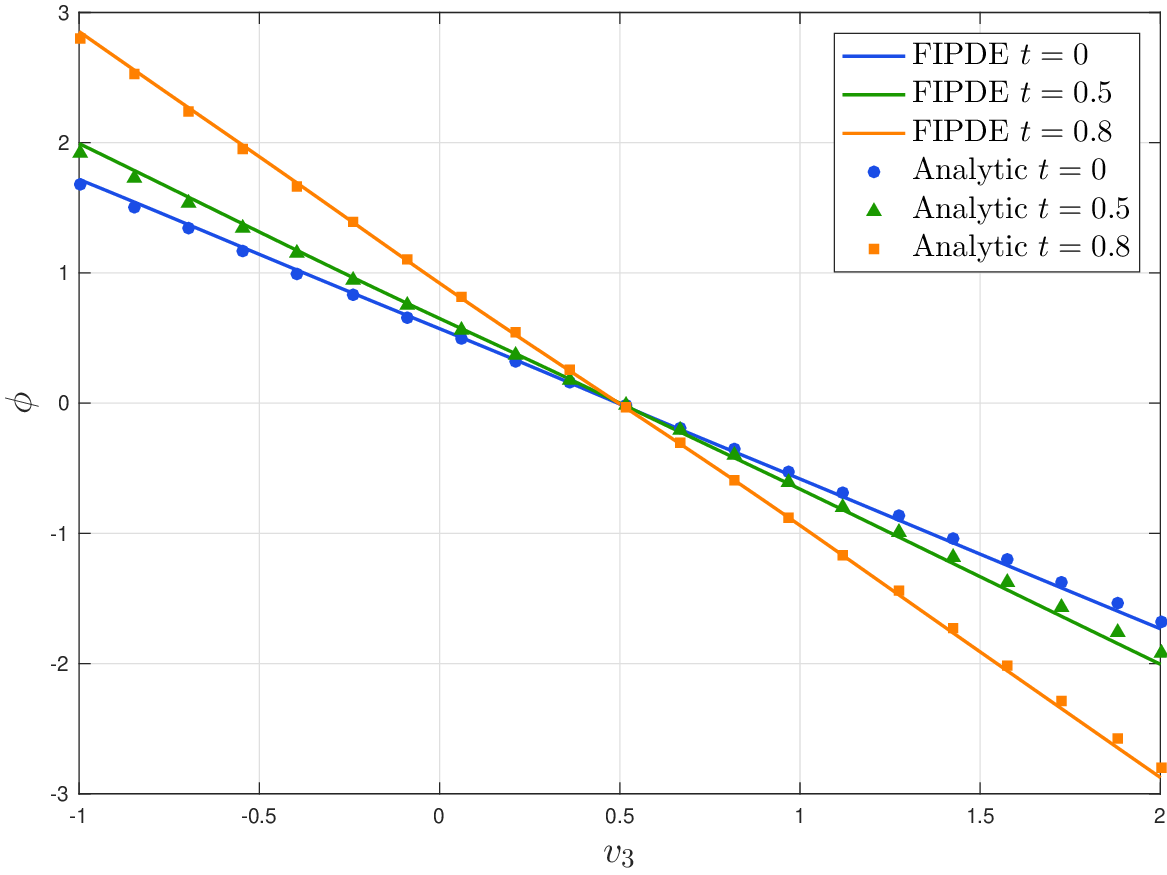}
	   \quad
       \includegraphics[  width=0.3\columnwidth,height=4.4cm]{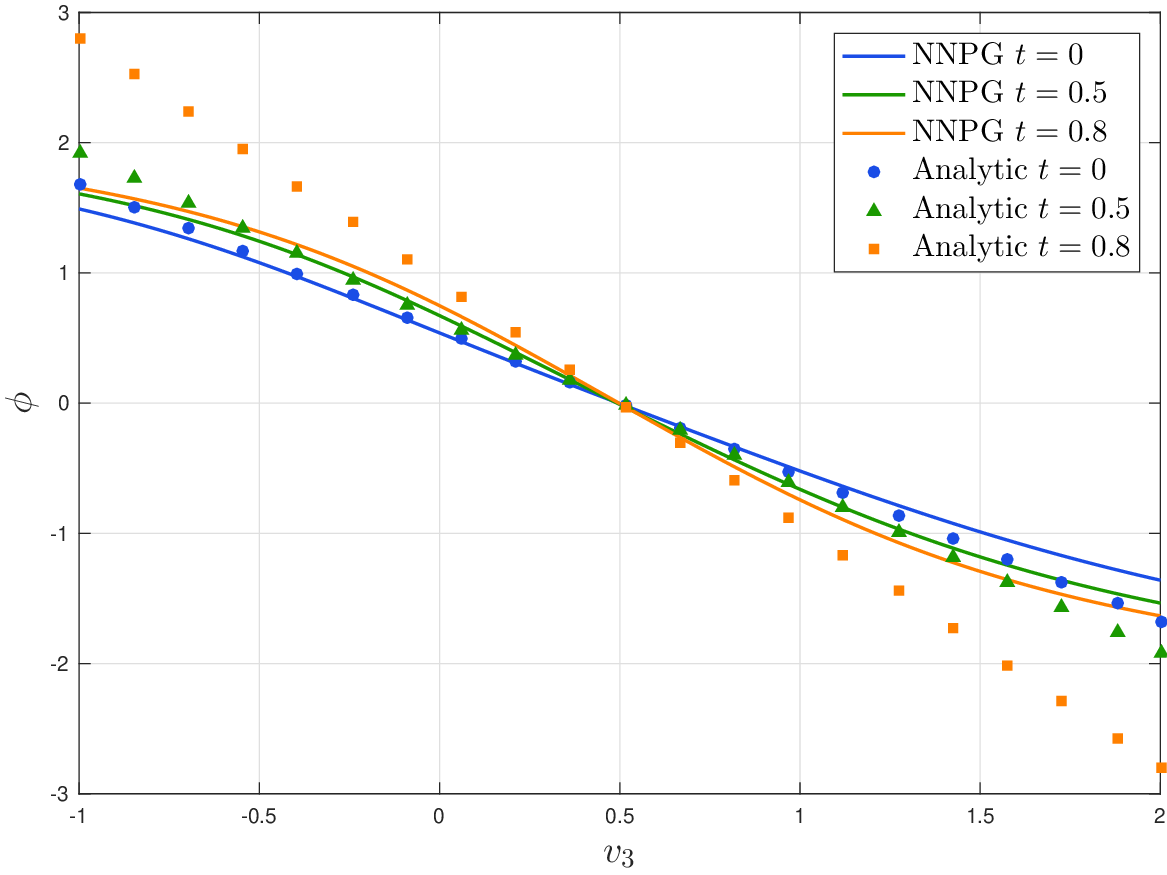}

    \caption{
  Feedback controls  for  
  six-dimensional  C-S model with 
$\beta=\gamma_2=0$; 
from left to right: 
feedback control  of the FIPDE scheme 
 at  $t=0$ and different $(x_3,v_3)$,
 feedback controls  of the FIPDE scheme 
 and the analytic solution at different $(t,v_3)$,
and  feedback controls  of the NNPG scheme 
 and the analytic solution at different $(t,v_3)$,
where the remaining components are set to  $0.5$.
 }
    \label{fig:cs_6d_beta0_fb}
 \end{figure}

We then examine the  accuracy of the FIPDE method for approximating the optimal feedback control $\phi^\star$,
which satisfies 
$\phi^\star(t,x,v)=A_t (v-\sE[v^\star_t])$ for all $(t,x,v)\in [0,T]\t \sR^3\t \sR^3$, 
with $A\in C([0,T];\sR^{3\t 3})$ being the solution to a differential Riccati equation,
and $v^\star$ being the optimal velocity process. 
As shown in  Figure \ref{fig:cs_6d_beta0_fb} (left),
although 
 the FIPDE method is implemented with
 feedforward neural networks taking inputs $(t,x,v)$,
the resulting 
 approximate feedback control correctly captures the important structures 
  of the optimal control $\phi^\star$,
  i.e., the 
     affineness  in $v$ and the independence of $x$. 
To further highlight the advantage of the PDE-based solver over pure  data-driven approaches,
we  implement the NNPG method with the same neural network architecture (see {Appendix \ref{appendix:cs}} for more details),
and compare the resulting  feedback control against the FIPDE method and the analytic solution. 
Comparing Figure \ref{fig:cs_6d_beta0_fb} (middle) and (right),
we see the FIPDE method  recovers 
 the exact feedback control accurately
on {(cross-sections of)} the entire space-time domain,
while the NNPG method only approximates the exact control 
on certain sub-domains (depending on the trajectories of the optimal state process), 
and fails to capture the time dependence of the exact control.

Now we proceed to apply the FIPDE method to the  nonsmooth MFC problem  \eqref{eq:cs}-\eqref{eq:cs_loss}
with $\beta=0$ and $\gamma_2\in \{0.1, 0.3\}$.
As initial feedback control $\phi^0$
for Algorithm \ref{alg:NAG2}
we use the approximate feedback control for the LQ problem (with $\gamma_2=0$) 
obtained by the FIPDE method,
and choose the stepsize $\tau=1/5$. 
At the $m$-th NAG iteration, 
 the residual approximation (RA) approach in Section \ref{SEC:NN}
is  implemented 
to solve  \eqref{eq:pde_cs}  on the domain $\cD=[-1,2]^6$.
In particular, 
we decompose  
the numerical solution 
into the form  $\bar{u}+\tilde{u}^m$,
where $\bar{u}:[0,T]\t \sR^6\to \sR^6$ is the approximate decoupling field of the above LQ problem
obtained by the FIPDE method,
 and $\tilde{u}^m:[0,T]\t \sR^6\to \sR^6$ is 
an unknown nonlinear residual correction.
We then approximate 
the components of $\tilde{u}^m$
by neural  networks with the sigmoid activation function, depth 4 (3 hidden layers),
and the dimensions of the input, output and hidden layers being 
7, 1 and  40, respectively,
and determine 
 the  optimal neural network representation 
 by minimizing  
 the  empirical  loss \eqref{eq:E_empirical}
 with  100 SGD   iterations (see Appendix \ref{appendix:cs}). 
 In the following, we refer to the method
 as ``RA with 100 SGD''.

To demonstrate the  efficiency  
of the RA approach, 
we  compare the performance of the above FIPDE method 
to those of the Direct Method and the NNPG method. 
In the Direct Method, we 
do not decompose  numerical solutions
 in a separable form   (i.e., we set $\bar{u}=0$ in \eqref{eq:decomposition}),
 and directly minimize  the  residuals of  \eqref{eq:pde_cs}
 over    neural networks.
For each NAG iteration,
we choose the same 
trial functions  (4-layer networks with hidden width 40),
numbers of training samples
and  learning rates of the SGD algorithm
as those of the FIPDE method, 
and perform $q$ SGD iterations
with different choices of $q\in \N$,
which will be referred to as 
``DM with $q$ SGD'' in the following discussion. 
In the NNPG method, 
we  minimize the objective \eqref{eq:cs_loss}
over all feedback controls parametrized by   4-layer  networks 
with the same architecture
 (see {Appendix \ref{appendix:cs}} for more details).

\begin{figure}[!ht]
    \centering
         \includegraphics[  width=0.35\columnwidth,height=5cm]{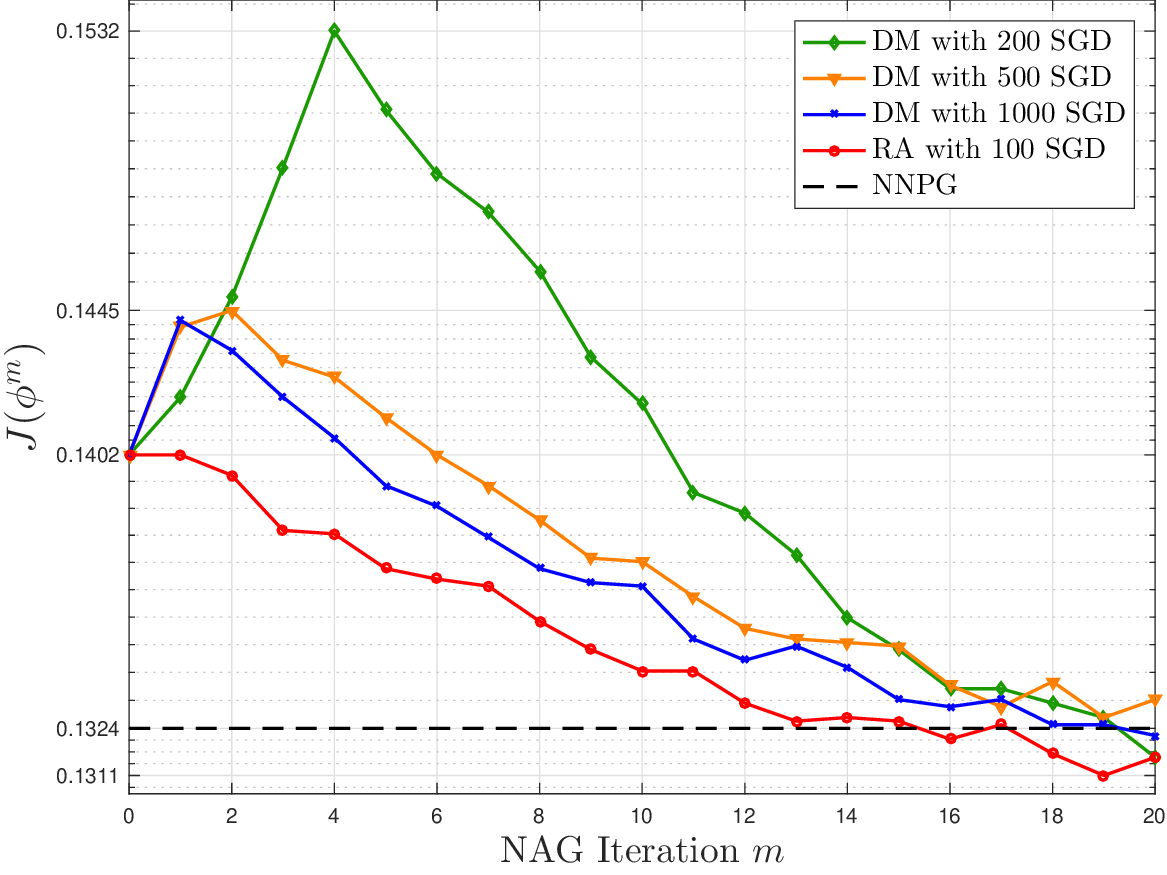}

    \caption{
   Efficiency improvement of the residual approximation approach  over  the Direct Method
   for  
  six-dimensional  C-S model with 
$\beta=0$ and $\gamma_2=0.1$.
    }
    \label{fig:value_conv_6d_l1_ra}
 \end{figure}

Figure \ref{fig:value_conv_6d_l1_ra} depicts the convergence of  value functions from  ``RA with 100 SGD''
and ``DM with $q$ SGD''  (with $q\in \{200,500,1000\}$)
for solving  \eqref{eq:cs}-\eqref{eq:cs_loss}
 with $\gamma_2=0.1$,
 where we take the value function from the NNPG method  as a reference value. 
It clearly demonstrates that
 the  residual approximation  requires significantly 
 less number of SGD iterations than the Direct Method
 for solving  the PDE system \eqref{eq:pde_cs} at each NAG iteration;
in particular, 
the approximate feedback controls of
 ``RA with 100 SGD'' (the red line)
at all NAG iterations
  consistently achieve lower costs than those of ``DM with $1000$ SGD'' (the blue line).
  {\color{black}
Observe that for each NAG iteration,
  the computation times of 
the residual approximation approach and      the Direct Method 
scale linearly with the number of SGD iterations.
Thus 
the  residual decomposition \eqref{eq:decomposition} 
reduces the total computational time 
of the Direct Method
by roughly a factor of 10.
}
  A similar efficiency enhancement of the residual approximation approach   
  has also been observed in the nonsmooth problem \eqref{eq:cs}-\eqref{eq:cs_loss}
with $\beta=0$ and $\gamma_2=0.3$.
  
\begin{figure}[!ht]
    \begin{subfigure}{.31\textwidth}
        \centering
        \includegraphics[width=\linewidth,height=4.4cm]{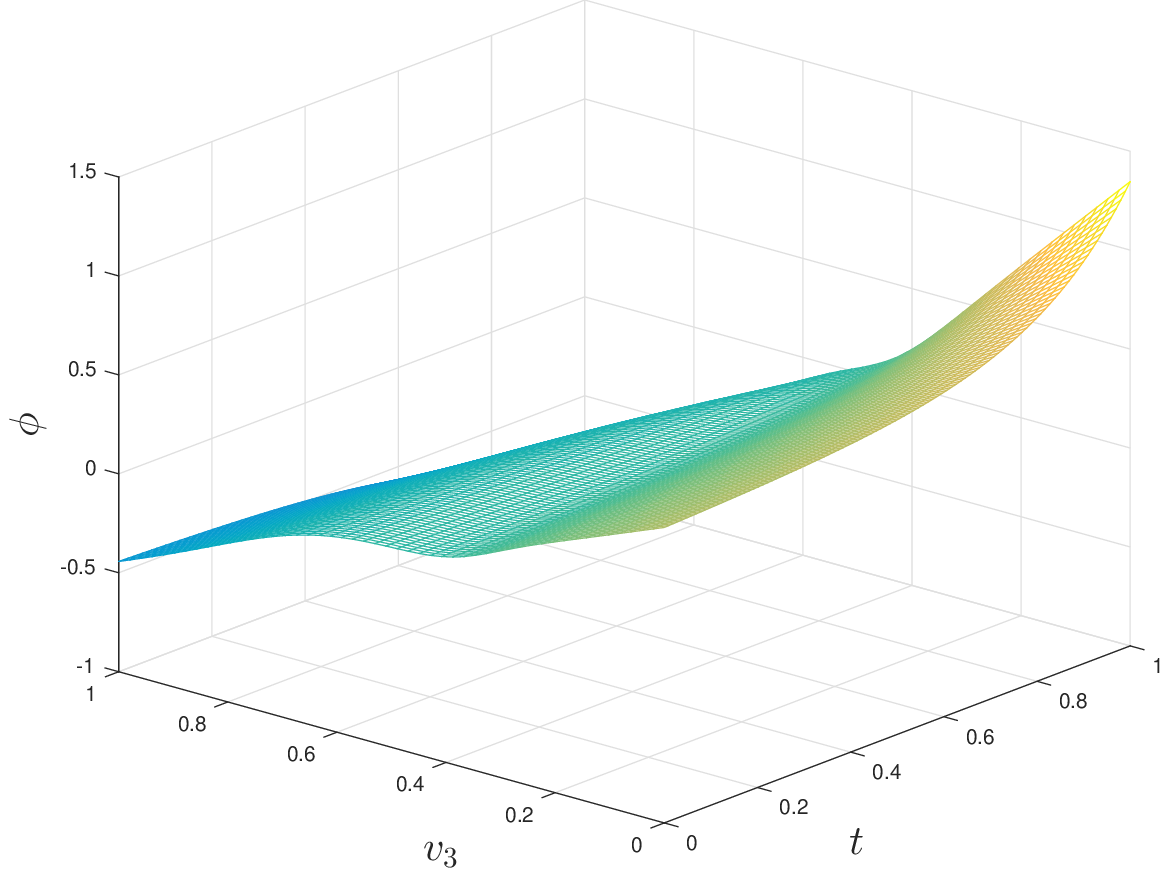}
        \caption{FIPDE for $\gamma_2=0.1$}
    \end{subfigure}\quad
    \begin{subfigure}{.31\textwidth}
        \centering
         \includegraphics[width=\linewidth,height=4.4cm]{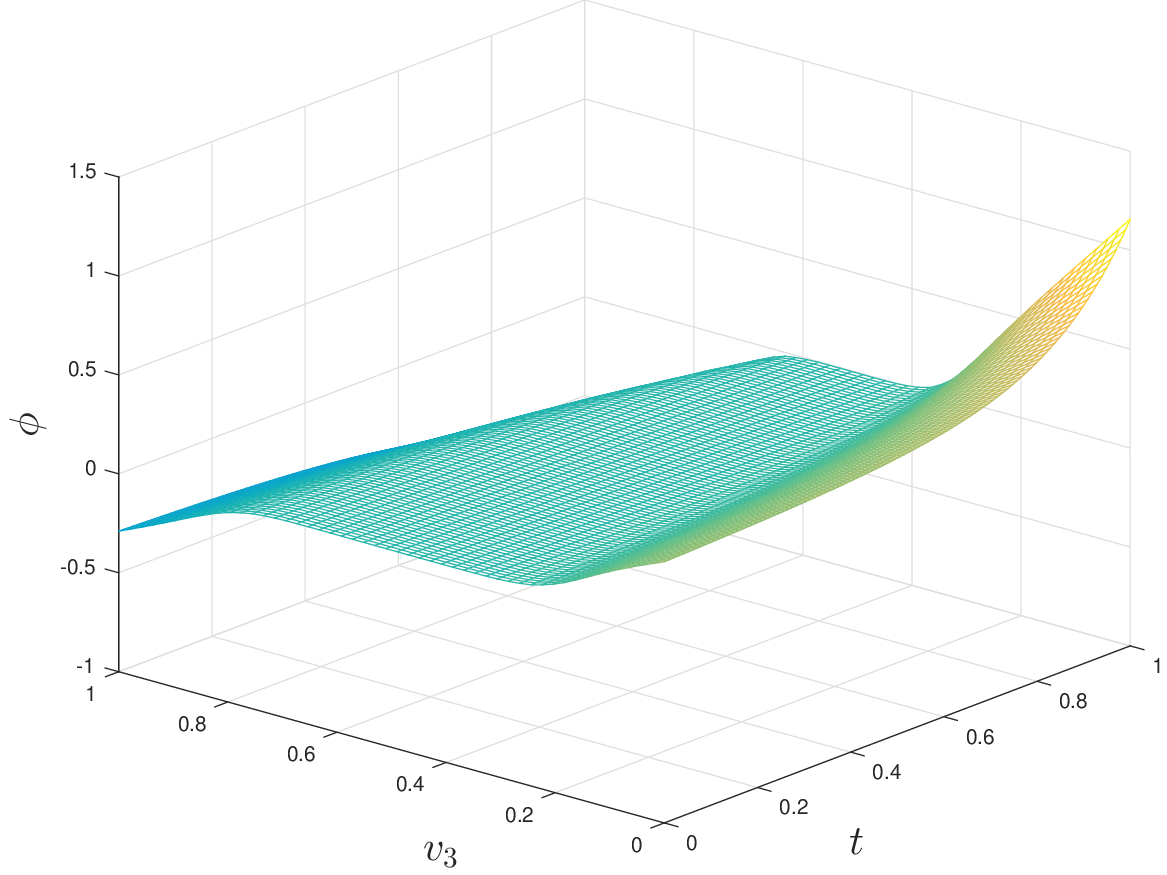}
        \caption{FIPDE for $\gamma_2=0.3$}
    \end{subfigure}\quad
    \begin{subfigure}{.31\textwidth}
        \centering
         \includegraphics[width=\linewidth,height=4.4cm]{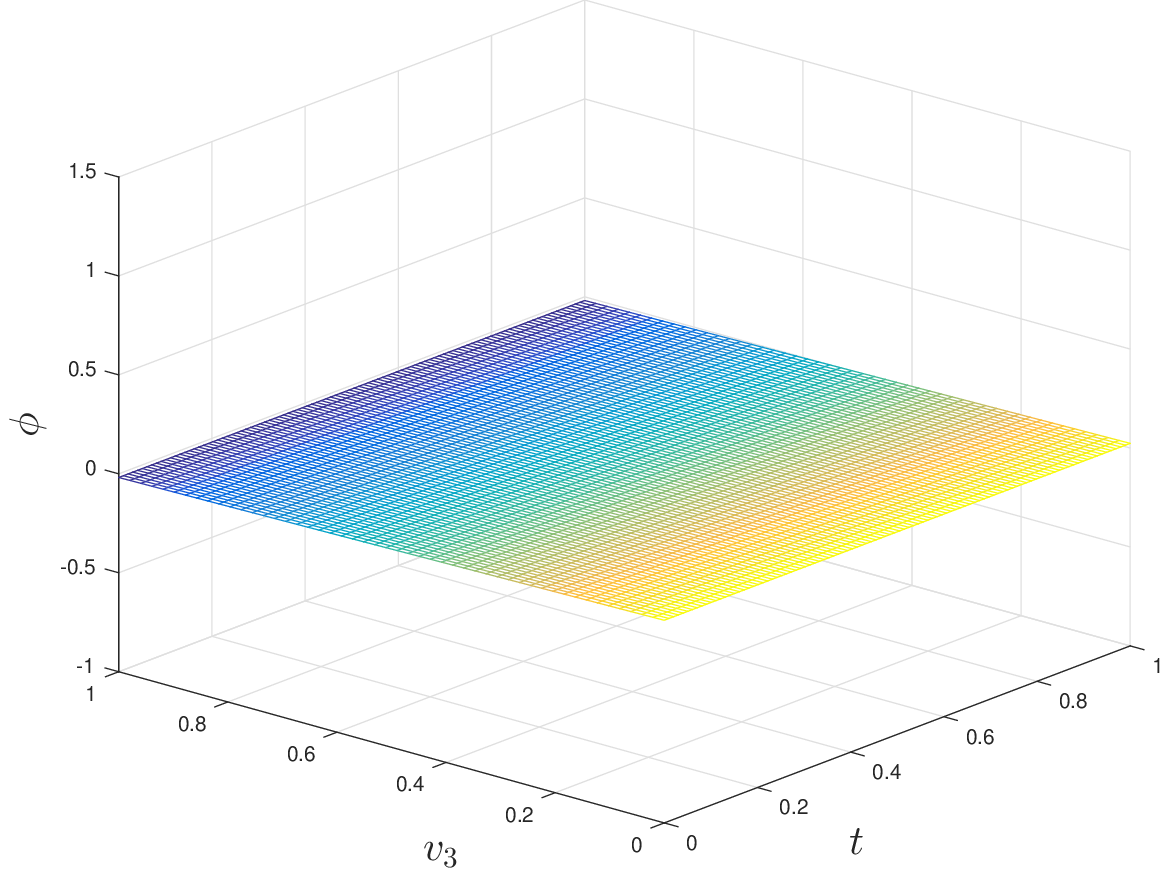}
        \caption{NNPG for $\gamma_2=0.3$}
    \end{subfigure}
    \caption{Feedback controls  for  
  six-dimensional  C-S models with 
$\beta=0$ and $\gamma_2>0$,
where the controls are evaluated at different $(t,v_3)$
with the remaining components being   $0.5$.
}
        \label{fig:value_conv_6d_l1}

\end{figure}

 Figure  \ref{fig:value_conv_6d_l1}
 presents the approximate feedback controls 
 of the nonsmooth problem   \eqref{eq:cs}-\eqref{eq:cs_loss}
 with different $\gamma_2>0$,
 obtained by using the FIPDE  and  NNPG methods.
 Due to the $\ell_1$-norm in \eqref{eq:cs_loss},
the optimal  feedback control is nonlinear  in $v$
and admits a sparse structure. That is, 
the optimal decision  is to 
 act only
on the particles whose velocities are far from
the mean and to steer them to consensus, 
without intervening with the particles near the consensus manifold
\cite{caponigro2013sparse}.
As shown in Figure \ref{fig:value_conv_6d_l1} (left) and (middle),
the approximate controls of 
the FIPDE method correctly capture the sparse features  of the optimal control,
where   the zero-control region shrinks  as   the terminal time approaches
and expands as the  parameter $\gamma_2$ increases.
On the other hand, the NNPG method results in   a  (suboptimal) non-sparse linear strategy, 
which suggests to control all agents with mild strength. 
It also fails to capture the time dependence of the optimal control, 
as already observed in the LQ setting (see Figure \ref{fig:cs_6d_beta0_fb}).
Consequently, for $\gamma_2=0.3$, 
 the  control of the NNPG method 
has
a
 $2.4\%$ higher
 expected cost \eqref{eq:cs_loss} 
  than that of the FIPDE method (in terms of the relative error);
  the expected costs for the NNPG and FIPDE methods are 
0.1635 and 0.1596, respectively. 

\begin{figure}[ht]
    \begin{subfigure}{.31\textwidth}
        \centering
        \includegraphics[width=\linewidth,height=4.4cm]{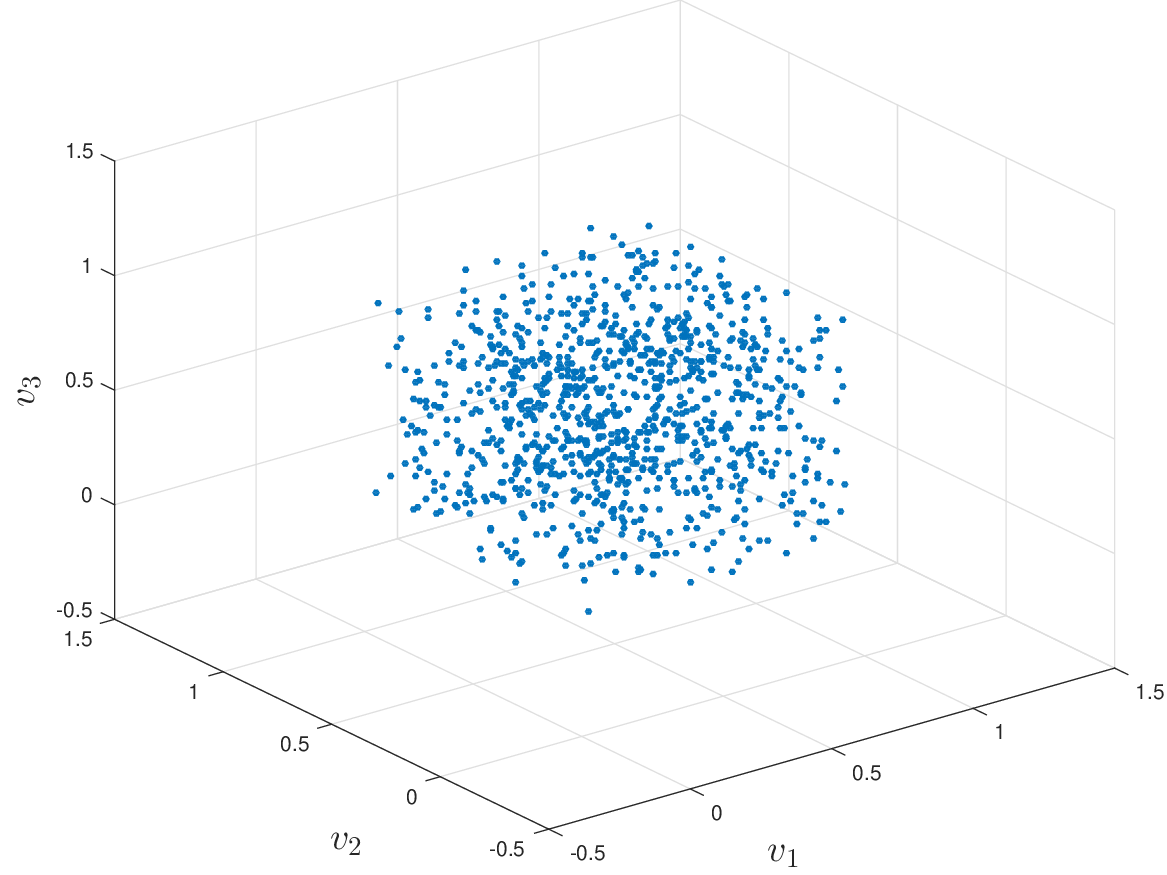}
        \caption{Uncontrolled velocity, $t=0$}
    \end{subfigure}\quad
    \begin{subfigure}{.31\textwidth}
        \centering
         \includegraphics[width=\linewidth,height=4.4cm]{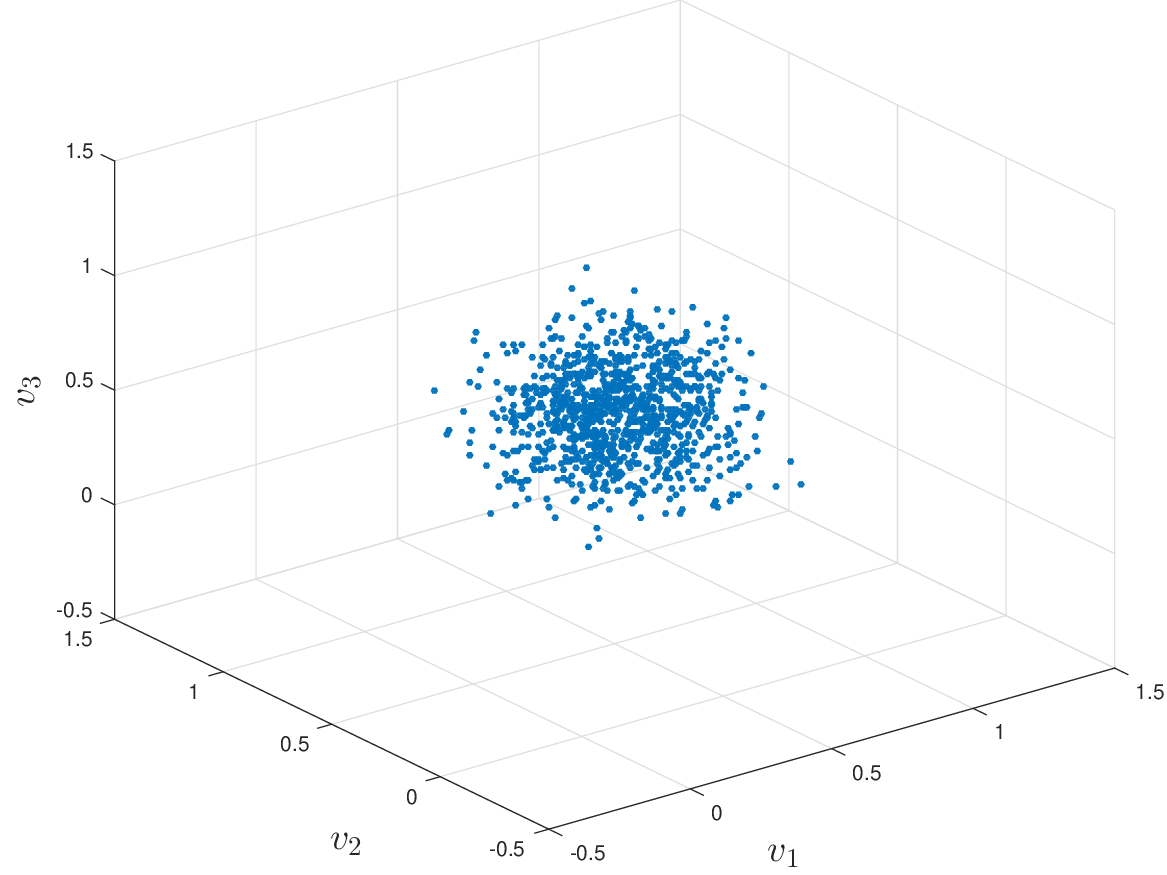}
        \caption{Uncontrolled velocity, $t=1$}
    \end{subfigure}\quad
    \begin{subfigure}{.31\textwidth}
        \centering
         \includegraphics[width=\linewidth,height=4.4cm]{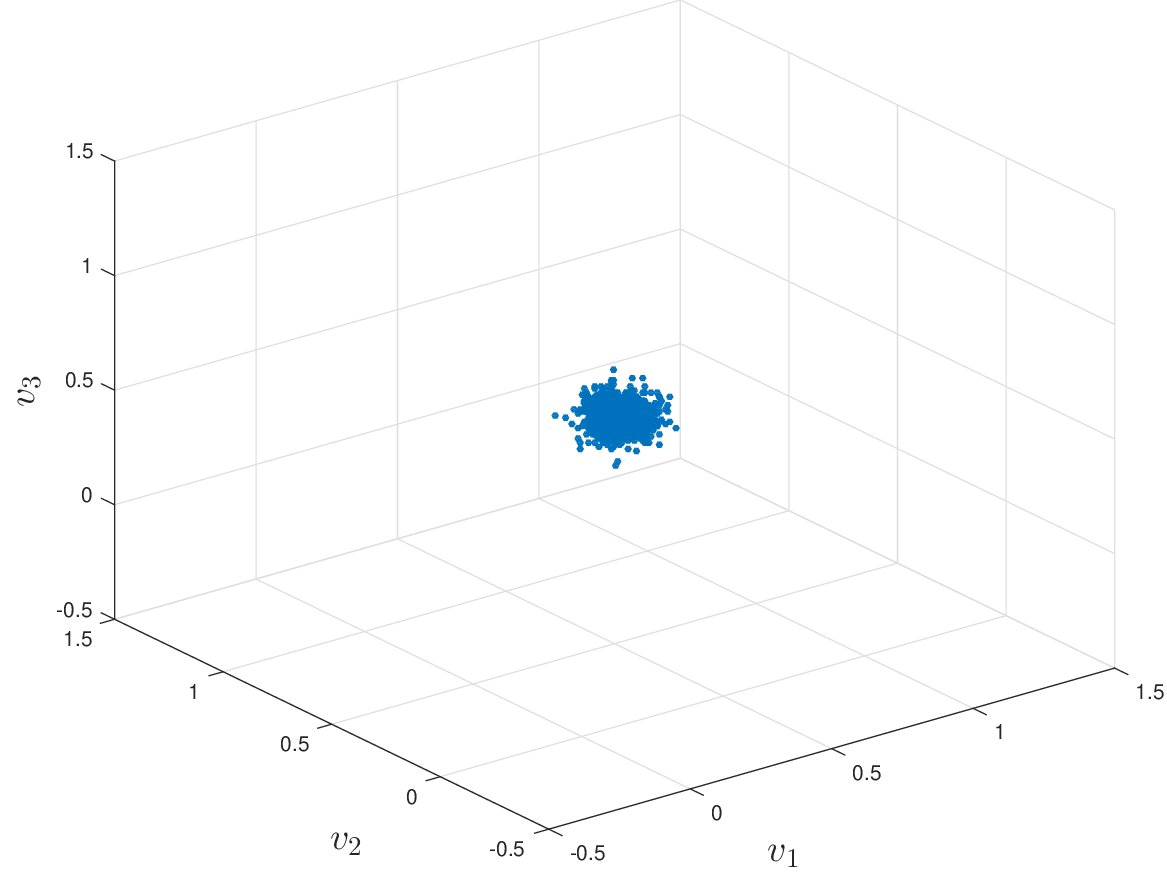}
        \caption{Controlled velocity, $t=1$}
        \end{subfigure}

         \caption{
Uncontrolled and controlled six-dimensional C-S models with $\beta=1$
  at different times.
 }
    \label{fig:cs_v_6d_beta1}
 \end{figure}

Finally, we apply the FIPDE method to solve the nonconvex  MFC problem 
  \eqref{eq:cs}-\eqref{eq:cs_loss}
 with $\beta=1$ and $\gamma_2=0$.
As in the two-dimensional setting studied in Section \ref{sec:CS_2d},
the uncontrolled velocity process does not form a consensus,
and the optimal feedback control in general depends 
 nonlinearly
on $(t,x,v)$. 
Similar to the above nonsmooth problem, 
we use the LQ feedback control  
as 
 the initial guess   $\phi^0$
in Algorithm \ref{alg:NAG2},
choose  the stepsize $\tau=1/5$,
and 
employ the residual approximation approach 
to solve  \eqref{eq:pde_cs} at all NAG iterations,
with the candidate solution $\bar{u}$ 
being the approximate decoupling field of the  LQ problem.
The components of the residual $\tilde{u}^m$
and the controls $\phi^m,\psi^m$
are approximated by 
 $7$-layer neural  networks with hidden width   60
 (see {Appendix \ref{appendix:cs}} for more details).
Figure \ref{fig:cs_v_6d_beta1} compares the uncontrolled velocity process 
and the  controlled velocity process obtained by the FIPDE method
(with $6$ NAG iterations), 
which clearly demonstrates the effectiveness of
   the FIPDE method on inducing consensus.
{
Compared with  the NNPG method (implemented with the same neural networks),
 the  FIPDE method 
results in  
 a similar value function but 
 a better consensus at the terminal time; 
 the value functions 
of the FIPDE and NNPG methods
are 0.111  and 0.114, respectively,
while 
 the final variances $\sE[|v_T-\sE[v_T]|^2]$
obtained by the FIPDE and NNPG methods
are 0.007  and 0.01, respectively.
Moreover, similar to  the above examples 
 in Figures \ref{fig:cs_6d_beta0_fb} and \ref{fig:value_conv_6d_l1},
the FIPDE method captures the nonlinear time dependence of the optimal control,
while the NNPG method leads to an approximate control that is  constant in $t$ (plots omitted).
}

\begin{appendices}

\section{Supplementary materials for Section \ref{SEC:NUM}}
\label{SEC:supplementary}

\subsection{Supplementary materials for Section \ref{SEC:PF}}
\label{appendix:portfolio}

\paragraph{Convergence of  the  NAG iteration.}

\begin{figure}[!ht]
    \centering
       \begin{subfigure}{.3\textwidth}
        \centering
        \includegraphics[width=\linewidth,height=4.4cm]{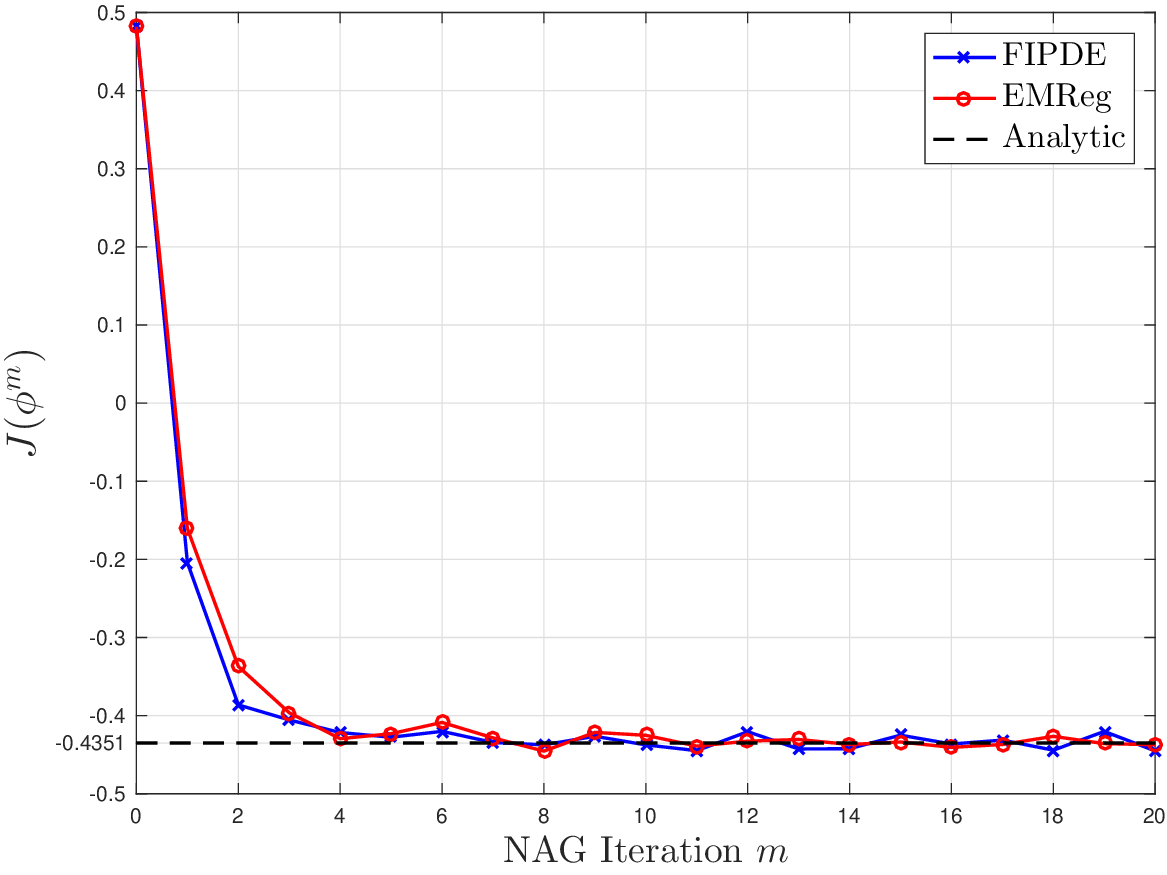}
        \caption{LQ case   with   $k_2=0$}
    \end{subfigure}     \hspace{1.5cm} 
    \begin{subfigure}{.3\textwidth}
        \centering
         \includegraphics[width=\linewidth,height=4.4cm]{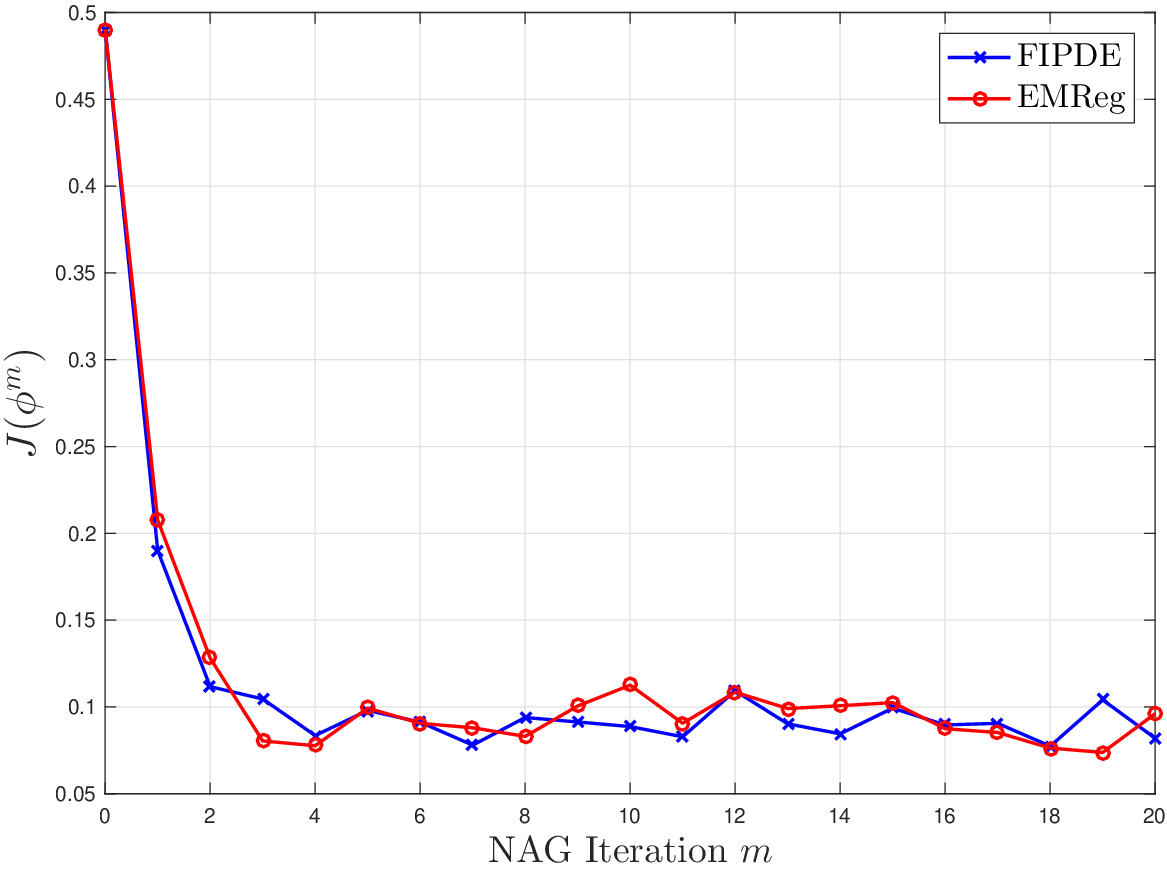}
        \caption{Nonsmooth case with  $k_2=1$}
    \end{subfigure}

    \caption{
  Convergence of
   the FIPDE and EMReg methods for approximating the value functions of 
optimal liquidation problems \eqref{eq:price_impact}
in terms of  NAG iteration.
 }
    \label{fig:price_value_conv}
 \end{figure}

Figure \ref{fig:price_value_conv} presents  the expected costs  
of
 the approximate feedback controls obtained by
all NAG iterations of 
  the FIPDE and EMReg methods
for the  optimal liquidation problem \eqref{eq:price_impact}
with $Q_0\in \cU(1,2)$, $k_1=1$ and $k_2\in \{0,1\}$.
To estimate the optimal expected cost for 
the LQ setting (with $k_2=0$),
we
implement  the exact feedback control \eqref{eq:optimal_fb_price}
with time stepsize $\Delta t =1/50$,
and replace the expectation in \eqref{eq:price_impact}
by the empirical average over $10^6$ sample trajectories of the state processes.
It shows that 
for both the LQ and nonsmooth MFC problems,  
 the FIPDE and EMReg methods  give convergent approximations 
to the optimal cost functions 
as the number of NAG iterations tends to infinity.
Moreover, 
a linear regression of the data  
indicates  that the FIPDE method 
approximates the value function
with an absolute error  
of the magnitude $\cO(m^{-2.7})$ for the LQ case.

\paragraph{Computational time.}

 The FIPDE method is  implemented by using  \textsc{Matlab} R2016b
 on a laptop with 
 2.2GHz 4-core Intel Core i7 processor
 and   16 GB   memory. 
 The computation takes around 1 second per NAG iteration. 
 The EMReg method  is  implemented by using  \textsc{Matlab} R2020b
 on a PC with 
 2.1GHz 6-core Intel Core i5 processor
 and   16 GB   memory. 
 The computation  takes around 5 minutes per NAG iteration based on our implementation.

\subsection{Supplementary materials for Section \ref{SEC:CS}}
\label{appendix:cs}

\paragraph{Effectiveness of the FIPDE method for  two-dimensional models.}

\begin{figure}[!ht]
 \begin{subfigure}{\textwidth}
        \centering
             \includegraphics[  width=0.31\columnwidth,height=4.5cm]{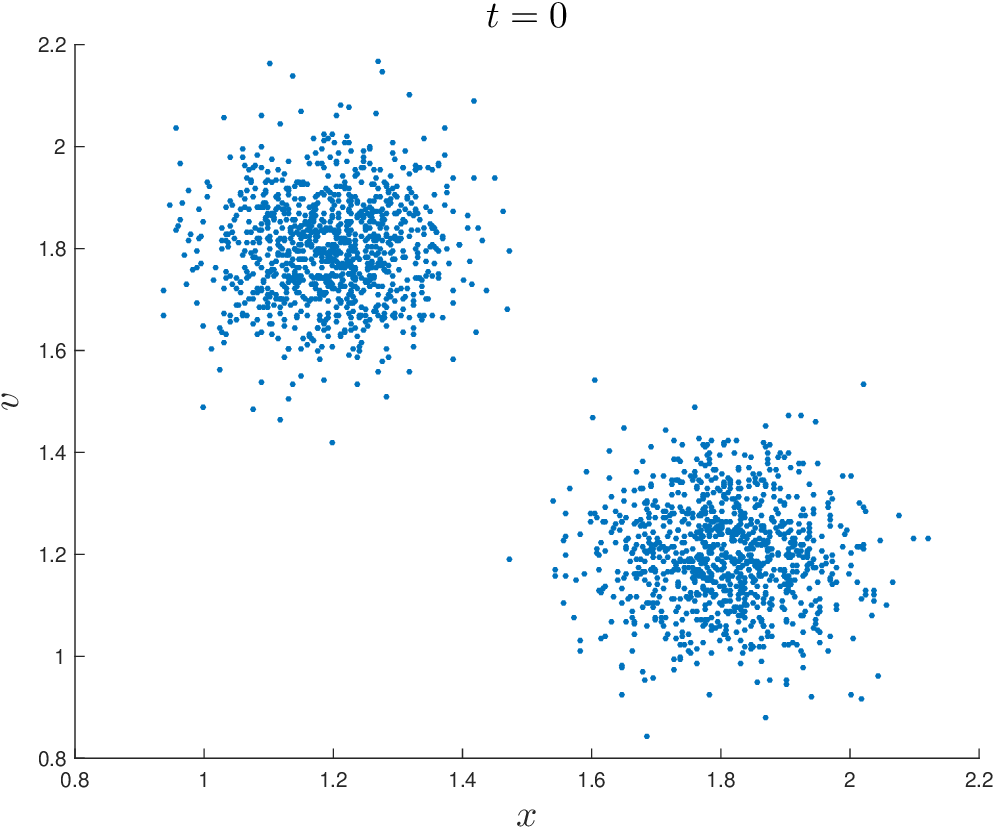}
         \quad
       \includegraphics[   width=0.31\columnwidth,height=4.5cm]{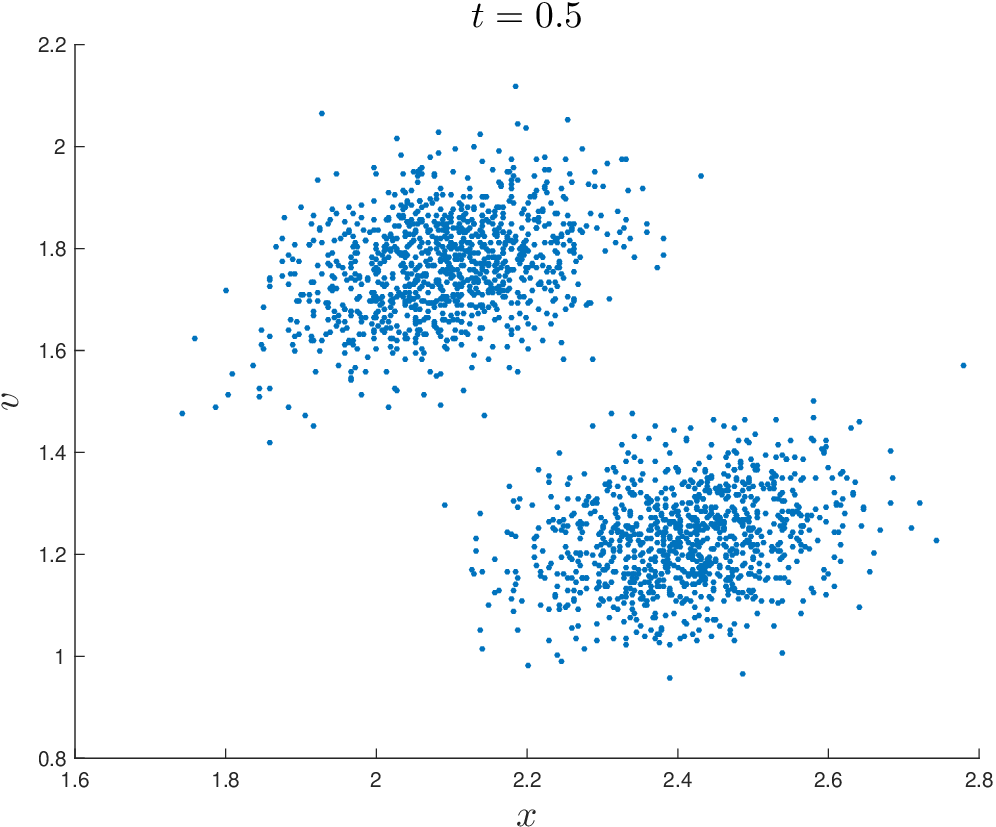}
                \quad
       \includegraphics[   width=0.31\columnwidth,height=4.5cm]{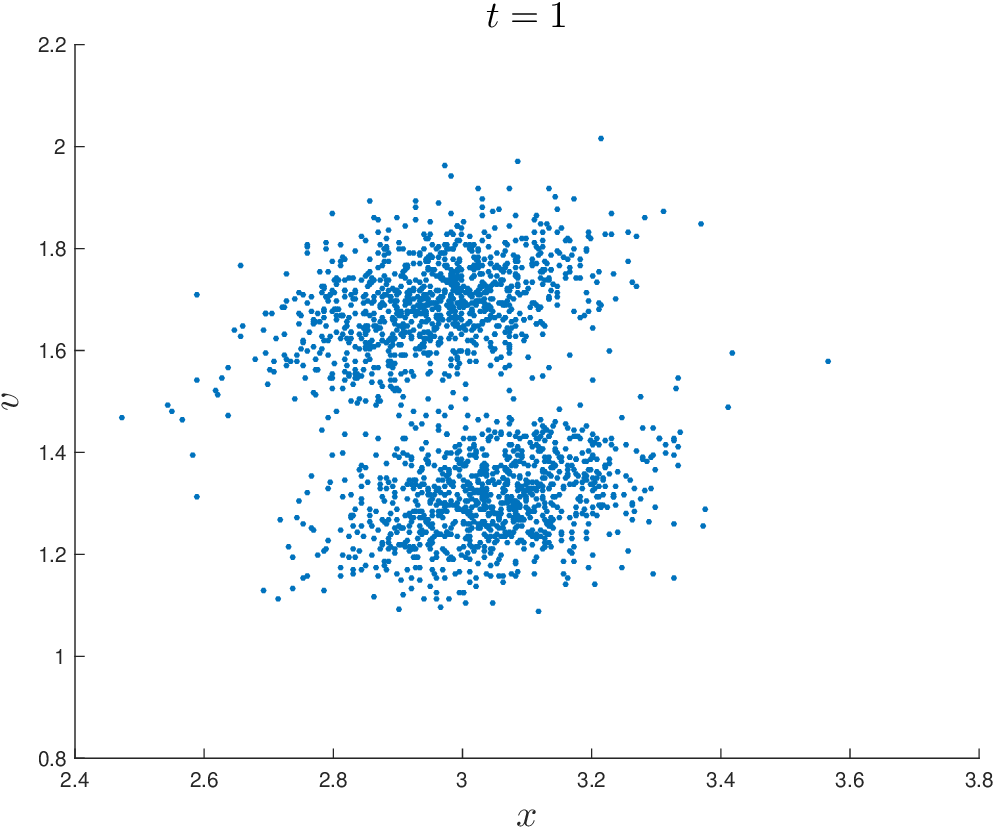}
      
        \caption{Uncontrolled position and velocity}
    \end{subfigure}
     \vspace{2mm}

 \begin{subfigure}{\textwidth}
        \centering
         \includegraphics[  width=0.31\columnwidth,height=4.5cm]{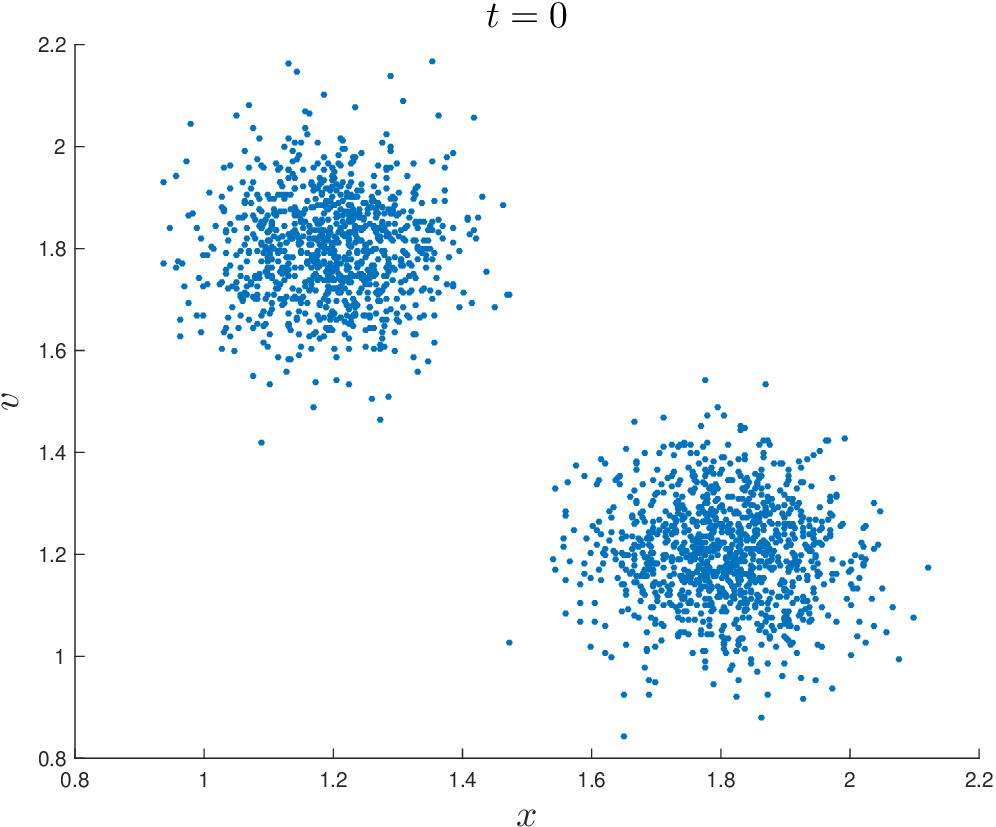}
         \quad
       \includegraphics[   width=0.31\columnwidth,height=4.5cm]{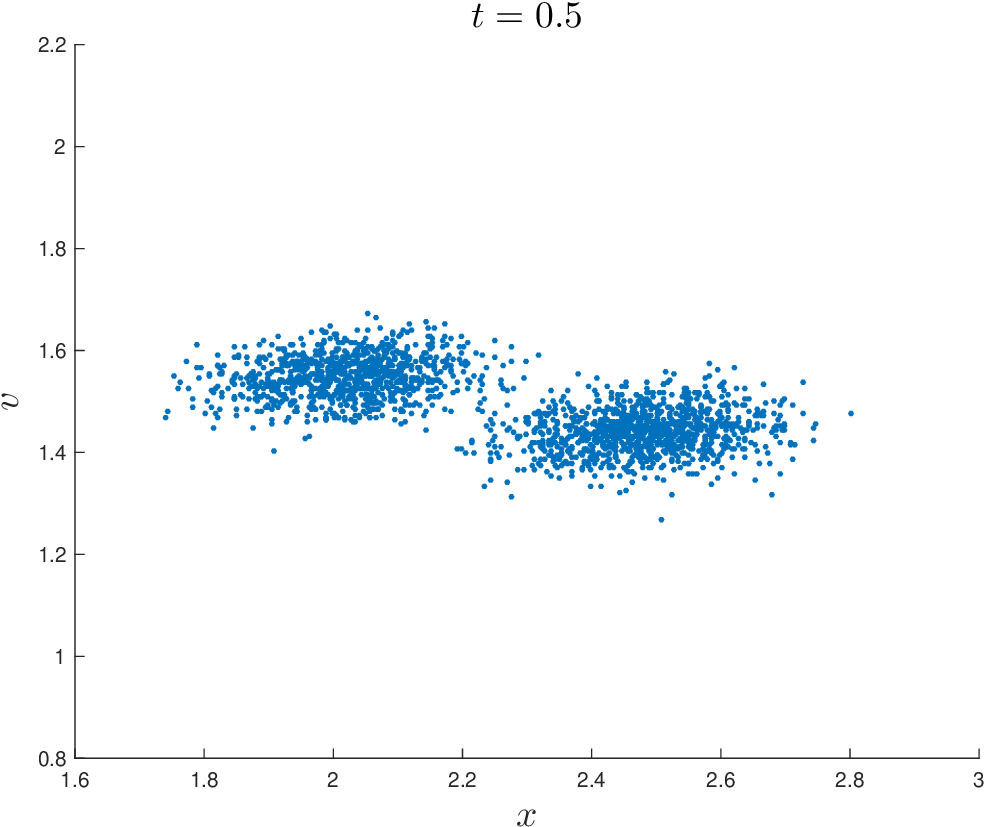}
       \quad
       \includegraphics[   width=0.31\columnwidth,height=4.5cm]{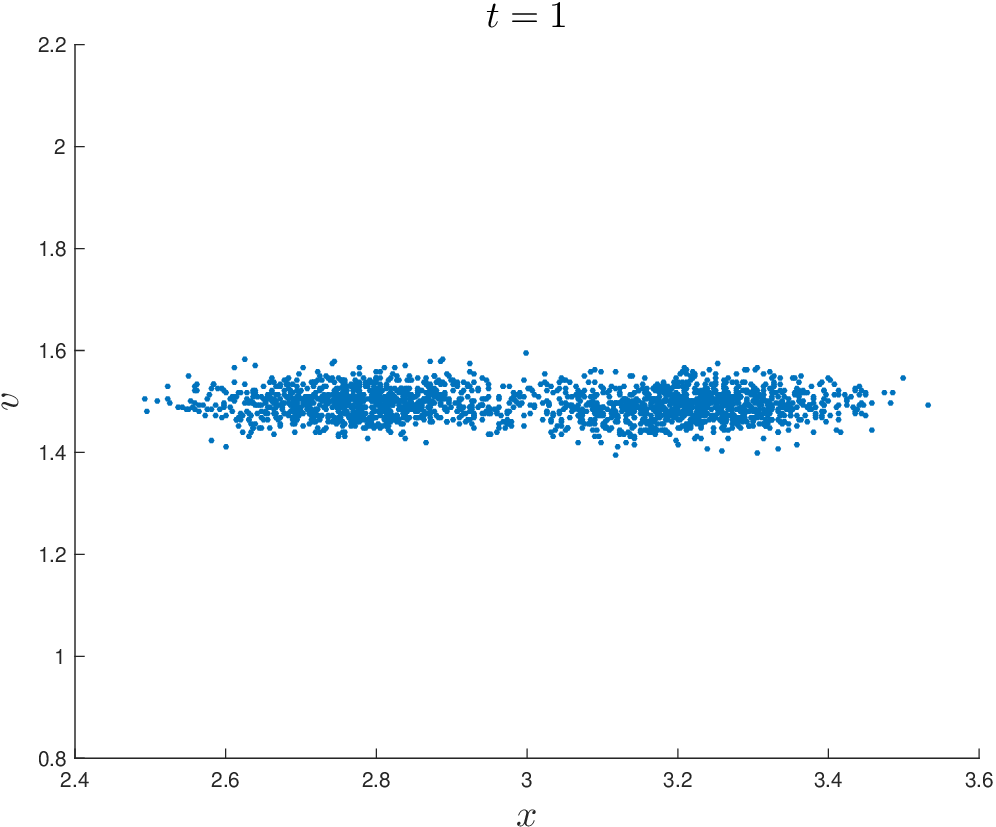}

        \caption{Controlled position and velocity obtained by the FIPDE method}
    \end{subfigure}
   
         \caption{
Uncontrolled and controlled two-dimensional C-S models with $\beta=10$
  at different times.
 }
    \label{fig:cs_state_2d}
 \end{figure}

 Figure \ref{fig:cs_state_2d} compares  the 
 uncontrolled  two-dimensional C-S model with $\beta=10$
 and the controlled model obtained by the FIPDE method,
 where  the  scatter plots are generated 
   based on $10^4$ simulated trajectories.
  One can clearly observe  that
   the uncontrolled velocity process does not admit a time-asymptotic flocking behaviour,
   and 
    the feedback control from the FIPDE method effectively induces the consensus of the velocity process.

\paragraph{Effectiveness of the FIPDE method for  six-dimensional models with $\beta=0$.}

 Figure \ref{fig:cs_state_2d} compares  the 
 uncontrolled  six-dimensional C-S model with $\beta=0$
 and the controlled model obtained by the FIPDE method
 with quadratic costs (i.e., $\gamma_2=0$ in \eqref{eq:cs_loss}),
 where  the  scatter plots are generated 
   based on $5\t 10^3$ simulated trajectories.
  One can clearly observe  that
    the feedback control from the FIPDE method  accelerates  the consensus of the velocity process.

\begin{figure}[!ht]
     \begin{subfigure}{.31\textwidth}
        \centering
        \includegraphics[width=\linewidth,height=4.4cm]{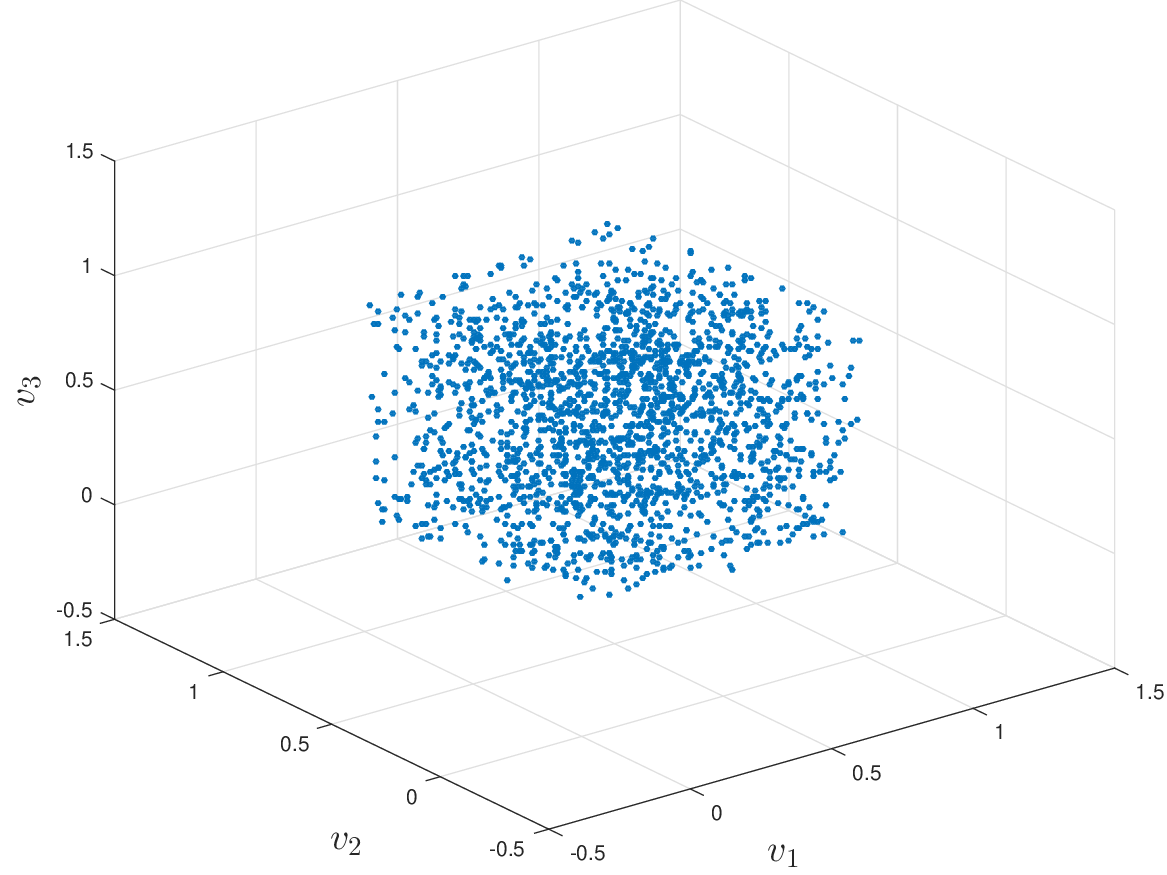}
        \caption{Uncontrolled velocity, $t=0$}
    \end{subfigure}\quad
    \begin{subfigure}{.31\textwidth}
        \centering
         \includegraphics[width=\linewidth,height=4.4cm]{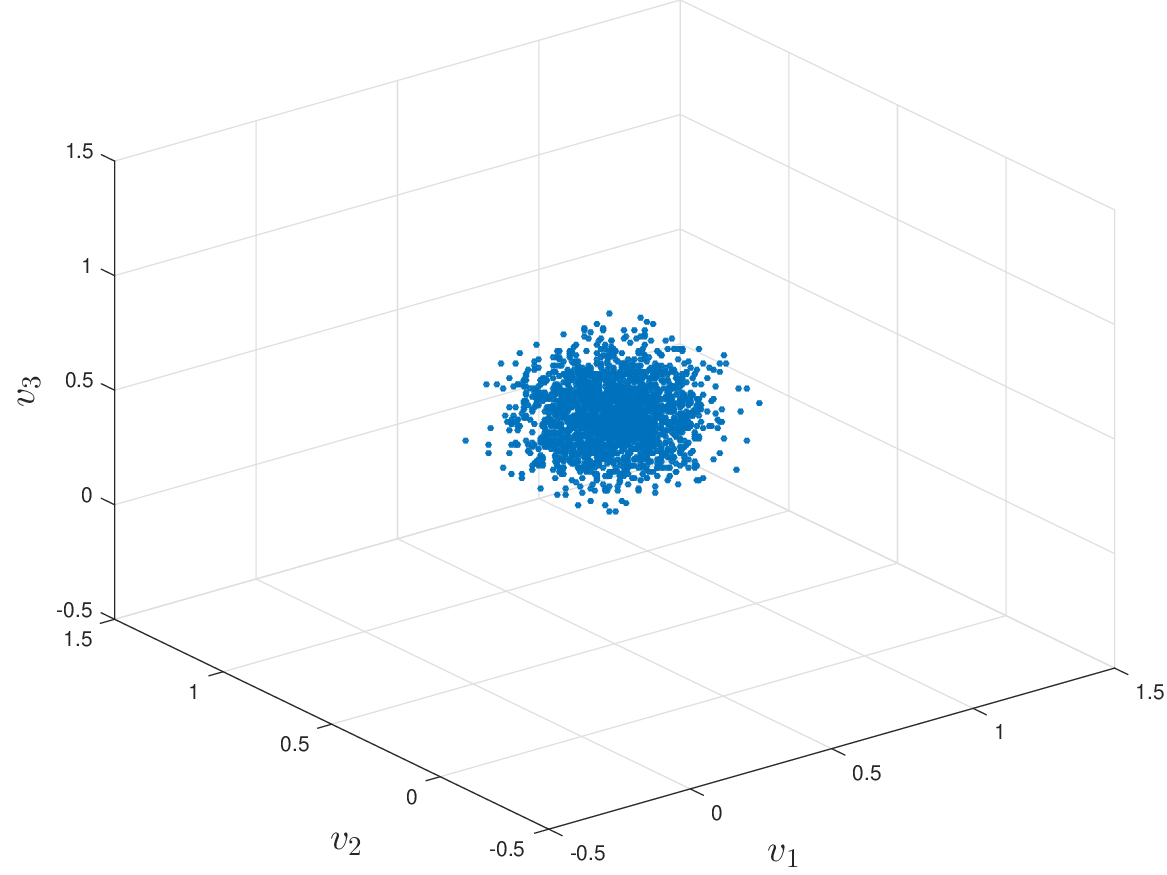}
        \caption{Uncontrolled velocity, $t=1$}
    \end{subfigure}\quad
    \begin{subfigure}{.31\textwidth}
        \centering
         \includegraphics[width=\linewidth,height=4.4cm]{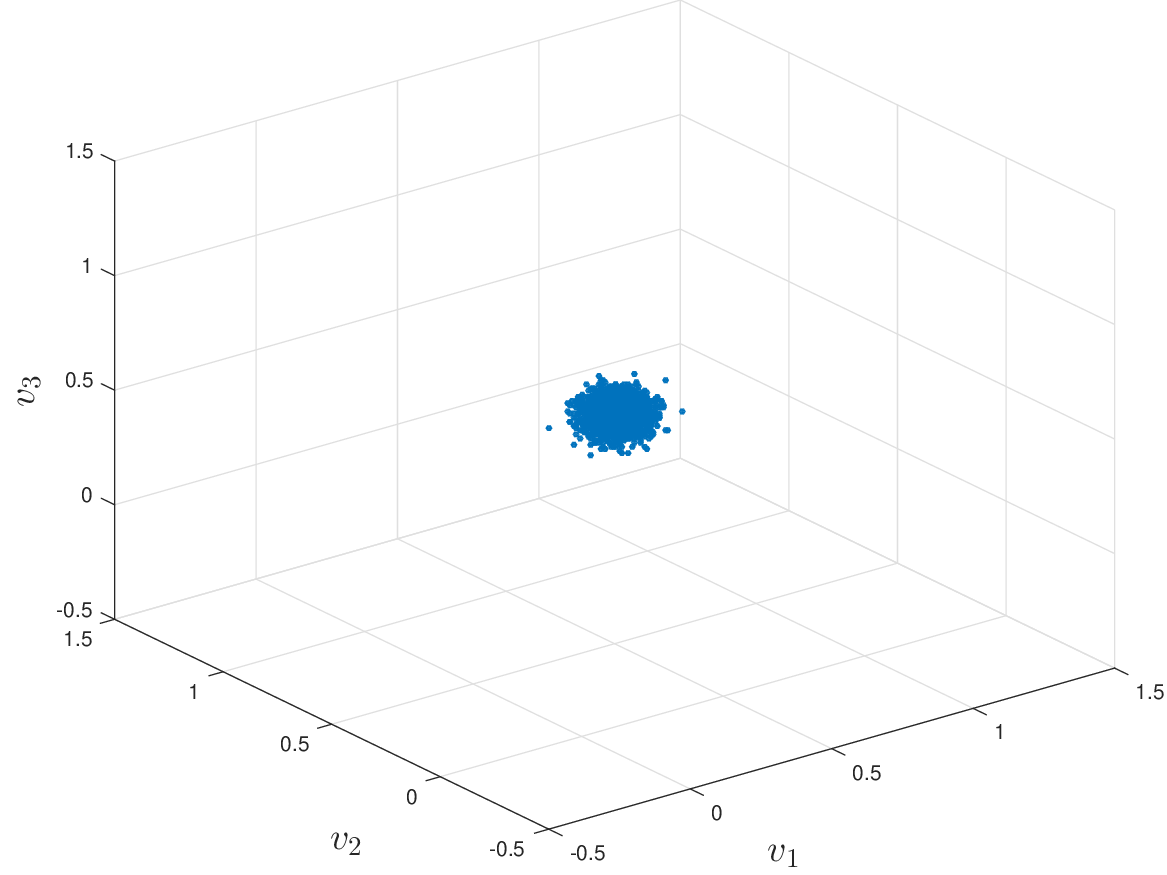}
        \caption{Controlled velocity, $t=1$}
        \end{subfigure}

         \caption{
Uncontrolled and controlled six-dimensional C-S models with $\beta=0$
  at different times.
 }
    \label{fig:cs_v_6d_lq}
 \end{figure}

\paragraph{Implementation of the NNPG method for two-dimensional models.}
 We implement the neural network-based policy gradient (NNPG) method using PyTorch
for the optimal control of two-dimensional C-S models with $\beta=10$ 
as follows.
The set of trial functions for feedback controls consists of $7$ layer fully-connected networks with the sigmoid activation function
and  the dimensions of input, hidden and output  layers being equal to $3$, $70$ and $1$, respectively.
 The  optimal weights $\theta$ of the network  are obtained by  running the   Adam   algorithm
 \cite{kingma2014adam} with $5\t 10^4$ iterations and  a decaying learning rate scheduler. 
At each  iteration, we 
 generate   $10^2$  trajectories of  \eqref{eq:cs} by using  \eqref{eq:forward_samples} with time stepsize $\Delta t=0.02$,
 approximate the functional $J$ in \eqref{eq:cs_loss} by the empirical average of simulated trajectories,
 and perform  gradient descent  based on the empirical loss. 
The learning rate is  initialized at $0.0005$ and 
will be reduced by 
 a factor of 0.3 
 if the loss does not decrease after 10 iterations. 
 The minimum learning rate is set to be $10^{-6}$.

\paragraph{Implementation of the NNPG method for six-dimensional models.}
We state  the  PyTorch implementation  of the NNPG method for six-dimensional C-S models
with different choices of $\beta,\gamma_2$ used  in Section \ref{SEC:CS}.
As  the chosen neural network architecture for each  example has been specified in Section \ref{SEC:CS},  it remains to discuss  the  implementation of the Adam algorithm. 
For $\beta=0$ and $\gamma_2\in \{0,0.1,0.3\}$,
at each SGD iteration,
we generate   $10^3$  trajectories of  \eqref{eq:cs} by using  \eqref{eq:forward_samples} with time stepsize $\Delta t=0.02$, 
and perform  gradient descent  based on an empirical approximation of  $J$ in \eqref{eq:cs_loss}.
We choose  
the same  ``reduce  on plateau" learning rate scheduler   as the above two-dimensional case
and perform a sufficiently  large number of iterations until the empirical loss stabilizes.
For $\beta=1$ and $\gamma_2=0$,
we reduce the  number of simulated trajectories per iteration to $10^2$ 
and keep the remaining configurations.
This helps to accommodate the increasing computational cost
of simulating  \eqref{eq:cs} with $\beta>0$.

\paragraph{Implementation of the FIPDE method for six-dimensional models.}

We state  the PyTorch implementation  of the FIPDE method for six-dimensional C-S models
with different choices of $\beta,\gamma_2$ used in Section \ref{SEC:CS}.
Before stating the configuration details, we remark that 
the algorithm's hyperparameters 
have not been optimally tuned
and hence 
the following choices   may not be the optimal combination. 
The neural network architectures for all  examples have been specified in Section \ref{SEC:CS}.

For all NAG iterations,
we sample $N$ trajectories  of \eqref{eq:cs}
by using  \eqref{eq:forward_samples} with  $\Delta t=0.01$,
where $N=4\t 10^4$ for $\beta=\gamma_2=0$,
and $N=5\t 10^3$ for the remaining cases.
Then, for all 
choices of $\beta,\gamma_2$,
we apply the Adam algorithm to 
 solve the corresponding risk minimization problem \eqref{eq:cE}
 with the initial  learning rate 0.005.
The learning rate  will be decreased  by employing the same  ``reduce  on plateau"
scheduler as the  above NNPG method 
for the examples with $\beta=0, \gamma_2\in\{0,0.1,0.3\}$,
 while 
for the example with $\beta=1, \gamma_2=0$
it will be decreased by a factor of 0.8
for every 5000 SGD iterations.
 
At each SGD iteration, 
a mini-batch of points with size 20 
is drawn by following uniform distributions from   
the interior and boundaries of the domain (i.e., $N_{\textrm{in}}=N_{\textrm{ter}}=N_{\textrm{bdy}}=20$ in \eqref{eq:E_empirical}).
 Based on these samples, we compute the empirical loss  \eqref{eq:E_empirical} with certain $\eta_1, \eta_2>0$
depending on the model parameters.   
In particular,
for $\beta=\gamma_2=0$,
 we choose 
 $\eta_1=5$, $\eta_2=10$ for the first 6 NAG iterations 
 and $\eta_1=2$, $\eta_2=5$ for the remaining NAG iterations
 to balance the interior and boundary losses,
 while for  $\beta=0$, $\gamma_2=\{0,1,0.3\}$ 
 and $\beta=1$, $\gamma_2=0$,
 we choose $\eta_1=1,\eta_2=1/400$ and  $\eta_1=1,\eta_2=1/800$, respectively, across all NAG iterations. 
 The total number of  Adam iterations 
is chosen as 9000 for $\beta=\gamma_2=0$, 
100 for $\beta=0,\gamma_2\in \{0.1,0.3\}$,
and 20000 for $\beta=1, \gamma_2=0$.
  
To solve the supervised learning problem \eqref{eq:supervised}, 
we carry out the Adam algorithm with initial learning rate 0.001,
which will be reduced by employing the same scheduler as that for solving \eqref{eq:cE}.
We shall     randomly draw a mini-batch of points with size 20 at each  iteration,
and run the algorithm with sufficiently many iterations until the loss stabilizes.

\paragraph{Computational time.}

For the two-dimensional C-S models, 
the FIPDE method is  implemented in  \textsc{Matlab} R2016b
 on a laptop with 
 2.2GHz 4-core Intel Core i7 processor
 and   16 GB   memory. 
 The computation with 15 NAG iterations takes around 30 seconds
 for $\beta=0$ and around $10$ minutes for $\beta=10$,
 due to an increased computational cost in evaluating the interaction kernel $\kappa$ 
 for general $\beta>0$. 
The NNPG method is   implemented in Python  3.8.5
 on a PC with 
 2.1GHz 6-core Intel Core i5 processor
 and   16 GB   memory. 
 The computation with  $5\t 10^4$ Adam iterations  takes around 2 hours.
 
 For the six-dimensional C-S models, 
both the FIPDE and NNPG methods are implemented by using  Python  3.8.5
 on a PC with 
 2.1GHz 6-core Intel Core i5 processor
 and   16 GB   memory. 
 The NNPG method typically takes around 20-40 minutes, 
 where the precise computation time depends on the complexity of the problem
 and the number of SGD iterations used in the simulation. 
 The FIPDE method 
takes around 20 minutes per NAG iteration for $\beta=0$,
 and around 8 hours per NAG iteration for  $\beta=1$.
 Note that these  times can be shortened 
 if we fine tune the hyperparameters  (e.g., the number of SGD iterations)
 and perform the computation on GPUs. 
  
\end{appendices}



\subsection*{Acknowledgements}

Wolfgang Stockinger is supported by a special Upper Austrian  Government grant.

\bibliographystyle{siam}
\bibliography{nag_mfc.bib}

\end{document}